
\ifx\shlhetal\undefinedcontrolsequence\let\shlhetal\relax\fi
\def\fmtname{AmS-TeX}

\def\fmtversion{2.2}
\catcode`\@=11
\ifx\amstexloaded@\relax\catcode`\@=\active
  \endinput\else\let\amstexloaded@\relax\fi
\newlinechar=`\^^J
\def\W@{\immediate\write\sixt@@n}
\def\CR@{\W@{^^J\fmtname - Version \fmtversion^^J}}
\CR@ \everyjob{\CR@}
\message{Loading definitions for}
\message{misc utility macros,}
\toksdef\toks@@=2
\long\def\rightappend@#1\to#2{\toks@{\\{#1}}\toks@@
 =\expandafter{#2}\xdef#2{\the\toks@@\the\toks@}\toks@{}\toks@@{}}
\def\alloclist@{}
\newif\ifalloc@
\def\showallocations{{\def\\{\immediate\write\m@ne}\alloclist@}\alloc@true}
\def\alloc@#1#2#3#4#5{\global\advance\count1#1by\@ne
 \ch@ck#1#4#2\allocationnumber=\count1#1
 \global#3#5=\allocationnumber
 \edef\next@{\string#5=\string#2\the\allocationnumber}%
 \expandafter\rightappend@\next@\to\alloclist@}
\newcount\count@@
\newcount\count@@@
\def\FN@{\futurelet\next}
\def\DN@{\def\next@}
\def\DNii@{\def\nextii@}
\def\RIfM@{\relax\ifmmode}
\def\RIfMIfI@{\relax\ifmmode\ifinner}
\def\setboxz@h{\setbox\z@\hbox}
\def\wdz@{\wd\z@}
\def\boxz@{\box\z@}
\def\setbox@ne{\setbox\@ne}
\def\wd@ne{\wd\@ne}
\def\iterate{\body\expandafter\iterate\else\fi}
\def\err@#1{\errmessage{AmS-TeX error: #1}}
\newhelp\defaulthelp@{Sorry, I already gave what help I could...^^J
Maybe you should try asking a human?^^J
An error might have occurred before I noticed any problems.^^J
``If all else fails, read the instructions.''}
\def\Err@{\errhelp\defaulthelp@\err@}
\def\eat@#1{}
\def\in@#1#2{\def\in@@##1#1##2##3\in@@{\ifx\in@##2\in@false\else\in@true\fi}%
 \in@@#2#1\in@\in@@}
\newif\ifin@
\def\space@.{\futurelet\space@\relax}
\space@. %
\newhelp\athelp@
{Only certain combinations beginning with @ make sense to me.^^J
Perhaps you wanted \string\@\space for a printed @?^^J
I've ignored the character or group after @.}
{\catcode`\~=\active 
 \lccode`\~=`\@ \lowercase{\gdef~{\FN@\at@}}}
\def\at@{\let\next@\at@@
 \ifcat\noexpand\next a\else\ifcat\noexpand\next0\else
 \ifcat\noexpand\next\relax\else
   \let\next\at@@@\fi\fi\fi
 \next@}
\def\at@@#1{\expandafter
 \ifx\csname\space @\string#1\endcsname\relax
  \expandafter\at@@@ \else
  \csname\space @\string#1\expandafter\endcsname\fi}
\def\at@@@#1{\errhelp\athelp@ \err@{\Invalid@@ @}}
\def\atdef@#1{\expandafter\def\csname\space @\string#1\endcsname}
\newhelp\defahelp@{If you typed \string\define\space cs instead of
\string\define\string\cs\space^^J
I've substituted an inaccessible control sequence so that your^^J
definition will be completed without mixing me up too badly.^^J
If you typed \string\define{\string\cs} the inaccessible control sequence^^J
was defined to be \string\cs, and the rest of your^^J
definition appears as input.}
\newhelp\defbhelp@{I've ignored your definition, because it might^^J
conflict with other uses that are important to me.}
\def\define{\FN@\define@}
\def\define@{\ifcat\noexpand\next\relax
 \expandafter\define@@\else\errhelp\defahelp@                               
 \err@{\string\define\space must be followed by a control
 sequence}\expandafter\def\expandafter\nextii@\fi}                          
\def\undefined@@@@@@@@@@{}
\def\preloaded@@@@@@@@@@{}
\def\next@@@@@@@@@@{}
\def\define@@#1{\ifx#1\relax\errhelp\defbhelp@                              
 \err@{\string#1\space is already defined}\DN@{\DNii@}\else
 \expandafter\ifx\csname\expandafter\eat@\string                            
 #1@@@@@@@@@@\endcsname\undefined@@@@@@@@@@\errhelp\defbhelp@
 \err@{\string#1\space can't be defined}\DN@{\DNii@}\else
 \expandafter\ifx\csname\expandafter\eat@\string#1\endcsname\relax          
 \global\let#1\undefined\DN@{\def#1}\else\errhelp\defbhelp@
 \err@{\string#1\space is already defined}\DN@{\DNii@}\fi
 \fi\fi\next@}

\def\predefine#1#2{\let#1#2}
\def\undefine#1{\let#1\undefined}
\message{page layout,}
\newdimen\captionwidth@
\captionwidth@\hsize
\advance\captionwidth@-1.5in
\def\pagewidth#1{\hsize#1\relax
 \captionwidth@\hsize\advance\captionwidth@-1.5in}
\def\pageheight#1{\vsize#1\relax}
\def\hcorrection#1{\advance\hoffset#1\relax}
\def\vcorrection#1{\advance\voffset#1\relax}
\message{accents/punctuation,}

\let\graveaccent\`
\let\acuteaccent\'
\let\tildeaccent\~
\let\hataccent\^
\let\underscore\_
\let\B\=
\let\D\.
\let\ic@\/
\def\/{\unskip\ic@}
\def\textfonti{\the\textfont\@ne}
\def\t#1#2{{\edef\next@{\the\font}\textfonti\accent"7F \next@#1#2}}
\def~{\unskip\nobreak\ \ignorespaces}
\def\.{.\spacefactor\@m}
\atdef@;{\leavevmode\null;}
\atdef@:{\leavevmode\null:}
\atdef@?{\leavevmode\null?}
\edef\@{\string @}
\def\relaxnext@{\let\next\relax}
\atdef@-{\relaxnext@\leavevmode
 \DN@{\ifx\next-\DN@-{\FN@\nextii@}\else
  \DN@{\leavevmode\hbox{-}}\fi\next@}%
 \DNii@{\ifx\next-\DN@-{\leavevmode\hbox{---}}\else
  \DN@{\leavevmode\hbox{--}}\fi\next@}%
 \FN@\next@}
\def\srdr@{\kern.16667em}
\def\drsr@{\kern.02778em}
\def\sldl@{\drsr@}
\def\dlsl@{\srdr@}
\atdef@"{\unskip\relaxnext@
 \DN@{\ifx\next\space@\DN@. {\FN@\nextii@}\else
  \DN@.{\FN@\nextii@}\fi\next@.}%
 \DNii@{\ifx\next`\DN@`{\FN@\nextiii@}\else
  \ifx\next\lq\DN@\lq{\FN@\nextiii@}\else
  \DN@####1{\FN@\nextiv@}\fi\fi\next@}%
 \def\nextiii@{\ifx\next`\DN@`{\sldl@``}\else\ifx\next\lq
  \DN@\lq{\sldl@``}\else\DN@{\dlsl@`}\fi\fi\next@}%
 \def\nextiv@{\ifx\next'\DN@'{\srdr@''}\else
  \ifx\next\rq\DN@\rq{\srdr@''}\else\DN@{\drsr@'}\fi\fi\next@}%
 \FN@\next@}

\def\textfontii{\the\textfont\tw@}
\def\lbrace@{\delimiter"4266308 }
\def\rbrace@{\delimiter"5267309 }
\def\{{\RIfM@\lbrace@\else{\textfontii f}\spacefactor\@m\fi}
\def\}{\RIfM@\rbrace@\else
 \let\@sf\empty\ifhmode\edef\@sf{\spacefactor\the\spacefactor}\fi
 {\textfontii g}\@sf\relax\fi}
\let\lbrace\{
\let\rbrace\}
\def\AmSTeX{{\textfontii A\kern-.1667em%
  \lower.5ex\hbox{M}\kern-.125emS}-\TeX\spacefactor1000 }
\message{line and page breaks,}
\def\vmodeerr@#1{\Err@{\string#1\space not allowed between paragraphs}}
\def\mathmodeerr@#1{\Err@{\string#1\space not allowed in math mode}}
\def\linebreak{\RIfM@\mathmodeerr@\linebreak\else
 \ifhmode\unskip\unkern\break\else\vmodeerr@\linebreak\fi\fi}

\newskip\saveskip@
\def\allowlinebreak{\RIfM@\mathmodeerr@\allowlinebreak\else
 \ifhmode\saveskip@\lastskip\unskip
 \allowbreak\ifdim\saveskip@>\z@\hskip\saveskip@\fi
 \else\vmodeerr@\allowlinebreak\fi\fi}
\def\nolinebreak{\RIfM@\mathmodeerr@\nolinebreak\else
 \ifhmode\saveskip@\lastskip\unskip
 \nobreak\ifdim\saveskip@>\z@\hskip\saveskip@\fi
 \else\vmodeerr@\nolinebreak\fi\fi}
\def\newline{\relaxnext@
 \DN@{\RIfM@\expandafter\mathmodeerr@\expandafter\newline\else
  \ifhmode\ifx\next\par\else
  \expandafter\unskip\expandafter\null\expandafter\hfill\expandafter\break\fi
  \else
  \expandafter\vmodeerr@\expandafter\newline\fi\fi}%
 \FN@\next@}
\def\dmatherr@#1{\Err@{\string#1\space not allowed in display math mode}}
\def\nondmatherr@#1{\Err@{\string#1\space not allowed in non-display math
 mode}}
\def\onlydmatherr@#1{\Err@{\string#1\space allowed only in display math mode}}
\def\nonmatherr@#1{\Err@{\string#1\space allowed only in math mode}}
\def\mathbreak{\RIfMIfI@\break\else
 \dmatherr@\mathbreak\fi\else\nonmatherr@\mathbreak\fi}
\def\nomathbreak{\RIfMIfI@\nobreak\else
 \dmatherr@\nomathbreak\fi\else\nonmatherr@\nomathbreak\fi}
\def\allowmathbreak{\RIfMIfI@\allowbreak\else
 \dmatherr@\allowmathbreak\fi\else\nonmatherr@\allowmathbreak\fi}
\def\pagebreak{\RIfM@
 \ifinner\nondmatherr@\pagebreak\else\postdisplaypenalty-\@M\fi
 \else\ifvmode\removelastskip\break\else\vadjust{\break}\fi\fi}
\def\nopagebreak{\RIfM@
 \ifinner\nondmatherr@\nopagebreak\else\postdisplaypenalty\@M\fi
 \else\ifvmode\nobreak\else\vadjust{\nobreak}\fi\fi}
\def\nonvmodeerr@#1{\Err@{\string#1\space not allowed within a paragraph
 or in math}}
\def\vnonvmode@#1#2{\relaxnext@\DNii@{\ifx\next\par\DN@{#1}\else
 \DN@{#2}\fi\next@}%
 \ifvmode\DN@{#1}\else
 \DN@{\FN@\nextii@}\fi\next@}
\def\newpage{\vnonvmode@{\vfill\break}{\nonvmodeerr@\newpage}}
\def\smallpagebreak{\vnonvmode@\smallbreak{\nonvmodeerr@\smallpagebreak}}
\def\medpagebreak{\vnonvmode@\medbreak{\nonvmodeerr@\medpagebreak}}
\def\bigpagebreak{\vnonvmode@\bigbreak{\nonvmodeerr@\bigpagebreak}}
\def\NoBlackBoxes{\global\overfullrule\z@}
\def\BlackBoxes{\global\overfullrule5\p@}
\def\Invalid@#1{\def#1{\Err@{\Invalid@@\string#1}}}
\def\Invalid@@{Invalid use of }
\message{figures,}
\Invalid@\caption
\Invalid@\captionwidth
\newdimen\smallcaptionwidth@
\def\topspace{\mid@false\ins@}
\def\midspace{\mid@true\ins@}
\newif\ifmid@
\def\captionfont@{}
\def\ins@#1{\relaxnext@\allowbreak
 \smallcaptionwidth@\captionwidth@\gdef\thespace@{#1}%
 \DN@{\ifx\next\space@\DN@. {\FN@\nextii@}\else
  \DN@.{\FN@\nextii@}\fi\next@.}%
 \DNii@{\ifx\next\caption\DN@\caption{\FN@\nextiii@}%
  \else\let\next@\nextiv@\fi\next@}%
 \def\nextiv@{\vnonvmode@
  {\ifmid@\expandafter\midinsert\else\expandafter\topinsert\fi
   \vbox to\thespace@{}\endinsert}
  {\ifmid@\nonvmodeerr@\midspace\else\nonvmodeerr@\topspace\fi}}%
 \def\nextiii@{\ifx\next\captionwidth\expandafter\nextv@
  \else\expandafter\nextvi@\fi}%
 \def\nextv@\captionwidth##1##2{\smallcaptionwidth@##1\relax\nextvi@{##2}}%
 \def\nextvi@##1{\def\thecaption@{\captionfont@##1}%
  \DN@{\ifx\next\space@\DN@. {\FN@\nextvii@}\else
   \DN@.{\FN@\nextvii@}\fi\next@.}%
  \FN@\next@}%
 \def\nextvii@{\vnonvmode@
  {\ifmid@\expandafter\midinsert\else
  \expandafter\topinsert\fi\vbox to\thespace@{}\nobreak\smallskip
  \setboxz@h{\noindent\ignorespaces\thecaption@\unskip}%
  \ifdim\wdz@>\smallcaptionwidth@\centerline{\vbox{\hsize\smallcaptionwidth@
   \noindent\ignorespaces\thecaption@\unskip}}%
  \else\centerline{\boxz@}\fi\endinsert}
  {\ifmid@\nonvmodeerr@\midspace
  \else\nonvmodeerr@\topspace\fi}}%
 \FN@\next@}
\message{comments,}
\def\newcodes@{\catcode`\\12\catcode`\{12\catcode`\}12\catcode`\#12%
 \catcode`\%12\relax}
\def\oldcodes@{\catcode`\\0\catcode`\{1\catcode`\}2\catcode`\#6%
 \catcode`\%14\relax}
\def\comment{\newcodes@\endlinechar=10 \comment@}
{\lccode`\0=`\\
\lowercase{\gdef\comment@#1^^J{\comment@@#10endcomment\comment@@@}%
\gdef\comment@@#10endcomment{\FN@\comment@@@}%
\gdef\comment@@@#1\comment@@@{\ifx\next\comment@@@\let\next\comment@
 \else\def\next{\oldcodes@\endlinechar=`\^^M\relax}%
 \fi\next}}}
\def\pr@m@s{\ifx'\next\DN@##1{\prim@s}\else\let\next@\egroup\fi\next@}
\def\prime{{\null\prime@\null}}
\mathchardef\prime@="0230
\let\dsize\displaystyle

\let\ssize\scriptstyle

\message{math spacing,}
\def\,{\RIfM@\mskip\thinmuskip\relax\else\kern.16667em\fi}
\def\!{\RIfM@\mskip-\thinmuskip\relax\else\kern-.16667em\fi}
\let\thinspace\,
\let\negthinspace\!
\def\medspace{\RIfM@\mskip\medmuskip\relax\else\kern.222222em\fi}
\def\negmedspace{\RIfM@\mskip-\medmuskip\relax\else\kern-.222222em\fi}
\def\thickspace{\RIfM@\mskip\thickmuskip\relax\else\kern.27777em\fi}
\let\;\thickspace
\def\negthickspace{\RIfM@\mskip-\thickmuskip\relax\else
 \kern-.27777em\fi}
\atdef@,{\RIfM@\mskip.1\thinmuskip\else\leavevmode\null,\fi}
\atdef@!{\RIfM@\mskip-.1\thinmuskip\else\leavevmode\null!\fi}
\atdef@.{\RIfM@&&\else\leavevmode.\spacefactor3000 \fi}
\def\and{\DOTSB\;\mathchar"3026 \;}

\message{fractions,}
\def\frac#1#2{{#1\over#2}}

\newdimen\ex@
\ex@.2326ex
\Invalid@\thickness
\def\thickfrac{\relaxnext@
 \DN@{\ifx\next\thickness\let\next@\nextii@\else
 \DN@{\nextii@\thickness1}\fi\next@}%
 \DNii@\thickness##1##2##3{{##2\above##1\ex@##3}}%
 \FN@\next@}

\def\thickfracwithdelims#1#2{\relaxnext@\def\ldelim@{#1}\def\rdelim@{#2}%
 \DN@{\ifx\next\thickness\let\next@\nextii@\else
 \DN@{\nextii@\thickness1}\fi\next@}%
 \DNii@\thickness##1##2##3{{##2\abovewithdelims
 \ldelim@\rdelim@##1\ex@##3}}%
 \FN@\next@}

\def\:{\nobreak\hskip.1111em\mathpunct{}\nonscript\mkern-\thinmuskip{:}\hskip
 .3333emplus.0555em\relax}
\def\snug{\unskip\kern-\mathsurround}
\message{smash commands,}
\def\topsmash{\top@true\bot@false\smash@}
\def\botsmash{\top@false\bot@true\smash@}
\newif\iftop@
\newif\ifbot@
\def\smash{\top@true\bot@true\smash@}
\def\smash@{\RIfM@\expandafter\mathpalette\expandafter\mathsm@sh\else
 \expandafter\makesm@sh\fi}
\def\finsm@sh{\iftop@\ht\z@\z@\fi\ifbot@\dp\z@\z@\fi\leavevmode\boxz@}
\message{large operator symbols,}
\def\LimitsOnSums{\global\let\slimits@\displaylimits}
\def\NoLimitsOnSums{\global\let\slimits@\nolimits}
\LimitsOnSums
\mathchardef\coprod@="1360       \def\coprod{\DOTSB\coprod@\slimits@}
\mathchardef\bigvee@="1357       \def\bigvee{\DOTSB\bigvee@\slimits@}
\mathchardef\bigwedge@="1356     \def\bigwedge{\DOTSB\bigwedge@\slimits@}
\mathchardef\biguplus@="1355     \def\biguplus{\DOTSB\biguplus@\slimits@}
\mathchardef\bigcap@="1354       \def\bigcap{\DOTSB\bigcap@\slimits@}
\mathchardef\bigcup@="1353       \def\bigcup{\DOTSB\bigcup@\slimits@}
\mathchardef\prod@="1351         \def\prod{\DOTSB\prod@\slimits@}
\mathchardef\sum@="1350          \def\sum{\DOTSB\sum@\slimits@}
\mathchardef\bigotimes@="134E    \def\bigotimes{\DOTSB\bigotimes@\slimits@}
\mathchardef\bigoplus@="134C     \def\bigoplus{\DOTSB\bigoplus@\slimits@}
\mathchardef\bigodot@="134A      \def\bigodot{\DOTSB\bigodot@\slimits@}
\mathchardef\bigsqcup@="1346     \def\bigsqcup{\DOTSB\bigsqcup@\slimits@}
\message{integrals,}
\def\LimitsOnInts{\global\let\ilimits@\displaylimits}
\def\NoLimitsOnInts{\global\let\ilimits@\nolimits}
\NoLimitsOnInts
\def\int{\DOTSI\intop\ilimits@}
\def\oint{\DOTSI\ointop\ilimits@}
\def\intic@{\mathchoice{\hskip.5em}{\hskip.4em}{\hskip.4em}{\hskip.4em}}
\def\negintic@{\mathchoice
 {\hskip-.5em}{\hskip-.4em}{\hskip-.4em}{\hskip-.4em}}
\def\intkern@{\mathchoice{\!\!\!}{\!\!}{\!\!}{\!\!}}
\def\intdots@{\mathchoice{\plaincdots@}
 {{\cdotp}\mkern1.5mu{\cdotp}\mkern1.5mu{\cdotp}}
 {{\cdotp}\mkern1mu{\cdotp}\mkern1mu{\cdotp}}
 {{\cdotp}\mkern1mu{\cdotp}\mkern1mu{\cdotp}}}
\newcount\intno@
\def\iint{\DOTSI\intno@\tw@\FN@\ints@}
\def\iiint{\DOTSI\intno@\thr@@\FN@\ints@}
\def\iiiint{\DOTSI\intno@4 \FN@\ints@}
\def\idotsint{\DOTSI\intno@\z@\FN@\ints@}
\def\ints@{\findlimits@\ints@@}
\newif\iflimtoken@
\newif\iflimits@
\def\findlimits@{\limtoken@true\ifx\next\limits\limits@true
 \else\ifx\next\nolimits\limits@false\else
 \limtoken@false\ifx\ilimits@\nolimits\limits@false\else
 \ifinner\limits@false\else\limits@true\fi\fi\fi\fi}
\def\multint@{\int\ifnum\intno@=\z@\intdots@                                
 \else\intkern@\fi                                                          
 \ifnum\intno@>\tw@\int\intkern@\fi                                         
 \ifnum\intno@>\thr@@\int\intkern@\fi                                       
 \int}                                                                      
\def\multintlimits@{\intop\ifnum\intno@=\z@\intdots@\else\intkern@\fi
 \ifnum\intno@>\tw@\intop\intkern@\fi
 \ifnum\intno@>\thr@@\intop\intkern@\fi\intop}
\def\ints@@{\iflimtoken@                                                    
 \def\ints@@@{\iflimits@\negintic@\mathop{\intic@\multintlimits@}\limits    
  \else\multint@\nolimits\fi                                                
  \eat@}                                                                    
 \else                                                                      
 \def\ints@@@{\iflimits@\negintic@
  \mathop{\intic@\multintlimits@}\limits\else
  \multint@\nolimits\fi}\fi\ints@@@}
\def\LimitsOnNames{\global\let\nlimits@\displaylimits}
\def\NoLimitsOnNames{\global\let\nlimits@\nolimits@}
\LimitsOnNames
\def\nolimits@{\relaxnext@
 \DN@{\ifx\next\limits\DN@\limits{\nolimits}\else
  \let\next@\nolimits\fi\next@}%
 \FN@\next@}
\message{operator names,}
\def\newmcodes@{\mathcode`\'"27\mathcode`\*"2A\mathcode`\."613A%
 \mathcode`\-"2D\mathcode`\/"2F\mathcode`\:"603A }
\def\operatorname#1{\mathop{\newmcodes@\kern\z@\fam\z@#1}\nolimits@}
\def\operatornamewithlimits#1{\mathop{\newmcodes@\kern\z@\fam\z@#1}\nlimits@}
\def\qopname@#1{\mathop{\fam\z@#1}\nolimits@}
\def\qopnamewl@#1{\mathop{\fam\z@#1}\nlimits@}
\def\arccos{\qopname@{arccos}}
\def\arcsin{\qopname@{arcsin}}
\def\arctan{\qopname@{arctan}}
\def\arg{\qopname@{arg}}
\def\cos{\qopname@{cos}}
\def\cosh{\qopname@{cosh}}
\def\cot{\qopname@{cot}}
\def\coth{\qopname@{coth}}
\def\csc{\qopname@{csc}}
\def\deg{\qopname@{deg}}
\def\det{\qopnamewl@{det}}
\def\dim{\qopname@{dim}}
\def\exp{\qopname@{exp}}
\def\gcd{\qopnamewl@{gcd}}
\def\hom{\qopname@{hom}}
\def\inf{\qopnamewl@{inf}}
\def\injlim{\qopnamewl@{inj\,lim}}
\def\ker{\qopname@{ker}}
\def\lg{\qopname@{lg}}
\def\lim{\qopnamewl@{lim}}
\def\liminf{\qopnamewl@{lim\,inf}}
\def\limsup{\qopnamewl@{lim\,sup}}
\def\ln{\qopname@{ln}}
\def\log{\qopname@{log}}
\def\max{\qopnamewl@{max}}
\def\min{\qopnamewl@{min}}
\def\Pr{\qopnamewl@{Pr}}
\def\projlim{\qopnamewl@{proj\,lim}}
\def\sec{\qopname@{sec}}
\def\sin{\qopname@{sin}}
\def\sinh{\qopname@{sinh}}
\def\sup{\qopnamewl@{sup}}
\def\tan{\qopname@{tan}}
\def\tanh{\qopname@{tanh}}
\def\varinjlim{\mathop{\vtop{\ialign{##\crcr
 \hfil\rm lim\hfil\crcr\noalign{\nointerlineskip}\rightarrowfill\crcr
 \noalign{\nointerlineskip\kern-\ex@}\crcr}}}}
\def\varprojlim{\mathop{\vtop{\ialign{##\crcr
 \hfil\rm lim\hfil\crcr\noalign{\nointerlineskip}\leftarrowfill\crcr
 \noalign{\nointerlineskip\kern-\ex@}\crcr}}}}
\def\varliminf{\mathop{\underline{\vrule height\z@ depth.2exwidth\z@
 \hbox{\rm lim}}}}

\newdimen\buffer@
\buffer@\fontdimen13 \tenex
\newdimen\buffer
\buffer\buffer@

\def\ResetBuffer{\fontdimen13 \tenex\buffer@\global\buffer\buffer@}
\def\shave#1{\mathop{\hbox{$\m@th\fontdimen13 \tenex\z@                     
 \displaystyle{#1}$}}\fontdimen13 \tenex\buffer}

\message{multilevel sub/superscripts,}
\Invalid@\\
\def\Let@{\relax\iffalse{\fi\let\\=\cr\iffalse}\fi}
\Invalid@\vspace
\def\vspace@{\def\vspace##1{\crcr\noalign{\vskip##1\relax}}}
\def\multilimits@{\bgroup\vspace@\Let@
 \baselineskip\fontdimen10 \scriptfont\tw@
 \advance\baselineskip\fontdimen12 \scriptfont\tw@
 \lineskip\thr@@\fontdimen8 \scriptfont\thr@@
 \lineskiplimit\lineskip
 \vbox\bgroup\ialign\bgroup\hfil$\m@th\scriptstyle{##}$\hfil\crcr}
\def\Sb{_\multilimits@}
\def\endSb{\crcr\egroup\egroup\egroup}
\def\Sp{^\multilimits@}

\def\spreadlines#1{\RIfMIfI@\onlydmatherr@\spreadlines\else
 \openup#1\relax\fi\else\onlydmatherr@\spreadlines\fi}
\def\Mathstrut@{\copy\Mathstrutbox@}
\newbox\Mathstrutbox@
\setbox\Mathstrutbox@\null
\setboxz@h{$\m@th($}
\ht\Mathstrutbox@\ht\z@
\dp\Mathstrutbox@\dp\z@
\message{matrices,}
\newdimen\spreadmlines@
\def\spreadmatrixlines#1{\RIfMIfI@
 \onlydmatherr@\spreadmatrixlines\else
 \spreadmlines@#1\relax\fi\else\onlydmatherr@\spreadmatrixlines\fi}
\def\matrix{\null\,\vcenter\bgroup\Let@\vspace@
 \normalbaselines\openup\spreadmlines@\ialign
 \bgroup\hfil$\m@th##$\hfil&&\quad\hfil$\m@th##$\hfil\crcr
 \Mathstrut@\crcr\noalign{\kern-\baselineskip}}
\def\endmatrix{\crcr\Mathstrut@\crcr\noalign{\kern-\baselineskip}\egroup
 \egroup\,}
\def\format{\crcr\egroup\iffalse{\fi\ifnum`}=0 \fi\format@}
\newtoks\hashtoks@
\hashtoks@{#}
\def\format@#1\\{\def\preamble@{#1}%
 \def\l{$\m@th\the\hashtoks@$\hfil}%
 \def\c{\hfil$\m@th\the\hashtoks@$\hfil}%
 \def\r{\hfil$\m@th\the\hashtoks@$}%
 \edef\preamble@@{\preamble@}\ifnum`{=0 \fi\iffalse}\fi
 \ialign\bgroup\span\preamble@@\crcr}
\def\smallmatrix{\null\,\vcenter\bgroup\vspace@\Let@
 \baselineskip9\ex@\lineskip\ex@
 \ialign\bgroup\hfil$\m@th\scriptstyle{##}$\hfil&&\thickspace\hfil
 $\m@th\scriptstyle{##}$\hfil\crcr}
\def\endsmallmatrix{\crcr\egroup\egroup\,}

\newmuskip\dotsspace@
\dotsspace@1.5mu
\def\strip@#1 {#1}
\def\spacehdots#1\for#2{\multispan{#2}\xleaders
 \hbox{$\m@th\mkern\strip@#1 \dotsspace@.\mkern\strip@#1 \dotsspace@$}\hfill}
\def\hdotsfor#1{\spacehdots\@ne\for{#1}}
\def\multispan@#1{\omit\mscount#1\unskip\loop\ifnum\mscount>\@ne\sp@n\repeat}
\def\spaceinnerhdots#1\for#2\after#3{\multispan@{\strip@#2 }#3\xleaders
 \hbox{$\m@th\mkern\strip@#1 \dotsspace@.\mkern\strip@#1 \dotsspace@$}\hfill}
\def\innerhdotsfor#1\after#2{\spaceinnerhdots\@ne\for#1\after{#2}}
\def\cases{\bgroup\spreadmlines@\jot\left\{\,\matrix\format\l&\quad\l\\}
\def\endcases{\endmatrix\right.\egroup}
\message{multiline displays,}
\newif\ifinany@
\newif\ifinalign@
\newif\ifingather@
\def\strut@{\copy\strutbox@}
\newbox\strutbox@
\setbox\strutbox@\hbox{\vrule height8\p@ depth3\p@ width\z@}
\def\topaligned{\null\,\vtop\aligned@}
\def\botaligned{\null\,\vbox\aligned@}
\def\aligned{\null\,\vcenter\aligned@}
\def\aligned@{\bgroup\vspace@\Let@
 \ifinany@\else\openup\jot\fi\ialign
 \bgroup\hfil\strut@$\m@th\displaystyle{##}$&
 $\m@th\displaystyle{{}##}$\hfil\crcr}
\def\endaligned{\crcr\egroup\egroup}

\def\alignedat#1{\null\,\vcenter\bgroup\doat@{#1}\vspace@\Let@
 \ifinany@\else\openup\jot\fi\ialign\bgroup\span\preamble@@\crcr}
\newcount\atcount@
\def\doat@#1{\toks@{\hfil\strut@$\m@th
 \displaystyle{\the\hashtoks@}$&$\m@th\displaystyle
 {{}\the\hashtoks@}$\hfil}
 \atcount@#1\relax\advance\atcount@\m@ne                                    
 \loop\ifnum\atcount@>\z@\toks@=\expandafter{\the\toks@&\hfil$\m@th
 \displaystyle{\the\hashtoks@}$&$\m@th
 \displaystyle{{}\the\hashtoks@}$\hfil}\advance
  \atcount@\m@ne\repeat                                                     
 \xdef\preamble@{\the\toks@}\xdef\preamble@@{\preamble@}}

\def\gathered{\null\,\vcenter\bgroup\vspace@\Let@
 \ifinany@\else\openup\jot\fi\ialign
 \bgroup\hfil\strut@$\m@th\displaystyle{##}$\hfil\crcr}
\def\endgathered{\crcr\egroup\egroup}
\newif\iftagsleft@
\def\TagsOnLeft{\global\tagsleft@true}
\def\TagsOnRight{\global\tagsleft@false}
\TagsOnLeft
\newif\ifmathtags@
\def\TagsAsMath{\global\mathtags@true}
\def\TagsAsText{\global\mathtags@false}
\TagsAsText
\def\tagform@#1{\hbox{\rm(\ignorespaces#1\unskip)}}
\def\thetag{\leavevmode\tagform@}
\def\tag#1$${\iftagsleft@\leqno\else\eqno\fi                                
 \maketag@#1\maketag@                                                       
 $$}                                                                        
\def\maketag@{\FN@\maketag@@}
\def\maketag@@{\ifx\next"\expandafter\maketag@@@\else\expandafter\maketag@@@@
 \fi}
\def\maketag@@@"#1"#2\maketag@{\hbox{\rm#1}}                                
\def\maketag@@@@#1\maketag@{\ifmathtags@\tagform@{$\m@th#1$}\else
 \tagform@{#1}\fi}
\interdisplaylinepenalty\@M
\def\allowdisplaybreaks{\RIfMIfI@
 \onlydmatherr@\allowdisplaybreaks\else
 \interdisplaylinepenalty\z@\fi\else\onlydmatherr@\allowdisplaybreaks\fi}
\Invalid@\allowdisplaybreak
\Invalid@\displaybreak
\Invalid@\intertext
\def\allowdisplaybreak@{\def\allowdisplaybreak{\crcr\noalign{\allowbreak}}}
\def\displaybreak@{\def\displaybreak{\crcr\noalign{\break}}}
\def\intertext@{\def\intertext##1{\crcr\noalign{%
 \penalty\postdisplaypenalty \vskip\belowdisplayskip
 \vbox{\normalbaselines\noindent##1}%
 \penalty\predisplaypenalty \vskip\abovedisplayskip}}}
\newskip\centering@
\centering@\z@ plus\@m\p@
\def\align{\relax\ifingather@\DN@{\csname align (in
  \string\gather)\endcsname}\else
 \ifmmode\ifinner\DN@{\onlydmatherr@\align}\else
  \let\next@\align@\fi
 \else\DN@{\onlydmatherr@\align}\fi\fi\next@}
\newhelp\andhelp@
{An extra & here is so disastrous that you should probably exit^^J
and fix things up.}
\newif\iftag@
\newcount\and@
\def\align@{\inalign@true\inany@true
 \vspace@\allowdisplaybreak@\displaybreak@\intertext@
 \def\tag{\global\tag@true\ifnum\and@=\z@\DN@{&&}\else
          \DN@{&}\fi\next@}%
 \iftagsleft@\DN@{\csname align \endcsname}\else
  \DN@{\csname align \space\endcsname}\fi\next@}
\def\Tag@{\iftag@\else\errhelp\andhelp@\err@{Extra & on this line}\fi}
\newdimen\lwidth@
\newdimen\rwidth@
\newdimen\maxlwidth@
\newdimen\maxrwidth@
\newdimen\totwidth@
\def\measure@#1\endalign{\lwidth@\z@\rwidth@\z@\maxlwidth@\z@\maxrwidth@\z@
 \global\and@\z@                                                            
 \setbox@ne\vbox                                                            
  {\everycr{\noalign{\global\tag@false\global\and@\z@}}\Let@                
  \halign{\setboxz@h{$\m@th\displaystyle{\@lign##}$}
   \global\lwidth@\wdz@                                                     
   \ifdim\lwidth@>\maxlwidth@\global\maxlwidth@\lwidth@\fi                  
   \global\advance\and@\@ne                                                 
   &\setboxz@h{$\m@th\displaystyle{{}\@lign##}$}\global\rwidth@\wdz@        
   \ifdim\rwidth@>\maxrwidth@\global\maxrwidth@\rwidth@\fi                  
   \global\advance\and@\@ne                                                
   &\Tag@
   \eat@{##}\crcr#1\crcr}}
 \totwidth@\maxlwidth@\advance\totwidth@\maxrwidth@}                       
\def\displ@y@{\global\dt@ptrue\openup\jot
 \everycr{\noalign{\global\tag@false\global\and@\z@\ifdt@p\global\dt@pfalse
 \vskip-\lineskiplimit\vskip\normallineskiplimit\else
 \penalty\interdisplaylinepenalty\fi}}}
\def\black@#1{\noalign{\ifdim#1>\displaywidth
 \dimen@\prevdepth\nointerlineskip                                          
 \vskip-\ht\strutbox@\vskip-\dp\strutbox@                                   
 \vbox{\noindent\hbox to#1{\strut@\hfill}}
 \prevdepth\dimen@                                                          
 \fi}}
\expandafter\def\csname align \space\endcsname#1\endalign
 {\measure@#1\endalign\global\and@\z@                                       
 \ifingather@\everycr{\noalign{\global\and@\z@}}\else\displ@y@\fi           
 \Let@\tabskip\centering@                                                   
 \halign to\displaywidth
  {\hfil\strut@\setboxz@h{$\m@th\displaystyle{\@lign##}$}
  \global\lwidth@\wdz@\boxz@\global\advance\and@\@ne                        
  \tabskip\z@skip                                                           
  &\setboxz@h{$\m@th\displaystyle{{}\@lign##}$}
  \global\rwidth@\wdz@\boxz@\hfill\global\advance\and@\@ne                  
  \tabskip\centering@                                                       
  &\setboxz@h{\@lign\strut@\maketag@##\maketag@}
  \dimen@\displaywidth\advance\dimen@-\totwidth@
  \divide\dimen@\tw@\advance\dimen@\maxrwidth@\advance\dimen@-\rwidth@     
  \ifdim\dimen@<\tw@\wdz@\llap{\vtop{\normalbaselines\null\boxz@}}
  \else\llap{\boxz@}\fi                                                    
  \tabskip\z@skip                                                          
  \crcr#1\crcr                                                             
  \black@\totwidth@}}                                                      
\newdimen\lineht@
\expandafter\def\csname align \endcsname#1\endalign{\measure@#1\endalign
 \global\and@\z@
 \ifdim\totwidth@>\displaywidth\let\displaywidth@\totwidth@\else
  \let\displaywidth@\displaywidth\fi                                        
 \ifingather@\everycr{\noalign{\global\and@\z@}}\else\displ@y@\fi
 \Let@\tabskip\centering@\halign to\displaywidth
  {\hfil\strut@\setboxz@h{$\m@th\displaystyle{\@lign##}$}%
  \global\lwidth@\wdz@\global\lineht@\ht\z@                                 
  \boxz@\global\advance\and@\@ne
  \tabskip\z@skip&\setboxz@h{$\m@th\displaystyle{{}\@lign##}$}%
  \global\rwidth@\wdz@\ifdim\ht\z@>\lineht@\global\lineht@\ht\z@\fi         
  \boxz@\hfil\global\advance\and@\@ne
  \tabskip\centering@&\kern-\displaywidth@                                  
  \setboxz@h{\@lign\strut@\maketag@##\maketag@}%
  \dimen@\displaywidth\advance\dimen@-\totwidth@
  \divide\dimen@\tw@\advance\dimen@\maxlwidth@\advance\dimen@-\lwidth@
  \ifdim\dimen@<\tw@\wdz@
   \rlap{\vbox{\normalbaselines\boxz@\vbox to\lineht@{}}}\else
   \rlap{\boxz@}\fi
  \tabskip\displaywidth@\crcr#1\crcr\black@\totwidth@}}
\expandafter\def\csname align (in \string\gather)\endcsname
  #1\endalign{\vcenter{\align@#1\endalign}}
\Invalid@\endalign
\newif\ifxat@
\def\alignat{\RIfMIfI@\DN@{\onlydmatherr@\alignat}\else
 \DN@{\csname alignat \endcsname}\fi\else
 \DN@{\onlydmatherr@\alignat}\fi\next@}
\newif\ifmeasuring@
\newbox\savealignat@
\expandafter\def\csname alignat \endcsname#1#2\endalignat                   
 {\inany@true\xat@false
 \def\tag{\global\tag@true\count@#1\relax\multiply\count@\tw@
  \xdef\tag@{}\loop\ifnum\count@>\and@\xdef\tag@{&\tag@}\advance\count@\m@ne
  \repeat\tag@}%
 \vspace@\allowdisplaybreak@\displaybreak@\intertext@
 \displ@y@\measuring@true                                                   
 \setbox\savealignat@\hbox{$\m@th\displaystyle\Let@
  \attag@{#1}
  \vbox{\halign{\span\preamble@@\crcr#2\crcr}}$}%
 \measuring@false                                                           
 \Let@\attag@{#1}
 \tabskip\centering@\halign to\displaywidth
  {\span\preamble@@\crcr#2\crcr                                             
  \black@{\wd\savealignat@}}}                                               
\Invalid@\endalignat
\def\xalignat{\RIfMIfI@
 \DN@{\onlydmatherr@\xalignat}\else
 \DN@{\csname xalignat \endcsname}\fi\else
 \DN@{\onlydmatherr@\xalignat}\fi\next@}
\expandafter\def\csname xalignat \endcsname#1#2\endxalignat
 {\inany@true\xat@true
 \def\tag{\global\tag@true\def\tag@{}\count@#1\relax\multiply\count@\tw@
  \loop\ifnum\count@>\and@\xdef\tag@{&\tag@}\advance\count@\m@ne\repeat\tag@}%
 \vspace@\allowdisplaybreak@\displaybreak@\intertext@
 \displ@y@\measuring@true\setbox\savealignat@\hbox{$\m@th\displaystyle\Let@
 \attag@{#1}\vbox{\halign{\span\preamble@@\crcr#2\crcr}}$}%
 \measuring@false\Let@
 \attag@{#1}\tabskip\centering@\halign to\displaywidth
 {\span\preamble@@\crcr#2\crcr\black@{\wd\savealignat@}}}
\def\attag@#1{\let\Maketag@\maketag@\let\TAG@\Tag@                          
 \let\Tag@=0\let\maketag@=0
 \ifmeasuring@\def\llap@##1{\setboxz@h{##1}\hbox to\tw@\wdz@{}}%
  \def\rlap@##1{\setboxz@h{##1}\hbox to\tw@\wdz@{}}\else
  \let\llap@\llap\let\rlap@\rlap\fi                                         
 \toks@{\hfil\strut@$\m@th\displaystyle{\@lign\the\hashtoks@}$\tabskip\z@skip
  \global\advance\and@\@ne&$\m@th\displaystyle{{}\@lign\the\hashtoks@}$\hfil
  \ifxat@\tabskip\centering@\fi\global\advance\and@\@ne}
 \iftagsleft@
  \toks@@{\tabskip\centering@&\Tag@\kern-\displaywidth
   \rlap@{\@lign\maketag@\the\hashtoks@\maketag@}%
   \global\advance\and@\@ne\tabskip\displaywidth}\else
  \toks@@{\tabskip\centering@&\Tag@\llap@{\@lign\maketag@
   \the\hashtoks@\maketag@}\global\advance\and@\@ne\tabskip\z@skip}\fi      
 \atcount@#1\relax\advance\atcount@\m@ne
 \loop\ifnum\atcount@>\z@
 \toks@=\expandafter{\the\toks@&\hfil$\m@th\displaystyle{\@lign
  \the\hashtoks@}$\global\advance\and@\@ne
  \tabskip\z@skip&$\m@th\displaystyle{{}\@lign\the\hashtoks@}$\hfil\ifxat@
  \tabskip\centering@\fi\global\advance\and@\@ne}\advance\atcount@\m@ne
 \repeat                                                                    
 \xdef\preamble@{\the\toks@\the\toks@@}
 \xdef\preamble@@{\preamble@}
 \let\maketag@\Maketag@\let\Tag@\TAG@}                                      
\Invalid@\endxalignat
\def\xxalignat{\RIfMIfI@
 \DN@{\onlydmatherr@\xxalignat}\else\DN@{\csname xxalignat
  \endcsname}\fi\else
 \DN@{\onlydmatherr@\xxalignat}\fi\next@}
\expandafter\def\csname xxalignat \endcsname#1#2\endxxalignat{\inany@true
 \vspace@\allowdisplaybreak@\displaybreak@\intertext@
 \displ@y\setbox\savealignat@\hbox{$\m@th\displaystyle\Let@
 \xxattag@{#1}\vbox{\halign{\span\preamble@@\crcr#2\crcr}}$}%
 \Let@\xxattag@{#1}\tabskip\z@skip\halign to\displaywidth
 {\span\preamble@@\crcr#2\crcr\black@{\wd\savealignat@}}}
\def\xxattag@#1{\toks@{\tabskip\z@skip\hfil\strut@
 $\m@th\displaystyle{\the\hashtoks@}$&%
 $\m@th\displaystyle{{}\the\hashtoks@}$\hfil\tabskip\centering@&}%
 \atcount@#1\relax\advance\atcount@\m@ne\loop\ifnum\atcount@>\z@
 \toks@=\expandafter{\the\toks@&\hfil$\m@th\displaystyle{\the\hashtoks@}$%
  \tabskip\z@skip&$\m@th\displaystyle{{}\the\hashtoks@}$\hfil
  \tabskip\centering@}\advance\atcount@\m@ne\repeat
 \xdef\preamble@{\the\toks@\tabskip\z@skip}\xdef\preamble@@{\preamble@}}
\Invalid@\endxxalignat
\newdimen\gwidth@
\newdimen\gmaxwidth@
\def\gmeasure@#1\endgather{\gwidth@\z@\gmaxwidth@\z@\setbox@ne\vbox{\Let@
 \halign{\setboxz@h{$\m@th\displaystyle{##}$}\global\gwidth@\wdz@
 \ifdim\gwidth@>\gmaxwidth@\global\gmaxwidth@\gwidth@\fi
 &\eat@{##}\crcr#1\crcr}}}
\def\gather{\RIfMIfI@\DN@{\onlydmatherr@\gather}\else
 \ingather@true\inany@true\def\tag{&}%
 \vspace@\allowdisplaybreak@\displaybreak@\intertext@
 \displ@y\Let@
 \iftagsleft@\DN@{\csname gather \endcsname}\else
  \DN@{\csname gather \space\endcsname}\fi\fi
 \else\DN@{\onlydmatherr@\gather}\fi\next@}
\expandafter\def\csname gather \space\endcsname#1\endgather
 {\gmeasure@#1\endgather\tabskip\centering@
 \halign to\displaywidth{\hfil\strut@\setboxz@h{$\m@th\displaystyle{##}$}%
 \global\gwidth@\wdz@\boxz@\hfil&
 \setboxz@h{\strut@{\maketag@##\maketag@}}%
 \dimen@\displaywidth\advance\dimen@-\gwidth@
 \ifdim\dimen@>\tw@\wdz@\llap{\boxz@}\else
 \llap{\vtop{\normalbaselines\null\boxz@}}\fi
 \tabskip\z@skip\crcr#1\crcr\black@\gmaxwidth@}}
\newdimen\glineht@
\expandafter\def\csname gather \endcsname#1\endgather{\gmeasure@#1\endgather
 \ifdim\gmaxwidth@>\displaywidth\let\gdisplaywidth@\gmaxwidth@\else
 \let\gdisplaywidth@\displaywidth\fi\tabskip\centering@\halign to\displaywidth
 {\hfil\strut@\setboxz@h{$\m@th\displaystyle{##}$}%
 \global\gwidth@\wdz@\global\glineht@\ht\z@\boxz@\hfil&\kern-\gdisplaywidth@
 \setboxz@h{\strut@{\maketag@##\maketag@}}%
 \dimen@\displaywidth\advance\dimen@-\gwidth@
 \ifdim\dimen@>\tw@\wdz@\rlap{\boxz@}\else
 \rlap{\vbox{\normalbaselines\boxz@\vbox to\glineht@{}}}\fi
 \tabskip\gdisplaywidth@\crcr#1\crcr\black@\gmaxwidth@}}
\newif\ifctagsplit@
\def\CenteredTagsOnSplits{\global\ctagsplit@true}
\def\TopOrBottomTagsOnSplits{\global\ctagsplit@false}
\TopOrBottomTagsOnSplits
\def\split{\relax\ifinany@\let\next@\insplit@\else
 \ifmmode\ifinner\def\next@{\onlydmatherr@\split}\else
 \let\next@\outsplit@\fi\else
 \def\next@{\onlydmatherr@\split}\fi\fi\next@}
\def\insplit@{\global\setbox\z@\vbox\bgroup\vspace@\Let@\ialign\bgroup
 \hfil\strut@$\m@th\displaystyle{##}$&$\m@th\displaystyle{{}##}$\hfill\crcr}
\def\endsplit{\crcr\egroup\egroup\iftagsleft@\expandafter\lendsplit@\else
 \expandafter\rendsplit@\fi}
\def\rendsplit@{\global\setbox9 \vbox
 {\unvcopy\z@\global\setbox8 \lastbox\unskip}
 \setbox@ne\hbox{\unhcopy8 \unskip\global\setbox\tw@\lastbox
 \unskip\global\setbox\thr@@\lastbox}
 \global\setbox7 \hbox{\unhbox\tw@\unskip}
 \ifinalign@\ifctagsplit@                                                   
  \gdef\split@{\hbox to\wd\thr@@{}&
   \vcenter{\vbox{\moveleft\wd\thr@@\boxz@}}}
 \else\gdef\split@{&\vbox{\moveleft\wd\thr@@\box9}\crcr
  \box\thr@@&\box7}\fi                                                      
 \else                                                                      
  \ifctagsplit@\gdef\split@{\vcenter{\boxz@}}\else
  \gdef\split@{\box9\crcr\hbox{\box\thr@@\box7}}\fi
 \fi
 \split@}                                                                   
\def\lendsplit@{\global\setbox9\vtop{\unvcopy\z@}
 \setbox@ne\vbox{\unvcopy\z@\global\setbox8\lastbox}
 \setbox@ne\hbox{\unhcopy8\unskip\setbox\tw@\lastbox
  \unskip\global\setbox\thr@@\lastbox}
 \ifinalign@\ifctagsplit@                                                   
  \gdef\split@{\hbox to\wd\thr@@{}&
  \vcenter{\vbox{\moveleft\wd\thr@@\box9}}}
  \else                                                                     
  \gdef\split@{\hbox to\wd\thr@@{}&\vbox{\moveleft\wd\thr@@\box9}}\fi
 \else
  \ifctagsplit@\gdef\split@{\vcenter{\box9}}\else
  \gdef\split@{\box9}\fi
 \fi\split@}
\def\outsplit@#1$${\align\insplit@#1\endalign$$}
\newdimen\multlinegap@
\multlinegap@1em
\newdimen\multlinetaggap@
\multlinetaggap@1em
\def\MultlineGap#1{\global\multlinegap@#1\relax}
\def\multlinegap#1{\RIfMIfI@\onlydmatherr@\multlinegap\else
 \multlinegap@#1\relax\fi\else\onlydmatherr@\multlinegap\fi}
\def\nomultlinegap{\multlinegap{\z@}}
\def\multline{\RIfMIfI@
 \DN@{\onlydmatherr@\multline}\else
 \DN@{\multline@}\fi\else
 \DN@{\onlydmatherr@\multline}\fi\next@}
\newif\iftagin@
\def\tagin@#1{\tagin@false\in@\tag{#1}\ifin@\tagin@true\fi}
\def\multline@#1$${\inany@true\vspace@\allowdisplaybreak@\displaybreak@
 \tagin@{#1}\iftagsleft@\DN@{\multline@l#1$$}\else
 \DN@{\multline@r#1$$}\fi\next@}
\newdimen\mwidth@
\def\rmmeasure@#1\endmultline{%
 \def\shoveleft##1{##1}\def\shoveright##1{##1}
 \setbox@ne\vbox{\Let@\halign{\setboxz@h
  {$\m@th\@lign\displaystyle{}##$}\global\mwidth@\wdz@
  \crcr#1\crcr}}}
\newdimen\mlineht@
\newif\ifzerocr@
\newif\ifonecr@
\def\lmmeasure@#1\endmultline{\global\zerocr@true\global\onecr@false
 \everycr{\noalign{\ifonecr@\global\onecr@false\fi
  \ifzerocr@\global\zerocr@false\global\onecr@true\fi}}
  \def\shoveleft##1{##1}\def\shoveright##1{##1}%
 \setbox@ne\vbox{\Let@\halign{\setboxz@h
  {$\m@th\@lign\displaystyle{}##$}\ifonecr@\global\mwidth@\wdz@
  \global\mlineht@\ht\z@\fi\crcr#1\crcr}}}
\newbox\mtagbox@
\newdimen\ltwidth@
\newdimen\rtwidth@
\def\multline@l#1$${\iftagin@\DN@{\lmultline@@#1$$}\else
 \DN@{\setbox\mtagbox@\null\ltwidth@\z@\rtwidth@\z@
  \lmultline@@@#1$$}\fi\next@}
\def\lmultline@@#1\endmultline\tag#2$${%
 \setbox\mtagbox@\hbox{\maketag@#2\maketag@}
 \lmmeasure@#1\endmultline\dimen@\mwidth@\advance\dimen@\wd\mtagbox@
 \advance\dimen@\multlinetaggap@                                            
 \ifdim\dimen@>\displaywidth\ltwidth@\z@\else\ltwidth@\wd\mtagbox@\fi       
 \lmultline@@@#1\endmultline$$}
\def\lmultline@@@{\displ@y
 \def\shoveright##1{##1\hfilneg\hskip\multlinegap@}%
 \def\shoveleft##1{\setboxz@h{$\m@th\displaystyle{}##1$}%
  \setbox@ne\hbox{$\m@th\displaystyle##1$}%
  \hfilneg
  \iftagin@
   \ifdim\ltwidth@>\z@\hskip\ltwidth@\hskip\multlinetaggap@\fi
  \else\hskip\multlinegap@\fi\hskip.5\wd@ne\hskip-.5\wdz@##1}
  \halign\bgroup\Let@\hbox to\displaywidth
   {\strut@$\m@th\displaystyle\hfil{}##\hfil$}\crcr
   \hfilneg                                                                 
   \iftagin@                                                                
    \ifdim\ltwidth@>\z@                                                     
     \box\mtagbox@\hskip\multlinetaggap@                                    
    \else
     \rlap{\vbox{\normalbaselines\hbox{\strut@\box\mtagbox@}%
     \vbox to\mlineht@{}}}\fi                                               
   \else\hskip\multlinegap@\fi}                                             
\def\multline@r#1$${\iftagin@\DN@{\rmultline@@#1$$}\else
 \DN@{\setbox\mtagbox@\null\ltwidth@\z@\rtwidth@\z@
  \rmultline@@@#1$$}\fi\next@}
\def\rmultline@@#1\endmultline\tag#2$${\ltwidth@\z@
 \setbox\mtagbox@\hbox{\maketag@#2\maketag@}%
 \rmmeasure@#1\endmultline\dimen@\mwidth@\advance\dimen@\wd\mtagbox@
 \advance\dimen@\multlinetaggap@
 \ifdim\dimen@>\displaywidth\rtwidth@\z@\else\rtwidth@\wd\mtagbox@\fi
 \rmultline@@@#1\endmultline$$}
\def\rmultline@@@{\displ@y
 \def\shoveright##1{##1\hfilneg\iftagin@\ifdim\rtwidth@>\z@
  \hskip\rtwidth@\hskip\multlinetaggap@\fi\else\hskip\multlinegap@\fi}%
 \def\shoveleft##1{\setboxz@h{$\m@th\displaystyle{}##1$}%
  \setbox@ne\hbox{$\m@th\displaystyle##1$}%
  \hfilneg\hskip\multlinegap@\hskip.5\wd@ne\hskip-.5\wdz@##1}%
 \halign\bgroup\Let@\hbox to\displaywidth
  {\strut@$\m@th\displaystyle\hfil{}##\hfil$}\crcr
 \hfilneg\hskip\multlinegap@}
\def\endmultline{\iftagsleft@\expandafter\lendmultline@\else
 \expandafter\rendmultline@\fi}
\def\lendmultline@{\hfilneg\hskip\multlinegap@\crcr\egroup}
\def\rendmultline@{\iftagin@                                                
 \ifdim\rtwidth@>\z@                                                        
  \hskip\multlinetaggap@\box\mtagbox@                                       
 \else\llap{\vtop{\normalbaselines\null\hbox{\strut@\box\mtagbox@}}}\fi     
 \else\hskip\multlinegap@\fi                                                
 \hfilneg\crcr\egroup}
\def\bmod{\mskip-\medmuskip\mkern5mu\mathbin{\fam\z@ mod}\penalty900
 \mkern5mu\mskip-\medmuskip}
\def\pmod#1{\allowbreak\ifinner\mkern8mu\else\mkern18mu\fi
 ({\fam\z@ mod}\,\,#1)}
\def\pod#1{\allowbreak\ifinner\mkern8mu\else\mkern18mu\fi(#1)}
\def\mod#1{\allowbreak\ifinner\mkern12mu\else\mkern18mu\fi{\fam\z@ mod}\,\,#1}
\message{continued fractions,}
\newcount\cfraccount@
\def\cfrac{\bgroup\bgroup\advance\cfraccount@\@ne\strut
 \iffalse{\fi\def\\{\over\displaystyle}\iffalse}\fi}
\def\lcfrac{\bgroup\bgroup\advance\cfraccount@\@ne\strut
 \iffalse{\fi\def\\{\hfill\over\displaystyle}\iffalse}\fi}
\def\rcfrac{\bgroup\bgroup\advance\cfraccount@\@ne\strut\hfill
 \iffalse{\fi\def\\{\over\displaystyle}\iffalse}\fi}
\def\gloop@#1\repeat{\gdef\body{#1}\iterate}
\def\endcfrac{\gloop@\ifnum\cfraccount@>\z@\global\advance\cfraccount@\m@ne
 \egroup\hskip-\nulldelimiterspace\egroup\repeat}
\message{compound symbols,}
\def\binrel@#1{\setboxz@h{\thinmuskip0mu
  \medmuskip\m@ne mu\thickmuskip\@ne mu$#1\m@th$}%
 \setbox@ne\hbox{\thinmuskip0mu\medmuskip\m@ne mu\thickmuskip
  \@ne mu${}#1{}\m@th$}%
 \setbox\tw@\hbox{\hskip\wd@ne\hskip-\wdz@}}
\def\overset#1\to#2{\binrel@{#2}\ifdim\wd\tw@<\z@
 \mathbin{\mathop{\kern\z@#2}\limits^{#1}}\else\ifdim\wd\tw@>\z@
 \mathrel{\mathop{\kern\z@#2}\limits^{#1}}\else
 {\mathop{\kern\z@#2}\limits^{#1}}{}\fi\fi}
\def\underset#1\to#2{\binrel@{#2}\ifdim\wd\tw@<\z@
 \mathbin{\mathop{\kern\z@#2}\limits_{#1}}\else\ifdim\wd\tw@>\z@
 \mathrel{\mathop{\kern\z@#2}\limits_{#1}}\else
 {\mathop{\kern\z@#2}\limits_{#1}}{}\fi\fi}
\def\oversetbrace#1\to#2{\overbrace{#2}^{#1}}
\def\undersetbrace#1\to#2{\underbrace{#2}_{#1}}
\def\sideset#1\and#2\to#3{%
 \setbox@ne\hbox{$\dsize{\vphantom{#3}}#1{#3}\m@th$}%
 \setbox\tw@\hbox{$\dsize{#3}#2\m@th$}%
 \hskip\wd@ne\hskip-\wd\tw@\mathop{\hskip\wd\tw@\hskip-\wd@ne
  {\vphantom{#3}}#1{#3}#2}}
\def\rightarrowfill@#1{\setboxz@h{$#1-\m@th$}\ht\z@\z@
  $#1\m@th\copy\z@\mkern-6mu\cleaders
  \hbox{$#1\mkern-2mu\box\z@\mkern-2mu$}\hfill
  \mkern-6mu\mathord\rightarrow$}
\def\leftarrowfill@#1{\setboxz@h{$#1-\m@th$}\ht\z@\z@
  $#1\m@th\mathord\leftarrow\mkern-6mu\cleaders
  \hbox{$#1\mkern-2mu\copy\z@\mkern-2mu$}\hfill
  \mkern-6mu\box\z@$}
\def\leftrightarrowfill@#1{\setboxz@h{$#1-\m@th$}\ht\z@\z@
  $#1\m@th\mathord\leftarrow\mkern-6mu\cleaders
  \hbox{$#1\mkern-2mu\box\z@\mkern-2mu$}\hfill
  \mkern-6mu\mathord\rightarrow$}
\def\overrightarrow{\mathpalette\overrightarrow@}
\def\overrightarrow@#1#2{\vbox{\ialign{##\crcr\rightarrowfill@#1\crcr
 \noalign{\kern-\ex@\nointerlineskip}$\m@th\hfil#1#2\hfil$\crcr}}}

\def\overleftarrow{\mathpalette\overleftarrow@}
\def\overleftarrow@#1#2{\vbox{\ialign{##\crcr\leftarrowfill@#1\crcr
 \noalign{\kern-\ex@\nointerlineskip}$\m@th\hfil#1#2\hfil$\crcr}}}
\def\overleftrightarrow{\mathpalette\overleftrightarrow@}
\def\overleftrightarrow@#1#2{\vbox{\ialign{##\crcr\leftrightarrowfill@#1\crcr
 \noalign{\kern-\ex@\nointerlineskip}$\m@th\hfil#1#2\hfil$\crcr}}}
\def\underrightarrow{\mathpalette\underrightarrow@}
\def\underrightarrow@#1#2{\vtop{\ialign{##\crcr$\m@th\hfil#1#2\hfil$\crcr
 \noalign{\nointerlineskip}\rightarrowfill@#1\crcr}}}

\def\underleftarrow{\mathpalette\underleftarrow@}
\def\underleftarrow@#1#2{\vtop{\ialign{##\crcr$\m@th\hfil#1#2\hfil$\crcr
 \noalign{\nointerlineskip}\leftarrowfill@#1\crcr}}}
\def\underleftrightarrow{\mathpalette\underleftrightarrow@}
\def\underleftrightarrow@#1#2{\vtop{\ialign{##\crcr$\m@th\hfil#1#2\hfil$\crcr
 \noalign{\nointerlineskip}\leftrightarrowfill@#1\crcr}}}
\message{various kinds of dots,}
\let\DOTSI\relax
\let\DOTSB\relax

\newif\ifmath@
{\uccode`7=`\\ \uccode`8=`m \uccode`9=`a \uccode`0=`t \uccode`!=`h
 \uppercase{\gdef\math@#1#2#3#4#5#6\math@{\global\math@false\ifx 7#1\ifx 8#2%
 \ifx 9#3\ifx 0#4\ifx !#5\xdef\meaning@{#6}\global\math@true\fi\fi\fi\fi\fi}}}
\newif\ifmathch@
{\uccode`7=`c \uccode`8=`h \uccode`9=`\"
 \uppercase{\gdef\mathch@#1#2#3#4#5#6\mathch@{\global\mathch@false
  \ifx 7#1\ifx 8#2\ifx 9#5\global\mathch@true\xdef\meaning@{9#6}\fi\fi\fi}}}
\newcount\classnum@
\def\getmathch@#1.#2\getmathch@{\classnum@#1 \divide\classnum@4096
 \ifcase\number\classnum@\or\or\gdef\thedots@{\dotsb@}\or
 \gdef\thedots@{\dotsb@}\fi}
\newif\ifmathbin@
{\uccode`4=`b \uccode`5=`i \uccode`6=`n
 \uppercase{\gdef\mathbin@#1#2#3{\relaxnext@
  \DNii@##1\mathbin@{\ifx\space@\next\global\mathbin@true\fi}%
 \global\mathbin@false\DN@##1\mathbin@{}%
 \ifx 4#1\ifx 5#2\ifx 6#3\DN@{\FN@\nextii@}\fi\fi\fi\next@}}}
\newif\ifmathrel@
{\uccode`4=`r \uccode`5=`e \uccode`6=`l
 \uppercase{\gdef\mathrel@#1#2#3{\relaxnext@
  \DNii@##1\mathrel@{\ifx\space@\next\global\mathrel@true\fi}%
 \global\mathrel@false\DN@##1\mathrel@{}%
 \ifx 4#1\ifx 5#2\ifx 6#3\DN@{\FN@\nextii@}\fi\fi\fi\next@}}}
\newif\ifmacro@
{\uccode`5=`m \uccode`6=`a \uccode`7=`c
 \uppercase{\gdef\macro@#1#2#3#4\macro@{\global\macro@false
  \ifx 5#1\ifx 6#2\ifx 7#3\global\macro@true
  \xdef\meaning@{\macro@@#4\macro@@}\fi\fi\fi}}}
\def\macro@@#1->#2\macro@@{#2}
\newif\ifDOTS@
\newcount\DOTSCASE@
{\uccode`6=`\\ \uccode`7=`D \uccode`8=`O \uccode`9=`T \uccode`0=`S
 \uppercase{\gdef\DOTS@#1#2#3#4#5{\global\DOTS@false\DN@##1\DOTS@{}%
  \ifx 6#1\ifx 7#2\ifx 8#3\ifx 9#4\ifx 0#5\let\next@\DOTS@@\fi\fi\fi\fi\fi
  \next@}}}
{\uccode`3=`B \uccode`4=`I \uccode`5=`X
 \uppercase{\gdef\DOTS@@#1{\relaxnext@
  \DNii@##1\DOTS@{\ifx\space@\next\global\DOTS@true\fi}%
  \DN@{\FN@\nextii@}%
  \ifx 3#1\global\DOTSCASE@\z@\else
  \ifx 4#1\global\DOTSCASE@\@ne\else
  \ifx 5#1\global\DOTSCASE@\tw@\else\DN@##1\DOTS@{}%
  \fi\fi\fi\next@}}}
\newif\ifnot@
{\uccode`5=`\\ \uccode`6=`n \uccode`7=`o \uccode`8=`t
 \uppercase{\gdef\not@#1#2#3#4{\relaxnext@
  \DNii@##1\not@{\ifx\space@\next\global\not@true\fi}%
 \global\not@false\DN@##1\not@{}%
 \ifx 5#1\ifx 6#2\ifx 7#3\ifx 8#4\DN@{\FN@\nextii@}\fi\fi\fi
 \fi\next@}}}
\newif\ifkeybin@
\def\keybin@{\keybin@true
 \ifx\next+\else\ifx\next=\else\ifx\next<\else\ifx\next>\else\ifx\next-\else
 \ifx\next*\else\ifx\next:\else\keybin@false\fi\fi\fi\fi\fi\fi\fi}
\def\dots{\RIfM@\expandafter\mdots@\else\expandafter\tdots@\fi}
\def\tdots@{\unskip\relaxnext@
 \DN@{$\m@th\mathinner{\ldotp\ldotp\ldotp}\,
   \ifx\next,\,$\else\ifx\next.\,$\else\ifx\next;\,$\else\ifx\next:\,$\else
   \ifx\next?\,$\else\ifx\next!\,$\else$ \fi\fi\fi\fi\fi\fi}%
 \ \FN@\next@}
\def\mdots@{\FN@\mdots@@}
\def\mdots@@{\gdef\thedots@{\dotso@}
 \ifx\next\boldkey\gdef\thedots@\boldkey{\boldkeydots@}\else                
 \ifx\next\boldsymbol\gdef\thedots@\boldsymbol{\boldsymboldots@}\else       
 \ifx,\next\gdef\thedots@{\dotsc}
 \else\ifx\not\next\gdef\thedots@{\dotsb@}
 \else\keybin@
 \ifkeybin@\gdef\thedots@{\dotsb@}
 \else\xdef\meaning@{\meaning\next..........}\xdef\meaning@@{\meaning@}
  \expandafter\math@\meaning@\math@
  \ifmath@
   \expandafter\mathch@\meaning@\mathch@
   \ifmathch@\expandafter\getmathch@\meaning@\getmathch@\fi                 
  \else\expandafter\macro@\meaning@@\macro@                                 
  \ifmacro@                                                                
   \expandafter\not@\meaning@\not@\ifnot@\gdef\thedots@{\dotsb@}
  \else\expandafter\DOTS@\meaning@\DOTS@
  \ifDOTS@
   \ifcase\number\DOTSCASE@\gdef\thedots@{\dotsb@}%
    \or\gdef\thedots@{\dotsi}\else\fi                                      
  \else\expandafter\math@\meaning@\math@                                   
  \ifmath@\expandafter\mathbin@\meaning@\mathbin@
  \ifmathbin@\gdef\thedots@{\dotsb@}
  \else\expandafter\mathrel@\meaning@\mathrel@
  \ifmathrel@\gdef\thedots@{\dotsb@}
  \fi\fi\fi\fi\fi\fi\fi\fi\fi\fi\fi\fi
 \thedots@}
\def\plainldots@{\mathinner{\ldotp\ldotp\ldotp}}
\def\plaincdots@{\mathinner{\cdotp\cdotp\cdotp}}
\def\dotsi{\!\plaincdots@}
\let\dotsb@\plaincdots@
\newif\ifextra@
\newif\ifrightdelim@
\def\rightdelim@{\global\rightdelim@true                                    
 \ifx\next)\else                                                            
 \ifx\next]\else
 \ifx\next\rbrack\else
 \ifx\next\}\else
 \ifx\next\rbrace\else
 \ifx\next\rangle\else
 \ifx\next\rceil\else
 \ifx\next\rfloor\else
 \ifx\next\rgroup\else
 \ifx\next\rmoustache\else
 \ifx\next\right\else
 \ifx\next\bigr\else
 \ifx\next\biggr\else
 \ifx\next\Bigr\else                                                        
 \ifx\next\Biggr\else\global\rightdelim@false
 \fi\fi\fi\fi\fi\fi\fi\fi\fi\fi\fi\fi\fi\fi\fi}
\def\extra@{%
 \global\extra@false\rightdelim@\ifrightdelim@\global\extra@true            
 \else\ifx\next$\global\extra@true                                          
 \else\xdef\meaning@{\meaning\next..........}
 \expandafter\macro@\meaning@\macro@\ifmacro@                               
 \expandafter\DOTS@\meaning@\DOTS@
 \ifDOTS@
 \ifnum\DOTSCASE@=\tw@\global\extra@true                                    
 \fi\fi\fi\fi\fi}
\newif\ifbold@
\def\dotso@{\relaxnext@
 \ifbold@
  \let\next\delayed@
  \DNii@{\extra@\plainldots@\ifextra@\,\fi}%
 \else
  \DNii@{\DN@{\extra@\plainldots@\ifextra@\,\fi}\FN@\next@}%
 \fi
 \nextii@}
\def\extrap@#1{%
 \ifx\next,\DN@{#1\,}\else
 \ifx\next;\DN@{#1\,}\else
 \ifx\next.\DN@{#1\,}\else\extra@
 \ifextra@\DN@{#1\,}\else
 \let\next@#1\fi\fi\fi\fi\next@}
\def\ldots{\DN@{\extrap@\plainldots@}%
 \FN@\next@}
\def\cdots{\DN@{\extrap@\plaincdots@}%
 \FN@\next@}

\def\dotsc{\relaxnext@
 \DN@{\ifx\next;\plainldots@\,\else
  \ifx\next.\plainldots@\,\else\extra@\plainldots@
  \ifextra@\,\fi\fi\fi}%
 \FN@\next@}
\def\cdot{\mathchar"2201 }

\def\mapsto{\DOTSB\mapstochar\rightarrow}

\message{special superscripts,}
\def\dddot#1{{\mathop{#1}\limits^{\vbox to-1.4\ex@{\kern-\tw@\ex@
 \hbox{\rm...}\vss}}}}
\def\ddddot#1{{\mathop{#1}\limits^{\vbox to-1.4\ex@{\kern-\tw@\ex@
 \hbox{\rm....}\vss}}}}
\def\sphat{^{\mathchoice{}{}%
 {\,\,\botsmash{\hbox{\lower4\ex@\hbox{$\m@th\widehat{\null}$}}}}%
 {\,\botsmash{\hbox{\lower3\ex@\hbox{$\m@th\hat{\null}$}}}}}}

\def\spacute{^{\!\botsmash{\hbox{\lower\@ne ex\hbox{\'{}}}}}}
\def\spgrave{^{\mathchoice{}{}{}{\!}%
 \botsmash{\hbox{\lower\@ne ex\hbox{\`{}}}}}}
\def\spdot{^{\hbox{\raise\ex@\hbox{\rm.}}}}
\def\spddot{^{\hbox{\raise\ex@\hbox{\rm..}}}}
\def\spdddot{^{\hbox{\raise\ex@\hbox{\rm...}}}}
\def\spddddot{^{\hbox{\raise\ex@\hbox{\rm....}}}}
\def\spbreve{^{\!\botsmash{\hbox{\lower4\ex@\hbox{\u{}}}}}}

\message{\string\text,}
\def\textonlyfont@#1#2{\def#1{\RIfM@
 \Err@{Use \string#1\space only in text}\else#2\fi}}
\textonlyfont@\rm\tenrm
\textonlyfont@\it\tenit
\textonlyfont@\sl\tensl
\textonlyfont@\bf\tenbf
\def\oldnos#1{\RIfM@{\mathcode`\,="013B \fam\@ne#1}\else
 \leavevmode\hbox{$\m@th\mathcode`\,="013B \fam\@ne#1$}\fi}
\def\text{\RIfM@\expandafter\text@\else\expandafter\text@@\fi}
\def\text@@#1{\leavevmode\hbox{#1}}
\def\mathhexbox@#1#2#3{\text{$\m@th\mathchar"#1#2#3$}}
\def\dag{{\mathhexbox@279}}
\def\ddag{{\mathhexbox@27A}}
\def\S{{\mathhexbox@278}}
\def\P{{\mathhexbox@27B}}
\newif\iffirstchoice@
\firstchoice@true
\def\text@#1{\mathchoice
 {\hbox{\everymath{\displaystyle}\def\textfonti{\the\textfont\@ne}%
  \def\textfontii{\the\textfont\tw@}\textdef@@ T#1}}
 {\hbox{\firstchoice@false
  \everymath{\textstyle}\def\textfonti{\the\textfont\@ne}%
  \def\textfontii{\the\textfont\tw@}\textdef@@ T#1}}
 {\hbox{\firstchoice@false
  \everymath{\scriptstyle}\def\textfonti{\the\scriptfont\@ne}%
  \def\textfontii{\the\scriptfont\tw@}\textdef@@ S\rm#1}}
 {\hbox{\firstchoice@false
  \everymath{\scriptscriptstyle}\def\textfonti
  {\the\scriptscriptfont\@ne}%
  \def\textfontii{\the\scriptscriptfont\tw@}\textdef@@ s\rm#1}}}
\def\textdef@@#1{\textdef@#1\rm\textdef@#1\bf\textdef@#1\sl\textdef@#1\it}
\def\rmfam{0}
\def\textdef@#1#2{%
 \DN@{\csname\expandafter\eat@\string#2fam\endcsname}%
 \if S#1\edef#2{\the\scriptfont\next@\relax}%
 \else\if s#1\edef#2{\the\scriptscriptfont\next@\relax}%
 \else\edef#2{\the\textfont\next@\relax}\fi\fi}
\scriptfont\itfam\tenit \scriptscriptfont\itfam\tenit
\scriptfont\slfam\tensl \scriptscriptfont\slfam\tensl
\newif\iftopfolded@
\newif\ifbotfolded@
\def\topfoldedtext{\topfolded@true\botfolded@false\foldedtext@}
\def\botfoldedtext{\botfolded@true\topfolded@false\foldedtext@}
\def\foldedtext{\topfolded@false\botfolded@false\foldedtext@}
\Invalid@\foldedwidth
\def\foldedtext@{\relaxnext@
 \DN@{\ifx\next\foldedwidth\let\next@\nextii@\else
  \DN@{\nextii@\foldedwidth{.3\hsize}}\fi\next@}%
 \DNii@\foldedwidth##1##2{\setbox\z@\vbox
  {\normalbaselines\hsize##1\relax
  \tolerance1600 \noindent\ignorespaces##2}\ifbotfolded@\boxz@\else
  \iftopfolded@\vtop{\unvbox\z@}\else\vcenter{\boxz@}\fi\fi}%
 \FN@\next@}
\message{math font commands,}
\def\bold{\RIfM@\expandafter\bold@\else
 \expandafter\nonmatherr@\expandafter\bold\fi}
\def\bold@#1{{\bold@@{#1}}}
\def\bold@@#1{\fam\bffam\relax#1}
\def\slanted{\RIfM@\expandafter\slanted@\else
 \expandafter\nonmatherr@\expandafter\slanted\fi}
\def\slanted@#1{{\slanted@@{#1}}}
\def\slanted@@#1{\fam\slfam\relax#1}
\def\roman{\RIfM@\expandafter\roman@\else
 \expandafter\nonmatherr@\expandafter\roman\fi}
\def\roman@#1{{\roman@@{#1}}}
\def\roman@@#1{\fam\rmfam\relax#1}
\def\italic{\RIfM@\expandafter\italic@\else
 \expandafter\nonmatherr@\expandafter\italic\fi}
\def\italic@#1{{\italic@@{#1}}}
\def\italic@@#1{\fam\itfam\relax#1}
\def\Cal{\RIfM@\expandafter\Cal@\else
 \expandafter\nonmatherr@\expandafter\Cal\fi}
\def\Cal@#1{{\Cal@@{#1}}}
\def\Cal@@#1{\noaccents@\fam\tw@#1}
\mathchardef\Gamma="0000
\mathchardef\Delta="0001
\mathchardef\Theta="0002
\mathchardef\Lambda="0003
\mathchardef\Xi="0004
\mathchardef\Pi="0005
\mathchardef\Sigma="0006
\mathchardef\Upsilon="0007
\mathchardef\Phi="0008
\mathchardef\Psi="0009
\mathchardef\Omega="000A
\mathchardef\varGamma="0100
\mathchardef\varDelta="0101
\mathchardef\varTheta="0102
\mathchardef\varLambda="0103
\mathchardef\varXi="0104
\mathchardef\varPi="0105
\mathchardef\varSigma="0106
\mathchardef\varUpsilon="0107
\mathchardef\varPhi="0108
\mathchardef\varPsi="0109
\mathchardef\varOmega="010A
\let\alloc@@\alloc@
\def\hexnumber@#1{\ifcase#1 0\or 1\or 2\or 3\or 4\or 5\or 6\or 7\or 8\or
 9\or A\or B\or C\or D\or E\or F\fi}
\def\loadmsam{%
 \font@\tenmsa=msam10
 \font@\sevenmsa=msam7
 \font@\fivemsa=msam5
 \alloc@@8\fam\chardef\sixt@@n\msafam
 \textfont\msafam=\tenmsa
 \scriptfont\msafam=\sevenmsa
 \scriptscriptfont\msafam=\fivemsa
 \edef\next{\hexnumber@\msafam}%
 \mathchardef\dabar@"0\next39
 \edef\dashrightarrow{\mathrel{\dabar@\dabar@\mathchar"0\next4B}}%
 \edef\dashleftarrow{\mathrel{\mathchar"0\next4C\dabar@\dabar@}}%
 \let\dasharrow\dashrightarrow
 \edef\ulcorner{\delimiter"4\next70\next70 }%
 \edef\urcorner{\delimiter"5\next71\next71 }%
 \edef\llcorner{\delimiter"4\next78\next78 }%
 \edef\lrcorner{\delimiter"5\next79\next79 }%
 \edef\yen{{\noexpand\mathhexbox@\next55}}%
 \edef\checkmark{{\noexpand\mathhexbox@\next58}}%
 \edef\circledR{{\noexpand\mathhexbox@\next72}}%
 \edef\maltese{{\noexpand\mathhexbox@\next7A}}%
 \global\let\loadmsam\empty}%
\def\loadmsbm{%
 \font@\tenmsb=msbm10 \font@\sevenmsb=msbm7 \font@\fivemsb=msbm5
 \alloc@@8\fam\chardef\sixt@@n\msbfam
 \textfont\msbfam=\tenmsb
 \scriptfont\msbfam=\sevenmsb \scriptscriptfont\msbfam=\fivemsb
 \global\let\loadmsbm\empty
 }
\def\widehat#1{\ifx\undefined\msbfam \DN@{362}%
  \else \setboxz@h{$\m@th#1$}%
    \edef\next@{\ifdim\wdz@>\tw@ em%
        \hexnumber@\msbfam 5B%
      \else 362\fi}\fi
  \mathaccent"0\next@{#1}}
\def\widetilde#1{\ifx\undefined\msbfam \DN@{365}%
  \else \setboxz@h{$\m@th#1$}%
    \edef\next@{\ifdim\wdz@>\tw@ em%
        \hexnumber@\msbfam 5D%
      \else 365\fi}\fi
  \mathaccent"0\next@{#1}}
\message{\string\newsymbol,}
\def\newsymbol#1#2#3#4#5{\define#1{}%
  \count@#2\relax \advance\count@\m@ne 
 \ifcase\count@
   \ifx\undefined\msafam\loadmsam\fi \let\next@\msafam
 \or \ifx\undefined\msbfam\loadmsbm\fi \let\next@\msbfam
 \else  \Err@{\Invalid@@\string\newsymbol}\let\next@\tw@\fi
 \mathchardef#1="#3\hexnumber@\next@#4#5\space}

\def\Bbb{\RIfM@\expandafter\Bbb@\else
 \expandafter\nonmatherr@\expandafter\Bbb\fi}
\def\Bbb@#1{{\Bbb@@{#1}}}
\def\Bbb@@#1{\noaccents@\fam\msbfam\relax#1}
\message{bold Greek and bold symbols,}
\def\loadbold{%
 \font@\tencmmib=cmmib10 \font@\sevencmmib=cmmib7 \font@\fivecmmib=cmmib5
 \skewchar\tencmmib'177 \skewchar\sevencmmib'177 \skewchar\fivecmmib'177
 \alloc@@8\fam\chardef\sixt@@n\cmmibfam
 \textfont\cmmibfam\tencmmib
 \scriptfont\cmmibfam\sevencmmib \scriptscriptfont\cmmibfam\fivecmmib
 \font@\tencmbsy=cmbsy10 \font@\sevencmbsy=cmbsy7 \font@\fivecmbsy=cmbsy5
 \skewchar\tencmbsy'60 \skewchar\sevencmbsy'60 \skewchar\fivecmbsy'60
 \alloc@@8\fam\chardef\sixt@@n\cmbsyfam
 \textfont\cmbsyfam\tencmbsy
 \scriptfont\cmbsyfam\sevencmbsy \scriptscriptfont\cmbsyfam\fivecmbsy
 \let\loadbold\empty
}
\def\boldnotloaded#1{\Err@{\ifcase#1\or First\else Second\fi
       bold symbol font not loaded}}
\def\mathchari@#1#2#3{\ifx\undefined\cmmibfam
    \boldnotloaded@\@ne
  \else\mathchar"#1\hexnumber@\cmmibfam#2#3\space \fi}
\def\mathcharii@#1#2#3{\ifx\undefined\cmbsyfam
    \boldnotloaded\tw@
  \else \mathchar"#1\hexnumber@\cmbsyfam#2#3\space\fi}
\edef\bffam@{\hexnumber@\bffam}
\def\boldkey#1{\ifcat\noexpand#1A%
  \ifx\undefined\cmmibfam \boldnotloaded\@ne
  \else {\fam\cmmibfam#1}\fi
 \else
 \ifx#1!\mathchar"5\bffam@21 \else
 \ifx#1(\mathchar"4\bffam@28 \else\ifx#1)\mathchar"5\bffam@29 \else
 \ifx#1+\mathchar"2\bffam@2B \else\ifx#1:\mathchar"3\bffam@3A \else
 \ifx#1;\mathchar"6\bffam@3B \else\ifx#1=\mathchar"3\bffam@3D \else
 \ifx#1?\mathchar"5\bffam@3F \else\ifx#1[\mathchar"4\bffam@5B \else
 \ifx#1]\mathchar"5\bffam@5D \else
 \ifx#1,\mathchari@63B \else
 \ifx#1-\mathcharii@200 \else
 \ifx#1.\mathchari@03A \else
 \ifx#1/\mathchari@03D \else
 \ifx#1<\mathchari@33C \else
 \ifx#1>\mathchari@33E \else
 \ifx#1*\mathcharii@203 \else
 \ifx#1|\mathcharii@06A \else
 \ifx#10\bold0\else\ifx#11\bold1\else\ifx#12\bold2\else\ifx#13\bold3\else
 \ifx#14\bold4\else\ifx#15\bold5\else\ifx#16\bold6\else\ifx#17\bold7\else
 \ifx#18\bold8\else\ifx#19\bold9\else
  \Err@{\string\boldkey\space can't be used with #1}%
 \fi\fi\fi\fi\fi\fi\fi\fi\fi\fi\fi\fi\fi\fi\fi
 \fi\fi\fi\fi\fi\fi\fi\fi\fi\fi\fi\fi\fi\fi}
\def\boldsymbol#1{%
 \DN@{\Err@{You can't use \string\boldsymbol\space with \string#1}#1}%
 \ifcat\noexpand#1A%
   \let\next@\relax
   \ifx\undefined\cmmibfam \boldnotloaded\@ne
   \else {\fam\cmmibfam#1}\fi
 \else
  \xdef\meaning@{\meaning#1.........}%
  \expandafter\math@\meaning@\math@
  \ifmath@
   \expandafter\mathch@\meaning@\mathch@
   \ifmathch@
    \expandafter\boldsymbol@@\meaning@\boldsymbol@@
   \fi
  \else
   \expandafter\macro@\meaning@\macro@
   \expandafter\delim@\meaning@\delim@
   \ifdelim@
    \expandafter\delim@@\meaning@\delim@@
   \else
    \boldsymbol@{#1}%
   \fi
  \fi
 \fi
 \next@}
\def\mathhexboxii@#1#2{\ifx\undefined\cmbsyfam
    \boldnotloaded\tw@
  \else \mathhexbox@{\hexnumber@\cmbsyfam}{#1}{#2}\fi}
\def\boldsymbol@#1{\let\next@\relax\let\next#1%
 \ifx\next\cdot\mathcharii@201 \else
 \ifx\next\prime{{\null\mathcharii@030 \null}}\else
 \ifx\next\lbrack\mathchar"4\bffam@5B \else
 \ifx\next\rbrack\mathchar"5\bffam@5D \else
 \ifx\next\{\mathcharii@466 \else
 \ifx\next\lbrace\mathcharii@466 \else
 \ifx\next\}\mathcharii@567 \else
 \ifx\next\rbrace\mathcharii@567 \else
 \ifx\next\surd{{\mathcharii@170}}\else
 \ifx\next\S{{\mathhexboxii@78}}\else
 \ifx\next\P{{\mathhexboxii@7B}}\else
 \ifx\next\dag{{\mathhexboxii@79}}\else
 \ifx\next\ddag{{\mathhexboxii@7A}}\else
 \DN@{\Err@{You can't use \string\boldsymbol\space with \string#1}#1}%
 \fi\fi\fi\fi\fi\fi\fi\fi\fi\fi\fi\fi\fi}
\def\boldsymbol@@#1.#2\boldsymbol@@{\classnum@#1 \count@@@\classnum@        
 \divide\classnum@4096 \count@\classnum@                                    
 \multiply\count@4096 \advance\count@@@-\count@ \count@@\count@@@           
 \divide\count@@@\@cclvi \count@\count@@                                    
 \multiply\count@@@\@cclvi \advance\count@@-\count@@@                       
 \divide\count@@@\@cclvi                                                    
 \multiply\classnum@4096 \advance\classnum@\count@@                         
 \ifnum\count@@@=\z@                                                        
  \count@"\bffam@ \multiply\count@\@cclvi
  \advance\classnum@\count@
  \DN@{\mathchar\number\classnum@}%
 \else
  \ifnum\count@@@=\@ne                                                      
   \ifx\undefined\cmmibfam \DN@{\boldnotloaded\@ne}%
   \else \count@\cmmibfam \multiply\count@\@cclvi
     \advance\classnum@\count@
     \DN@{\mathchar\number\classnum@}\fi
  \else
   \ifnum\count@@@=\tw@                                                    
     \ifx\undefined\cmbsyfam
       \DN@{\boldnotloaded\tw@}%
     \else
       \count@\cmbsyfam \multiply\count@\@cclvi
       \advance\classnum@\count@
       \DN@{\mathchar\number\classnum@}%
     \fi
  \fi
 \fi
\fi}
\newif\ifdelim@
\newcount\delimcount@
{\uccode`6=`\\ \uccode`7=`d \uccode`8=`e \uccode`9=`l
 \uppercase{\gdef\delim@#1#2#3#4#5\delim@
  {\delim@false\ifx 6#1\ifx 7#2\ifx 8#3\ifx 9#4\delim@true
   \xdef\meaning@{#5}\fi\fi\fi\fi}}}
\def\delim@@#1"#2#3#4#5#6\delim@@{\if#32%
\let\next@\relax
 \ifx\undefined\cmbsyfam \boldnotloaded\@ne
 \else \mathcharii@#2#4#5\space \fi\fi}
\def\vert{\delimiter"026A30C }
\def\Vert{\delimiter"026B30D }
\let\|\Vert
\def\backslash{\delimiter"026E30F }
\def\boldkeydots@#1{\bold@true\let\next=#1\let\delayed@=#1\mdots@@
 \boldkey#1\bold@false}  
\def\boldsymboldots@#1{\bold@true\let\next#1\let\delayed@#1\mdots@@
 \boldsymbol#1\bold@false}
\message{Euler fonts,}

\def\frak{\mathfont@\frak}

\def\loadmathfont#1{%
   \expandafter\font@\csname ten#1\endcsname=#110
   \expandafter\font@\csname seven#1\endcsname=#17
   \expandafter\font@\csname five#1\endcsname=#15
   \edef\next{\noexpand\alloc@@8\fam\chardef\sixt@@n
     \expandafter\noexpand\csname#1fam\endcsname}%
   \next
   \textfont\csname#1fam\endcsname \csname ten#1\endcsname
   \scriptfont\csname#1fam\endcsname \csname seven#1\endcsname
   \scriptscriptfont\csname#1fam\endcsname \csname five#1\endcsname
   \expandafter\def\csname #1\expandafter\endcsname\expandafter{%
      \expandafter\mathfont@\csname#1\endcsname}%
 \expandafter\gdef\csname load#1\endcsname{}%
}
\def\mathfont@#1{\RIfM@\expandafter\mathfont@@\expandafter#1\else
  \expandafter\nonmatherr@\expandafter#1\fi}
\def\mathfont@@#1#2{{\mathfont@@@#1{#2}}}
\def\mathfont@@@#1#2{\noaccents@
   \fam\csname\expandafter\eat@\string#1fam\endcsname
   \relax#2}
\message{math accents,}
\def\accentclass@{7}
\def\noaccents@{\def\accentclass@{0}}
\def\makeacc@#1#2{\def#1{\mathaccent"\accentclass@#2 }}
\makeacc@\hat{05E}
\makeacc@\check{014}
\makeacc@\tilde{07E}
\makeacc@\acute{013}
\makeacc@\grave{012}
\makeacc@\dot{05F}
\makeacc@\ddot{07F}
\makeacc@\breve{015}
\makeacc@\bar{016}

\newcount\skewcharcount@
\newcount\familycount@
\def\theskewchar@{\familycount@\@ne
 \global\skewcharcount@\the\skewchar\textfont\@ne                           
 \ifnum\fam>\m@ne\ifnum\fam<16
  \global\familycount@\the\fam\relax
  \global\skewcharcount@\the\skewchar\textfont\the\fam\relax\fi\fi          
 \ifnum\skewcharcount@>\m@ne
  \ifnum\skewcharcount@<128
  \multiply\familycount@256
  \global\advance\skewcharcount@\familycount@
  \global\advance\skewcharcount@28672
  \mathchar\skewcharcount@\else
  \global\skewcharcount@\m@ne\fi\else
 \global\skewcharcount@\m@ne\fi}                                            
\newcount\pointcount@
\def\getpoints@#1.#2\getpoints@{\pointcount@#1 }
\newdimen\accentdimen@
\newcount\accentmu@
\def\dimentomu@{\multiply\accentdimen@ 100
 \expandafter\getpoints@\the\accentdimen@\getpoints@
 \multiply\pointcount@18
 \divide\pointcount@\@m
 \global\accentmu@\pointcount@}
\def\Makeacc@#1#2{\def#1{\RIfM@\DN@{\mathaccent@
 {"\accentclass@#2 }}\else\DN@{\nonmatherr@{#1}}\fi\next@}}
\def\unbracefonts@{\let\Cal@\Cal@@\let\roman@\roman@@\let\bold@\bold@@
 \let\slanted@\slanted@@}
\def\mathaccent@#1#2{\ifnum\fam=\m@ne\xdef\thefam@{1}\else
 \xdef\thefam@{\the\fam}\fi                                                 
 \accentdimen@\z@                                                           
 \setboxz@h{\unbracefonts@$\m@th\fam\thefam@\relax#2$}
 \ifdim\accentdimen@=\z@\DN@{\mathaccent#1{#2}}
  \setbox@ne\hbox{\unbracefonts@$\m@th\fam\thefam@\relax#2\theskewchar@$}
  \setbox\tw@\hbox{$\m@th\ifnum\skewcharcount@=\m@ne\else
   \mathchar\skewcharcount@\fi$}
  \global\accentdimen@\wd@ne\global\advance\accentdimen@-\wdz@
  \global\advance\accentdimen@-\wd\tw@                                     
  \global\multiply\accentdimen@\tw@
  \dimentomu@\global\advance\accentmu@\@ne                                 
 \else\DN@{{\mathaccent#1{#2\mkern\accentmu@ mu}%
    \mkern-\accentmu@ mu}{}}\fi                                             
 \next@}\Makeacc@\Hat{05E}
\Makeacc@\Check{014}
\Makeacc@\Tilde{07E}
\Makeacc@\Acute{013}
\Makeacc@\Grave{012}
\Makeacc@\Dot{05F}
\Makeacc@\Ddot{07F}
\Makeacc@\Breve{015}
\Makeacc@\Bar{016}
\def\Vec{\RIfM@\DN@{\mathaccent@{"017E }}\else
 \DN@{\nonmatherr@\Vec}\fi\next@}
\def\accentedsymbol#1#2{\csname newbox\expandafter\endcsname
  \csname\expandafter\eat@\string#1@box\endcsname
 \expandafter\setbox\csname\expandafter\eat@
  \string#1@box\endcsname\hbox{$\m@th#2$}\define
  #1{\copy\csname\expandafter\eat@\string#1@box\endcsname{}}}
\message{roots,}
\def\sqrt#1{\radical"270370 {#1}}
\let\underline@\underline
\let\overline@\overline
\def\underline#1{\underline@{#1}}
\def\overline#1{\overline@{#1}}
\Invalid@\leftroot
\Invalid@\uproot
\newcount\uproot@
\newcount\leftroot@
\def\root{\relaxnext@
  \DN@{\ifx\next\uproot\let\next@\nextii@\else
   \ifx\next\leftroot\let\next@\nextiii@\else
   \let\next@\plainroot@\fi\fi\next@}%
  \DNii@\uproot##1{\uproot@##1\relax\FN@\nextiv@}%
  \def\nextiv@{\ifx\next\space@\DN@. {\FN@\nextv@}\else
   \DN@.{\FN@\nextv@}\fi\next@.}%
  \def\nextv@{\ifx\next\leftroot\let\next@\nextvi@\else
   \let\next@\plainroot@\fi\next@}%
  \def\nextvi@\leftroot##1{\leftroot@##1\relax\plainroot@}%
   \def\nextiii@\leftroot##1{\leftroot@##1\relax\FN@\nextvii@}%
  \def\nextvii@{\ifx\next\space@
   \DN@. {\FN@\nextviii@}\else
   \DN@.{\FN@\nextviii@}\fi\next@.}%
  \def\nextviii@{\ifx\next\uproot\let\next@\nextix@\else
   \let\next@\plainroot@\fi\next@}%
  \def\nextix@\uproot##1{\uproot@##1\relax\plainroot@}%
  \bgroup\uproot@\z@\leftroot@\z@\FN@\next@}
\def\plainroot@#1\of#2{\setbox\rootbox\hbox{$\m@th\scriptscriptstyle{#1}$}%
 \mathchoice{\r@@t\displaystyle{#2}}{\r@@t\textstyle{#2}}
 {\r@@t\scriptstyle{#2}}{\r@@t\scriptscriptstyle{#2}}\egroup}
\def\r@@t#1#2{\setboxz@h{$\m@th#1\sqrt{#2}$}%
 \dimen@\ht\z@\advance\dimen@-\dp\z@
 \setbox@ne\hbox{$\m@th#1\mskip\uproot@ mu$}\advance\dimen@ 1.667\wd@ne
 \mkern-\leftroot@ mu\mkern5mu\raise.6\dimen@\copy\rootbox
 \mkern-10mu\mkern\leftroot@ mu\boxz@}
\def\boxed#1{\setboxz@h{$\m@th\displaystyle{#1}$}\dimen@.4\ex@
 \advance\dimen@3\ex@\advance\dimen@\dp\z@
 \hbox{\lower\dimen@\hbox{%
 \vbox{\hrule height.4\ex@
 \hbox{\vrule width.4\ex@\hskip3\ex@\vbox{\vskip3\ex@\boxz@\vskip3\ex@}%
 \hskip3\ex@\vrule width.4\ex@}\hrule height.4\ex@}%
 }}}
\message{commutative diagrams,}
\let\ampersand@\relax
\newdimen\minaw@
\minaw@11.11128\ex@
\newdimen\minCDaw@
\minCDaw@2.5pc
\def\minCDarrowwidth#1{\RIfMIfI@\onlydmatherr@\minCDarrowwidth
 \else\minCDaw@#1\relax\fi\else\onlydmatherr@\minCDarrowwidth\fi}
\newif\ifCD@
\def\CD{\bgroup\vspace@\relax\let\ampersand@&\iffalse}\fi
 \CD@true\vcenter\bgroup\Let@\tabskip\z@skip\baselineskip20\ex@
 \lineskip3\ex@\lineskiplimit3\ex@\halign\bgroup
 &\hfill$\m@th##$\hfill\crcr}
\def\endCD{\crcr\egroup\egroup\egroup}
\newdimen\bigaw@
\atdef@>#1>#2>{\ampersand@                                                  
 \setboxz@h{$\m@th\ssize\;{#1}\;\;$}
 \setbox@ne\hbox{$\m@th\ssize\;{#2}\;\;$}
 \setbox\tw@\hbox{$\m@th#2$}
 \ifCD@\global\bigaw@\minCDaw@\else\global\bigaw@\minaw@\fi                 
 \ifdim\wdz@>\bigaw@\global\bigaw@\wdz@\fi
 \ifdim\wd@ne>\bigaw@\global\bigaw@\wd@ne\fi                                
 \ifCD@\enskip\fi                                                           
 \ifdim\wd\tw@>\z@
  \mathrel{\mathop{\hbox to\bigaw@{\rightarrowfill@\displaystyle}}%
    \limits^{#1}_{#2}}
 \else\mathrel{\mathop{\hbox to\bigaw@{\rightarrowfill@\displaystyle}}%
    \limits^{#1}}\fi                                                        
 \ifCD@\enskip\fi                                                          
 \ampersand@}                                                              
\atdef@<#1<#2<{\ampersand@\setboxz@h{$\m@th\ssize\;\;{#1}\;$}%
 \setbox@ne\hbox{$\m@th\ssize\;\;{#2}\;$}\setbox\tw@\hbox{$\m@th#2$}%
 \ifCD@\global\bigaw@\minCDaw@\else\global\bigaw@\minaw@\fi
 \ifdim\wdz@>\bigaw@\global\bigaw@\wdz@\fi
 \ifdim\wd@ne>\bigaw@\global\bigaw@\wd@ne\fi
 \ifCD@\enskip\fi
 \ifdim\wd\tw@>\z@
  \mathrel{\mathop{\hbox to\bigaw@{\leftarrowfill@\displaystyle}}%
       \limits^{#1}_{#2}}\else
  \mathrel{\mathop{\hbox to\bigaw@{\leftarrowfill@\displaystyle}}%
       \limits^{#1}}\fi
 \ifCD@\enskip\fi\ampersand@}
\begingroup
 \catcode`\~=\active \lccode`\~=`\@
 \lowercase{%
  \global\atdef@)#1)#2){~>#1>#2>}
  \global\atdef@(#1(#2({~<#1<#2<}}
\endgroup
\atdef@ A#1A#2A{\llap{$\m@th\vcenter{\hbox
 {$\ssize#1$}}$}\Big\uparrow\rlap{$\m@th\vcenter{\hbox{$\ssize#2$}}$}&&}
\atdef@ V#1V#2V{\llap{$\m@th\vcenter{\hbox
 {$\ssize#1$}}$}\Big\downarrow\rlap{$\m@th\vcenter{\hbox{$\ssize#2$}}$}&&}
\atdef@={&\enskip\mathrel
 {\vbox{\hrule width\minCDaw@\vskip3\ex@\hrule width
 \minCDaw@}}\enskip&}
\atdef@|{\Big\Vert&&}
\atdef@\vert{\Big\Vert&&}
\def\pretend#1\haswidth#2{\setboxz@h{$\m@th\scriptstyle{#2}$}\hbox
 to\wdz@{\hfill$\m@th\scriptstyle{#1}$\hfill}}
\message{poor man's bold,}
\def\pmb{\RIfM@\expandafter\mathpalette\expandafter\pmb@\else
 \expandafter\pmb@@\fi}
\def\pmb@@#1{\leavevmode\setboxz@h{#1}%
   \dimen@-\wdz@
   \kern-.5\ex@\copy\z@
   \kern\dimen@\kern.25\ex@\raise.4\ex@\copy\z@
   \kern\dimen@\kern.25\ex@\box\z@
}
\def\binrel@@#1{\ifdim\wd2<\z@\mathbin{#1}\else\ifdim\wd\tw@>\z@
 \mathrel{#1}\else{#1}\fi\fi}
\newdimen\pmbraise@
\def\pmb@#1#2{\setbox\thr@@\hbox{$\m@th#1{#2}$}%
 \setbox4\hbox{$\m@th#1\mkern.5mu$}\pmbraise@\wd4\relax
 \binrel@{#2}%
 \dimen@-\wd\thr@@
   \binrel@@{%
   \mkern-.8mu\copy\thr@@
   \kern\dimen@\mkern.4mu\raise\pmbraise@\copy\thr@@
   \kern\dimen@\mkern.4mu\box\thr@@
}}
\def\documentstyle#1{\W@{}\input #1.sty\relax}
\message{syntax check,}
\font\dummyft@=dummy
\fontdimen1 \dummyft@=\z@
\fontdimen2 \dummyft@=\z@
\fontdimen3 \dummyft@=\z@
\fontdimen4 \dummyft@=\z@
\fontdimen5 \dummyft@=\z@
\fontdimen6 \dummyft@=\z@
\fontdimen7 \dummyft@=\z@
\fontdimen8 \dummyft@=\z@
\fontdimen9 \dummyft@=\z@
\fontdimen10 \dummyft@=\z@
\fontdimen11 \dummyft@=\z@
\fontdimen12 \dummyft@=\z@
\fontdimen13 \dummyft@=\z@
\fontdimen14 \dummyft@=\z@
\fontdimen15 \dummyft@=\z@
\fontdimen16 \dummyft@=\z@
\fontdimen17 \dummyft@=\z@
\fontdimen18 \dummyft@=\z@
\fontdimen19 \dummyft@=\z@
\fontdimen20 \dummyft@=\z@
\fontdimen21 \dummyft@=\z@
\fontdimen22 \dummyft@=\z@
\def\fontlist@{\\{\tenrm}\\{\sevenrm}\\{\fiverm}\\{\teni}\\{\seveni}%
 \\{\fivei}\\{\tensy}\\{\sevensy}\\{\fivesy}\\{\tenex}\\{\tenbf}\\{\sevenbf}%
 \\{\fivebf}\\{\tensl}\\{\tenit}}
\def\font@#1=#2 {\rightappend@#1\to\fontlist@\font#1=#2 }
\def\dodummy@{{\def\\##1{\global\let##1\dummyft@}\fontlist@}}
\def\nopages@{\output{\setbox\z@\box\@cclv \deadcycles\z@}%
 \alloc@5\toks\toksdef\@cclvi\output}
\let\galleys\nopages@
\newif\ifsyntax@
\newcount\countxviii@
\def\syntax{\syntax@true\dodummy@\countxviii@\count18
 \loop\ifnum\countxviii@>\m@ne\textfont\countxviii@=\dummyft@
 \scriptfont\countxviii@=\dummyft@\scriptscriptfont\countxviii@=\dummyft@
 \advance\countxviii@\m@ne\repeat                                           
 \dummyft@\tracinglostchars\z@\nopages@\frenchspacing\hbadness\@M}
\def\first@#1#2\end{#1}
\def\printoptions{\W@{Do you want S(yntax check),
  G(alleys) or P(ages)?}%
 \message{Type S, G or P, followed by <return>: }%
 \begingroup 
 \endlinechar\m@ne 
 \read\m@ne to\ans@
 \edef\ans@{\uppercase{\def\noexpand\ans@{%
   \expandafter\first@\ans@ P\end}}}%
 \expandafter\endgroup\ans@
 \if\ans@ P
 \else \if\ans@ S\syntax
 \else \if\ans@ G\galleys
 \else\message{? Unknown option: \ans@; using the `pages' option.}%
 \fi\fi\fi}
\def\alloc@#1#2#3#4#5{\global\advance\count1#1by\@ne
 \ch@ck#1#4#2\allocationnumber=\count1#1
 \global#3#5=\allocationnumber
 \ifalloc@\wlog{\string#5=\string#2\the\allocationnumber}\fi}
\def\document{\def\alloclist@{}\def\fontlist@{}}
\let\enddocument\bye

\let\proclaim\undefined
\let\footnote\undefined
\let\=\undefined
\let\>\undefined

\catcode`\@=\active
\message{... finished}

\expandafter\ifx\csname mathdefs.tex\endcsname\relax
  \expandafter\gdef\csname mathdefs.tex\endcsname{}
\else \message{Hey!  Apparently you were trying to
  \string\input{mathdefs.tex} twice.   This does not make sense.} 
\errmessage{Please edit your file (probably \jobname.tex) and remove
any duplicate ``\string\input'' lines}\endinput\fi




\catcode`\X=12\catcode`\@=11

\def\n@wcount{\alloc@0\count\countdef\insc@unt}
\def\n@wwrite{\alloc@7\write\chardef\sixt@@n}
\def\n@wread{\alloc@6\read\chardef\sixt@@n}
\def\r@s@t{\relax}\def\v@idline{\par}\def\@mputate#1/{#1}
\def\l@c@l#1X{\firstpart.#1}\def\gl@b@l#1X{#1}\def\t@d@l#1X{{}}

\def\crossrefs#1{\ifx\all#1\let\tr@ce=\all\else\def\tr@ce{#1,}\fi
   \n@wwrite\cit@tionsout\openout\cit@tionsout=\jobname.cit 
   \write\cit@tionsout{\tr@ce}\expandafter\setfl@gs\tr@ce,}
\def\setfl@gs#1,{\def\@{#1}\ifx\@\empty\let\next=\relax
   \else\let\next=\setfl@gs\expandafter\xdef
   \csname#1tr@cetrue\endcsname{}\fi\next}
\def\m@ketag#1#2{\expandafter\n@wcount\csname#2tagno\endcsname
     \csname#2tagno\endcsname=0\let\tail=\all\xdef\all{\tail#2,}
   \ifx#1\l@c@l\let\tail=\r@s@t\xdef\r@s@t{\csname#2tagno\endcsname=0\tail}\fi
   \expandafter\gdef\csname#2cite\endcsname##1{\expandafter
     \ifx\csname#2tag##1\endcsname\relax?\else\csname#2tag##1\endcsname\fi
     \expandafter\ifx\csname#2tr@cetrue\endcsname\relax\else
     \write\cit@tionsout{#2tag ##1 cited on page \folio.}\fi}
   \expandafter\gdef\csname#2page\endcsname##1{\expandafter
     \ifx\csname#2page##1\endcsname\relax?\else\csname#2page##1\endcsname\fi
     \expandafter\ifx\csname#2tr@cetrue\endcsname\relax\else
     \write\cit@tionsout{#2tag ##1 cited on page \folio.}\fi}
   \expandafter\gdef\csname#2tag\endcsname##1{\expandafter
      \ifx\csname#2check##1\endcsname\relax
      \expandafter\xdef\csname#2check##1\endcsname{}%
      \else\immediate\write16{Warning: #2tag ##1 used more than once.}\fi
      \multit@g{#1}{#2}##1/X%
      \write\t@gsout{#2tag ##1 assigned number \csname#2tag##1\endcsname\space
      on page \number\count0.}%
   \csname#2tag##1\endcsname}}

\def\multit@g#1#2#3/#4X{\def\t@mp{#4}\ifx\t@mp\empty%
      \global\advance\csname#2tagno\endcsname by 1 
      \expandafter\xdef\csname#2tag#3\endcsname
      {#1\number\csname#2tagno\endcsnameX}%
   \else\expandafter\ifx\csname#2last#3\endcsname\relax
      \expandafter\n@wcount\csname#2last#3\endcsname
      \global\advance\csname#2tagno\endcsname by 1 
      \expandafter\xdef\csname#2tag#3\endcsname
      {#1\number\csname#2tagno\endcsnameX}
      \write\t@gsout{#2tag #3 assigned number \csname#2tag#3\endcsname\space
      on page \number\count0.}\fi
   \global\advance\csname#2last#3\endcsname by 1
   \def\t@mp{\expandafter\xdef\csname#2tag#3/}%
   \expandafter\t@mp\@mputate#4\endcsname
   {\csname#2tag#3\endcsname\lastpart{\csname#2last#3\endcsname}}\fi}
\def\t@gs#1{\def\all{}\m@ketag#1e\m@ketag#1s\m@ketag\t@d@l p
\let\realscite\scite
\let\realstag\stag
   \m@ketag\gl@b@l r \n@wread\t@gsin
   \openin\t@gsin=\jobname.tgs \re@der \closein\t@gsin
   \n@wwrite\t@gsout\openout\t@gsout=\jobname.tgs }
\outer\def\localtags{\t@gs\l@c@l}
\outer\def\globaltags{\t@gs\gl@b@l}
\outer\def\newlocaltag#1{\m@ketag\l@c@l{#1}}
\outer\def\newglobaltag#1{\m@ketag\gl@b@l{#1}}

\newif\ifpr@ 
\def\m@kecs #1tag #2 assigned number #3 on page #4.%
   {\expandafter\gdef\csname#1tag#2\endcsname{#3}
   \expandafter\gdef\csname#1page#2\endcsname{#4}
   \ifpr@\expandafter\xdef\csname#1check#2\endcsname{}\fi}
\def\re@der{\ifeof\t@gsin\let\next=\relax\else
   \read\t@gsin to\t@gline\ifx\t@gline\v@idline\else
   \expandafter\m@kecs \t@gline\fi\let \next=\re@der\fi\next}
\def\pretags#1{\pr@true\pret@gs#1,,}
\def\pret@gs#1,{\def\@{#1}\ifx\@\empty\let\n@xtfile=\relax
   \else\let\n@xtfile=\pret@gs \openin\t@gsin=#1.tgs \message{#1} \re@der 
   \closein\t@gsin\fi \n@xtfile}

\newcount\sectno\sectno=0\newcount\subsectno\subsectno=0
\newif\ifultr@local \def\ultralocal{\ultr@localtrue}
\def\firstpart{\number\sectno}
\def\lastpart#1{\ifcase#1 \or a\or b\or c\or d\or e\or f\or g\or h\or 
   i\or k\or l\or m\or n\or o\or p\or q\or r\or s\or t\or u\or v\or w\or 
   x\or y\or z \fi}

\def\resetall{\global\advance\sectno by 1\subsectno=0
   \gdef\firstpart{\number\sectno}\r@s@t}
\def\resetsub{\global\advance\subsectno by 1
   \gdef\firstpart{\number\sectno.\number\subsectno}\r@s@t}
\def\newsection#1\par{\resetall\vskip0pt plus.3\vsize\penalty-250
   \vskip0pt plus-.3\vsize\bigskip\bigskip
   \message{#1}\leftline{\bf#1}\nobreak\bigskip}
\def\subsection#1\par{\ifultr@local\resetsub\fi
   \vskip0pt plus.2\vsize\penalty-250\vskip0pt plus-.2\vsize
   \bigskip\smallskip\message{#1}\leftline{\bf#1}\nobreak\medskip}


\newdimen\marginshift

\newdimen\margindelta
\newdimen\marginmax
\newdimen\marginmin

\def\margininit{       
\marginmax=3 true cm                  
				      
\margindelta=0.1 true cm              
\marginmin=0.1true cm                 
\marginshift=\marginmin
}    

\def\t@gsjj#1,{\def\@{#1}\ifx\@\empty\let\next=\relax\else\let\next=\t@gsjj
   \def\@@{p}\ifx\@\@@\else
   \expandafter\gdef\csname#1cite\endcsname##1{\citejj{##1}}
   \expandafter\gdef\csname#1page\endcsname##1{?}
   \expandafter\gdef\csname#1tag\endcsname##1{\tagjj{##1}}\fi\fi\next}
\newif\ifshowstuffinmargin
\showstuffinmarginfalse
\def\jjtags{\ifx\shlhetal\relax 
  \else
\ifx\shlhetal\undefinedcontrolseq
\else
\showstuffinmargintrue
\ifx\all\relax\else\expandafter\t@gsjj\all,\fi\fi \fi
}

\def\tagjj#1{\realstag{#1}\oldmginpar{\zeigen{#1}}}
\def\citejj#1{\rechnen{#1}\mginpar{\zeigen{#1}}}     

\def\rechnen#1{\expandafter\ifx\csname stag#1\endcsname\relax ??\else
                           \csname stag#1\endcsname\fi}

\newdimen\theight

\def\marginfont{\sevenrm}

\def\trymarginbox#1{\setbox0=\hbox{\marginfont\hskip\marginshift #1}%
		\global\marginshift\wd0 
		\global\advance\marginshift\margindelta}

\def \oldmginpar#1{%
\ifvmode\setbox0\hbox to \hsize{\hfill\rlap{\marginfont\quad#1}}%
\ht0 0cm
\dp0 0cm
\box0\vskip-\baselineskip
\else 
             \vadjust{\trymarginbox{#1}%
		\ifdim\marginshift>\marginmax \global\marginshift\marginmin
			\trymarginbox{#1}%
                \fi
             \theight=\ht0
             \advance\theight by \dp0    \advance\theight by \lineskip
             \kern -\theight \vbox to \theight{\rightline{\rlap{\box0}}%
\vss}}\fi}

\newdimen\upordown
\global\upordown=8pt
\font\tinyfont=cmtt8 
\def\mginpar#1{\smash{\hbox to 0cm{\kern-10pt\raise7pt\hbox{\tinyfont #1}\hss}}}
\def\mginpar#1{{\hbox to 0cm{\kern-10pt\raise\upordown\hbox{\tinyfont #1}\hss}}\global\upordown-\upordown}


\def\t@gsoff#1,{\def\@{#1}\ifx\@\empty\let\next=\relax\else\let\next=\t@gsoff
   \def\@@{p}\ifx\@\@@\else
   \expandafter\gdef\csname#1cite\endcsname##1{\zeigen{##1}}
   \expandafter\gdef\csname#1page\endcsname##1{?}
   \expandafter\gdef\csname#1tag\endcsname##1{\zeigen{##1}}\fi\fi\next}
\def\verbatimtags{\showstuffinmarginfalse
\ifx\all\relax\else\expandafter\t@gsoff\all,\fi}
\def\zeigen#1{\hbox{$\scriptstyle\langle$}#1\hbox{$\scriptstyle\rangle$}}


\def\margintag#1{\ifshowstuffinmargin\oldmginpar{\zeigen{#1}}\fi}

\def\(#1){\edef\dot@g{\ifmmode\ifinner(\hbox{\noexpand\etag{#1}})
   \else\noexpand\eqno(\hbox{\noexpand\etag{#1}})\fi
   \else(\noexpand\ecite{#1})\fi}\dot@g}

\newif\ifbr@ck
\def\eat#1{}
\def\[#1]{\br@cktrue[\br@cket#1'X]}
\def\br@cket#1'#2X{\def\temp{#2}\ifx\temp\empty\let\next\eat
   \else\let\next\br@cket\fi
   \ifbr@ck\br@ckfalse\br@ck@t#1,X\else\br@cktrue#1\fi\next#2X}
\def\br@ck@t#1,#2X{\def\temp{#2}\ifx\temp\empty\let\neext\eat
   \else\let\neext\br@ck@t\def\temp{,}\fi
   \def\teemp{#1}\ifx\teemp\empty\else\rcite{#1}\fi\temp\neext#2X}
\def\resetbr@cket{\gdef\[##1]{[\rtag{##1}]}}
\def\references{\resetbr@cket\newsection References\par}

\newtoks\symb@ls\newtoks\s@mb@ls\newtoks\p@gelist\n@wcount\ftn@mber
    \ftn@mber=1\newif\ifftn@mbers\ftn@mbersfalse\newif\ifbyp@ge\byp@gefalse
\def\defm@rk{\ifftn@mbers\n@mberm@rk\else\symb@lm@rk\fi}
\def\n@mberm@rk{\xdef\m@rk{{\the\ftn@mber}}%
    \global\advance\ftn@mber by 1 }
\def\rot@te#1{\let\temp=#1\global#1=\expandafter\r@t@te\the\temp,X}
\def\r@t@te#1,#2X{{#2#1}\xdef\m@rk{{#1}}}
\def\b@@st#1{{$^{#1}$}}\def\str@p#1{#1}
\def\symb@lm@rk{\ifbyp@ge\rot@te\p@gelist\ifnum\expandafter\str@p\m@rk=1 
    \s@mb@ls=\symb@ls\fi\write\f@nsout{\number\count0}\fi \rot@te\s@mb@ls}
\def\byp@ge{\byp@getrue\n@wwrite\f@nsin\openin\f@nsin=\jobname.fns 
    \n@wcount\currentp@ge\currentp@ge=0\p@gelist={0}
    \re@dfns\closein\f@nsin\rot@te\p@gelist
    \n@wread\f@nsout\openout\f@nsout=\jobname.fns }
\def\m@kelist#1X#2{{#1,#2}}
\def\re@dfns{\ifeof\f@nsin\let\next=\relax\else\read\f@nsin to \f@nline
    \ifx\f@nline\v@idline\else\let\t@mplist=\p@gelist
    \ifnum\currentp@ge=\f@nline
    \global\p@gelist=\expandafter\m@kelist\the\t@mplistX0
    \else\currentp@ge=\f@nline
    \global\p@gelist=\expandafter\m@kelist\the\t@mplistX1\fi\fi
    \let\next=\re@dfns\fi\next}
\def\symbols#1{\symb@ls={#1}\s@mb@ls=\symb@ls} 
\def\bigsymbol{\textstyle}
\symbols{\bigsymbol\ast,\dagger,\ddagger,\sharp,\flat,\natural,\star}
\def\ftnumbers{\ftn@mberstrue} \def\ftsymbols{\ftn@mbersfalse}
\def\paginal{\byp@ge} \def\resetftnumbers{\ftn@mber=1}
\def\ftnote#1{\defm@rk\expandafter\expandafter\expandafter\footnote
    \expandafter\b@@st\m@rk{#1}}

\long\def\jump#1\endjump{}
\def\ssum{\mathop{\lower .1em\hbox{$\textstyle\Sigma$}}\nolimits}

\def\qed{\nobreak\kern 1em \vrule height .5em width .5em depth 0em}
\def\newneq{\hbox{\rlap{\hbox to 1\wd9{\hss$=$\hss}}\raise .1em 
   \hbox to 1\wd9{\hss$\scriptscriptstyle/$\hss}}}
\def\subsetne{\setbox9 = \hbox{$\subset$}\mathrel{\hbox{\rlap
   {\lower .4em \newneq}\raise .13em \hbox{$\subset$}}}}
\def\supsetne{\setbox9 = \hbox{$\subset$}\mathrel{\hbox{\rlap
   {\lower .4em \newneq}\raise .13em \hbox{$\supset$}}}}

\def\vbar{\mathchoice{\vrule height6.3ptdepth-.5ptwidth.8pt\kern-.8pt}
   {\vrule height6.3ptdepth-.5ptwidth.8pt\kern-.8pt}
   {\vrule height4.1ptdepth-.35ptwidth.6pt\kern-.6pt}
   {\vrule height3.1ptdepth-.25ptwidth.5pt\kern-.5pt}}
\def\f@dge{\mathchoice{}{}{\mkern.5mu}{\mkern.8mu}}
\def\b@c#1#2{{\rm \mkern#2mu\vbar\mkern-#2mu#1}}
\def\b@b#1{{\rm I\mkern-3.5mu #1}}
\def\b@a#1#2{{\rm #1\mkern-#2mu\f@dge #1}}
\def\bb#1{{\count4=`#1 \advance\count4by-64 \ifcase\count4\or\b@a A{11.5}\or
   \b@b B\or\b@c C{5}\or\b@b D\or\b@b E\or\b@b F \or\b@c G{5}\or\b@b H\or
   \b@b I\or\b@c J{3}\or\b@b K\or\b@b L \or\b@b M\or\b@b N\or\b@c O{5} \or
   \b@b P\or\b@c Q{5}\or\b@b R\or\b@a S{8}\or\b@a T{10.5}\or\b@c U{5}\or
   \b@a V{12}\or\b@a W{16.5}\or\b@a X{11}\or\b@a Y{11.7}\or\b@a Z{7.5}\fi}}

\catcode`\X=11 \catcode`\@=12




\let\thischap\jobname

\def\partof#1{\csname returnthe#1part\endcsname}
\def\chapof#1{\csname returnthe#1chap\endcsname}

\def\setchapter#1,#2,#3;{%
  \expandafter\def\csname returnthe#1part\endcsname{#2}%
  \expandafter\def\csname returnthe#1chap\endcsname{#3}%
}

\setchapter 300a,A,II.A;
\setchapter 300b,A,II.B;
\setchapter 300c,A,II.C;
\setchapter 300d,A,II.D;
\setchapter 300e,A,II.E;
\setchapter 300f,A,II.F;
\setchapter 300g,A,II.G;
\setchapter  E53,B,N;
\setchapter  88r,B,I;
\setchapter  600,B,III;
\setchapter  705,B,IV;
\setchapter  734,B,V;

\def\cprefix#1{
\edef\theotherpart{\partof{#1}}\edef\theotherchap{\chapof{#1}}%
\ifx\theotherpart\thispart
   \ifx\theotherchap\thischap 
    \else 
     \theotherchap%
    \fi
   \else 
     \theotherchap\fi}

\def\sectioncite[#1]#2{%
     \cprefix{#2}#1}

\edef\thispart{\partof{\thischap}}
\edef\thischap{\chapof{\thischap}}

\def\lastpage of '#1' is #2.{\expandafter\def\csname lastpage#1\endcsname{#2}}


\def\spuriousreset{}


\expandafter\ifx\csname citeadd.tex\endcsname\relax
\expandafter\gdef\csname citeadd.tex\endcsname{}
\else \message{Hey!  Apparently you were trying to
\string\input{citeadd.tex} twice.   This does not make sense.} 
\errmessage{Please edit your file (probably \jobname.tex) and remove
any duplicate ``\string\input'' lines}\endinput\fi

\sectno=-1   
\localtags
\jjtags
\NoBlackBoxes
\define\mr{\medskip\roster}
\define\sn{\smallskip\noindent}
\define\mn{\medskip\noindent}
\define\bn{\bigskip\noindent}
\define\ub{\underbar}
\define\wilog{\text{without loss of generality}}
\define\ermn{\endroster\medskip\noindent}

\define \nl{\newline}
\magnification=\magstep 1
\documentstyle{amsppt}

{    
\catcode`@11

\ifx\alicetwothousandloaded@\relax
  \endinput\else\global\let\alicetwothousandloaded@\relax\fi

\gdef\subjclass{\let\savedef@\subjclass
 \def\subjclass##1\endsubjclass{\let\subjclass\savedef@
   \toks@{\def\usualspace{{\rm\enspace}}\eightpoint}%
   \toks@@{##1\unskip.}%
   \edef\thesubjclass@{\the\toks@
     \frills@{{\noexpand\rm2000 {\noexpand\it Mathematics Subject
       Classification}.\noexpand\enspace}}%
     \the\toks@@}}%
  \nofrillscheck\subjclass}
} 


\expandafter\ifx\csname alice2jlem.tex\endcsname\relax
  \expandafter\xdef\csname alice2jlem.tex\endcsname{\the\catcode`@}
\else \message{Hey!  Apparently you were trying to
\string\input{alice2jlem.tex}  twice.   This does not make sense.}
\errmessage{Please edit your file (probably \jobname.tex) and remove
any duplicate ``\string\input'' lines}\endinput\fi

\expandafter\ifx\csname bib4plain.tex\endcsname\relax
  \expandafter\gdef\csname bib4plain.tex\endcsname{}
\else \message{Hey!  Apparently you were trying to \string\input
  bib4plain.tex twice.   This does not make sense.}
\errmessage{Please edit your file (probably \jobname.tex) and remove
any duplicate ``\string\input'' lines}\endinput\fi

\def\renewcommand{\newcommand}	       
\edef\cite{\the\catcode`@}%
\catcode`@ = 11
\let\@oldatcatcode = \cite
\chardef\@letter = 11
\chardef\@other = 12
%
%
%
%
\def\@innerdef#1#2{\edef#1{\expandafter\noexpand\csname #2\endcsname}}%
%
%
\@innerdef\@innernewcount{newcount}%
\@innerdef\@innernewdimen{newdimen}%
\@innerdef\@innernewif{newif}%
\@innerdef\@innernewwrite{newwrite}%
%
%
%
\def\@gobble#1{}%
%
%
%
\ifx\inputlineno\@undefined
   \let\@linenumber = \empty 
\else
   \def\@linenumber{\the\inputlineno:\space}%
\fi
%
%
%
\def\@futurenonspacelet#1{\def\cs{#1}%
   \afterassignment\@stepone\let\@nexttoken=
}%
\begingroup 
\def\\{\global\let\@stoken= }%
\\ 
\endgroup
\def\@stepone{\expandafter\futurelet\cs\@steptwo}%
\def\@steptwo{\expandafter\ifx\cs\@stoken\let\@@next=\@stepthree
   \else\let\@@next=\@nexttoken\fi \@@next}%
\def\@stepthree{\afterassignment\@stepone\let\@@next= }%
%
%
%
\def\@getoptionalarg#1{%
   \let\@optionaltemp = #1%
   \let\@optionalnext = \relax
   \@futurenonspacelet\@optionalnext\@bracketcheck
}%
%
%
\def\@bracketcheck{%
   \ifx [\@optionalnext
      \expandafter\@@getoptionalarg
   \else
      \let\@optionalarg = \empty
      \expandafter\@optionaltemp
   \fi
}%
\def\@@getoptionalarg[#1]{%
   \def\@optionalarg{#1}%
   \@optionaltemp
}%
%
%
%
\def\@nnil{\@nil}%
\def\@fornoop#1\@@#2#3{}%
\def\@for#1:=#2\do#3{%
   \edef\@fortmp{#2}%
   \ifx\@fortmp\empty \else
      \expandafter\@forloop#2,\@nil,\@nil\@@#1{#3}%
   \fi
}%
\def\@forloop#1,#2,#3\@@#4#5{\def#4{#1}\ifx #4\@nnil \else
       #5\def#4{#2}\ifx #4\@nnil \else#5\@iforloop #3\@@#4{#5}\fi\fi
}%
\def\@iforloop#1,#2\@@#3#4{\def#3{#1}\ifx #3\@nnil
       \let\@nextwhile=\@fornoop \else
      #4\relax\let\@nextwhile=\@iforloop\fi\@nextwhile#2\@@#3{#4}%
}%
%
%
%
\@innernewif\if@fileexists
\def\@testfileexistence{\@getoptionalarg\@finishtestfileexistence}%
\def\@finishtestfileexistence#1{%
   \begingroup
      \def\extension{#1}%
      \immediate\openin0 =
         \ifx\@optionalarg\empty\jobname\else\@optionalarg\fi
         \ifx\extension\empty \else .#1\fi
         \space
      \ifeof 0
         \global\@fileexistsfalse
      \else
         \global\@fileexiststrue
      \fi
      \immediate\closein0
   \endgroup
}%
%
%
%
%
\def\bibliographystyle#1{%
   \@readauxfile
   \@writeaux{\string\bibstyle{#1}}%
}%
\let\bibstyle = \@gobble
%
%
\let\bblfilebasename = \jobname
\def\bibliography#1{%
   \@readauxfile
   \@writeaux{\string\bibdata{#1}}%
   \@testfileexistence[\bblfilebasename]{bbl}%
   \if@fileexists
      \nobreak
      \@readbblfile
   \fi
}%
\let\bibdata = \@gobble
%
%
\def\nocite#1{%
   \@readauxfile
   \@writeaux{\string\citation{#1}}%
}%
\@innernewif\if@notfirstcitation
%
%
\def\cite{\@getoptionalarg\@cite}%
%
%
\def\@cite#1{%
   \let\@citenotetext = \@optionalarg
   \printcitestart
   \nocite{#1}%
   \@notfirstcitationfalse
   \@for \@citation :=#1\do
   {%
      \expandafter\@onecitation\@citation\@@
   }%
   \ifx\empty\@citenotetext\else
      \printcitenote{\@citenotetext}%
   \fi
   \printcitefinish
}%
\newif\ifweareinprivate
\weareinprivatetrue
\ifx\shlhetal\undefinedcontrolseq\weareinprivatefalse\fi
\ifx\shlhetal\relax\weareinprivatefalse\fi
\def\@onecitation#1\@@{%
   \if@notfirstcitation
      \printbetweencitations
   \fi
   \expandafter \ifx \csname\@citelabel{#1}\endcsname \relax
      \if@citewarning
         \message{\@linenumber Undefined citation `#1'.}%
      \fi
     \ifweareinprivate
      \expandafter\gdef\csname\@citelabel{#1}\endcsname{%
\strut 
\vadjust{\vskip-\dp\strutbox
\vbox to 0pt{\vss\parindent0cm \leftskip=\hsize 
\advance\leftskip3mm
\advance\hsize 4cm\strut\openup-4pt 
\rightskip 0cm plus 1cm minus 0.5cm ?  #1 ?\strut}}
         {\tt
            \escapechar = -1
            \nobreak\hskip0pt\pfeilsw
            \expandafter\string\csname#1\endcsname
             \pfeilso
            \nobreak\hskip0pt
         }%
      }%
     \else  
      \expandafter\gdef\csname\@citelabel{#1}\endcsname{%
            {\tt\expandafter\string\csname#1\endcsname}
      }%
     \fi  
   \fi
   \csname\@citelabel{#1}\endcsname
   \@notfirstcitationtrue
}%
%
%
\def\@citelabel#1{b@#1}%
%
%
\def\@citedef#1#2{\expandafter\gdef\csname\@citelabel{#1}\endcsname{#2}}%
%
%
%
\def\@readbblfile{%
   \ifx\@itemnum\@undefined
      \@innernewcount\@itemnum
   \fi
   \begingroup
      \def\begin##1##2{%
         \setbox0 = \hbox{\biblabelcontents{##2}}%
         \biblabelwidth = \wd0
      }%
      \def\end##1{}
      %
      %
      \@itemnum = 0
      \def\bibitem{\@getoptionalarg\@bibitem}%
      \def\@bibitem{%
         \ifx\@optionalarg\empty
            \expandafter\@numberedbibitem
         \else
            \expandafter\@alphabibitem
         \fi
      }%
      \def\@alphabibitem##1{%
         \expandafter \xdef\csname\@citelabel{##1}\endcsname {\@optionalarg}%
         \ifx\biblabelprecontents\@undefined
            \let\biblabelprecontents = \relax
         \fi
         \ifx\biblabelpostcontents\@undefined
            \let\biblabelpostcontents = \hss
         \fi
         \@finishbibitem{##1}%
      }%
      \def\@numberedbibitem##1{%
         \advance\@itemnum by 1
         \expandafter \xdef\csname\@citelabel{##1}\endcsname{\number\@itemnum}%
         \ifx\biblabelprecontents\@undefined
            \let\biblabelprecontents = \hss
         \fi
         \ifx\biblabelpostcontents\@undefined
            \let\biblabelpostcontents = \relax
         \fi
         \@finishbibitem{##1}%
      }%
      \def\@finishbibitem##1{%
         \biblabelprint{\csname\@citelabel{##1}\endcsname}%
         \@writeaux{\string\@citedef{##1}{\csname\@citelabel{##1}\endcsname}}%
         \ignorespaces
      }%
      %
      %
      \let\em = \bblem
      \let\newblock = \bblnewblock
      \let\sc = \bblsc
      \frenchspacing
      \clubpenalty = 4000 \widowpenalty = 4000
      \tolerance = 10000 \hfuzz = .5pt
      \everypar = {\hangindent = \biblabelwidth
                      \advance\hangindent by \biblabelextraspace}%
      \bblrm
      \parskip = 1.5ex plus .5ex minus .5ex
      \biblabelextraspace = .5em
      \bblhook
      \input \bblfilebasename.bbl
   \endgroup
}%
%
%
\@innernewdimen\biblabelwidth
\@innernewdimen\biblabelextraspace
%
%
%
\def\biblabelprint#1{%
   \noindent
   \hbox to \biblabelwidth{%
      \biblabelprecontents
      \biblabelcontents{#1}%
      \biblabelpostcontents
   }%
   \kern\biblabelextraspace
}%
%
%
%
\def\biblabelcontents#1{{\bblrm [#1]}}%
%
%
\def\bblrm{\rm}%
%
%
\def\bblem{\it}%
%
%
\def\bblsc{\ifx\@scfont\@undefined
              \font\@scfont = cmcsc10
           \fi
           \@scfont
}%
%
%
\def\bblnewblock{\hskip .11em plus .33em minus .07em }%
%
%
\let\bblhook = \empty
%
%
%
\def\printcitestart{[}
\def\printcitefinish{]}
\def\printbetweencitations{, }
\def\printcitenote#1{, #1}
%
%
%
\let\citation = \@gobble
%
%
%
\@innernewcount\@numparams
%
%
\def\newcommand#1{%
   \def\@commandname{#1}%
   \@getoptionalarg\@continuenewcommand
}%
%
%
\def\@continuenewcommand{%
   \@numparams = \ifx\@optionalarg\empty 0\else\@optionalarg \fi \relax
   \@newcommand
}%
%
%
\def\@newcommand#1{%
   \def\@startdef{\expandafter\edef\@commandname}%
   \ifnum\@numparams=0
      \let\@paramdef = \empty
   \else
      \ifnum\@numparams>9
         \errmessage{\the\@numparams\space is too many parameters}%
      \else
         \ifnum\@numparams<0
            \errmessage{\the\@numparams\space is too few parameters}%
         \else
            \edef\@paramdef{%
               \ifcase\@numparams
                  \empty  No arguments.
               \or ####1%
               \or ####1####2%
               \or ####1####2####3%
               \or ####1####2####3####4%
               \or ####1####2####3####4####5%
               \or ####1####2####3####4####5####6%
               \or ####1####2####3####4####5####6####7%
               \or ####1####2####3####4####5####6####7####8%
               \or ####1####2####3####4####5####6####7####8####9%
               \fi
            }%
         \fi
      \fi
   \fi
   \expandafter\@startdef\@paramdef{#1}%
}%
%
%
%
%
\def\@readauxfile{%
   \if@auxfiledone \else 
      \global\@auxfiledonetrue
      \@testfileexistence{aux}%
      \if@fileexists
         \begingroup
            \endlinechar = -1
            \catcode`@ = 11
            \input \jobname.aux
         \endgroup
      \else
         \message{\@undefinedmessage}%
         \global\@citewarningfalse
      \fi
      \immediate\openout\@auxfile = \jobname.aux
   \fi
}%
%
%
\newif\if@auxfiledone
\ifx\noauxfile\@undefined \else \@auxfiledonetrue\fi
%
%
%
%
\@innernewwrite\@auxfile
\def\@writeaux#1{\ifx\noauxfile\@undefined \write\@auxfile{#1}\fi}%
%
%
%
\ifx\@undefinedmessage\@undefined
   \def\@undefinedmessage{No .aux file; I won't give you warnings about
                          undefined citations.}%
\fi
%
%
\@innernewif\if@citewarning
\ifx\noauxfile\@undefined \@citewarningtrue\fi
%
%
%
\catcode`@ = \@oldatcatcode

\def\pfeilso{\leavevmode
            \vrule width 1pt height9pt depth 0pt\relax
           \vrule width 1pt height8.7pt depth 0pt\relax
           \vrule width 1pt height8.3pt depth 0pt\relax
           \vrule width 1pt height8.0pt depth 0pt\relax
           \vrule width 1pt height7.7pt depth 0pt\relax
            \vrule width 1pt height7.3pt depth 0pt\relax
            \vrule width 1pt height7.0pt depth 0pt\relax
            \vrule width 1pt height6.7pt depth 0pt\relax
            \vrule width 1pt height6.3pt depth 0pt\relax
            \vrule width 1pt height6.0pt depth 0pt\relax
            \vrule width 1pt height5.7pt depth 0pt\relax
            \vrule width 1pt height5.3pt depth 0pt\relax
            \vrule width 1pt height5.0pt depth 0pt\relax
            \vrule width 1pt height4.7pt depth 0pt\relax
            \vrule width 1pt height4.3pt depth 0pt\relax
            \vrule width 1pt height4.0pt depth 0pt\relax
            \vrule width 1pt height3.7pt depth 0pt\relax
            \vrule width 1pt height3.3pt depth 0pt\relax
            \vrule width 1pt height3.0pt depth 0pt\relax
            \vrule width 1pt height2.7pt depth 0pt\relax
            \vrule width 1pt height2.3pt depth 0pt\relax
            \vrule width 1pt height2.0pt depth 0pt\relax
            \vrule width 1pt height1.7pt depth 0pt\relax
            \vrule width 1pt height1.3pt depth 0pt\relax
            \vrule width 1pt height1.0pt depth 0pt\relax
            \vrule width 1pt height0.7pt depth 0pt\relax
            \vrule width 1pt height0.3pt depth 0pt\relax}

\def\pfeilsw{ \leavevmode 
            \vrule width 1pt height0.3pt depth 0pt\relax
            \vrule width 1pt height0.7pt depth 0pt\relax
            \vrule width 1pt height1.0pt depth 0pt\relax
            \vrule width 1pt height1.3pt depth 0pt\relax
            \vrule width 1pt height1.7pt depth 0pt\relax
            \vrule width 1pt height2.0pt depth 0pt\relax
            \vrule width 1pt height2.3pt depth 0pt\relax
            \vrule width 1pt height2.7pt depth 0pt\relax
            \vrule width 1pt height3.0pt depth 0pt\relax
            \vrule width 1pt height3.3pt depth 0pt\relax
            \vrule width 1pt height3.7pt depth 0pt\relax
            \vrule width 1pt height4.0pt depth 0pt\relax
            \vrule width 1pt height4.3pt depth 0pt\relax
            \vrule width 1pt height4.7pt depth 0pt\relax
            \vrule width 1pt height5.0pt depth 0pt\relax
            \vrule width 1pt height5.3pt depth 0pt\relax
            \vrule width 1pt height5.7pt depth 0pt\relax
            \vrule width 1pt height6.0pt depth 0pt\relax
            \vrule width 1pt height6.3pt depth 0pt\relax
            \vrule width 1pt height6.7pt depth 0pt\relax
            \vrule width 1pt height7.0pt depth 0pt\relax
            \vrule width 1pt height7.3pt depth 0pt\relax
            \vrule width 1pt height7.7pt depth 0pt\relax
            \vrule width 1pt height8.0pt depth 0pt\relax
            \vrule width 1pt height8.3pt depth 0pt\relax
            \vrule width 1pt height8.7pt depth 0pt\relax
            \vrule width 1pt height9pt depth 0pt\relax
      }


\def\widestnumber#1#2{}

\def\citewarning#1{\ifx\shlhetal\relax 
    \else
    \par{#1}\par
    \fi
}

\def\rm{\fam0 \tenrm}

\def\fakesubhead#1\endsubhead{\bigskip\noindent{\bf#1}\par}



%
%
%

%

\font\textrsfs=rsfs10
\font\scriptrsfs=rsfs7
\font\scriptscriptrsfs=rsfs5

\newfam\rsfsfam
\textfont\rsfsfam=\textrsfs
\scriptfont\rsfsfam=\scriptrsfs
\scriptscriptfont\rsfsfam=\scriptscriptrsfs

\edef\oldcatcodeofat{\the\catcode`\@}
\catcode`\@11

\def\Cal@@#1{\noaccents@ \fam \rsfsfam #1}

\catcode`\@\oldcatcodeofat


\expandafter\ifx \csname margininit\endcsname \relax\else\margininit\fi

\long\def\red#1\endred{}
\long\def\green#1\endgreen{}
\long\def\blue#1\endblue{}
\long\def\private#1\endprivate{}

\def\endred{ \unmatched endred! }
\def\endgreen{ \unmatched endgreen! }
\def\endblue{ \unmatched endblue! }
\def\endprivate{ \unmatched endprivate! }

\ifx\latexcolors\undefinedcs\def\latexcolors{}\fi

\def\emptycs{}
\def\evaluatelatexcolors{%
        \ifx\latexcolors\emptycs\else
        \expandafter\xxevaluate\latexcolors\xxfertig\evaluatelatexcolors\fi}
\def\xxevaluate#1,#2\xxfertig{\setupthiscolor{#1}%
        \def\latexcolors{#2}}


\font\smallfont=cmsl7
\def\rutgerscolor{\ifmmode\else\endgraf\fi\smallfont
\advance\leftskip0.5cm\relax}
\def\setupthiscolor#1{\edef\tmptmpcs{\noexpand\bgroup\noexpand\rutgerscolor
\noexpand\def\noexpand\currentcolor{#1}%
\noexpand}%
\expandafter\let\csname#1\endcsname\tmptmpcs
\def\tmptmpcs{\checkColorUnmatched{#1}\popthecolor}
\expandafter\let\csname end#1\endcsname\tmptmpcs}

\def\checkColorUnmatched#1{\def\expectcolor{#1}%
    \ifx\expectcolor\currentcolor   
    \else \edef\failhere{\noexpand\tryingToClose '\currentcolor' with end\expectcolor}\failhere\fi}

\def\currentcolor{???}

\def\popthecolor{\ifmmode\else\endgraf\fi\egroup}

\expandafter\def\csname#1\endcsname{}

\evaluatelatexcolors

 \let\outerhead\head
 \def\head{\innerhead}
 \let\innerhead\outerhead

 \let\outersubhead\subhead
 \def\subhead{\innersubhead}
 \let\innersubhead\outersubhead

 \let\outersubsubhead\subsubhead
 \def\subsubhead{\innersubsubhead}
 \let\innersubsubhead\outersubsubhead

 \def\proclaim{\innerproclaim}
 \let\innerproclaim\outerproclaim

 %
 %
 %
 %

\def\demo#1{\medskip\noindent{\it #1.\/}}
\def\enddemo{\smallskip}

\def\remark#1{\medskip\noindent{\it #1.\/}}
\def\endremark{\smallskip}

\pageheight{8.5truein}
\topmatter
\title{The height of the automorphism tower of a group} \endtitle
\author {Saharon Shelah \thanks {\null\newline I would like to thank 
Alice Leonhardt for the beautiful typing. \null\newline
Research partially supported by NSF Grant No. NSF-DMS 0100794. Publication 810 
} \endthanks} \endauthor 

\affil{Institute of Mathematics\\
 The Hebrew University of Jerusalem \\
 Einstein Institute of Mathematics \\
  Edmond J. Safra Campus, Givat Ram \\
 Jerusalem 91904, Israel
 \medskip
 Department of Mathematics \\
 Hill Center-Busch Campus \\
  Rutgers, The State University of New Jersey \\
 110 Frelinghuysen Road \\
 Piscataway, NJ 08854-8019 USA} \endaffil

\abstract    For a group $G$ with trivial center there is a natural
embedding of $G$ into its automorphism group, so we can look at the
latter as an extension of the group.  So an increasing continuous
sequence of groups, the automorphism tower, is defined, the height is
the ordinal where this becomes fixed, arriving to a complete group.
We show that for many such $\kappa$ there is a group of
height $> 2^\kappa$, so proving that the upper bound essentially
cannot be improved. \endabstract
\endtopmatter
\document

\newpage

\head {\S0 Introduction} \endhead  \resetall \sectno=0
 \spuriousreset
\bigskip

For a group $G$ with trivial center there is a natural embedding of
$G$ into its automorphism group Aut$(G)$ where $g \in G$ is mapped to the
inner automorphism $x \mapsto g x g^{-1}$ which is defined and is not
the identity for $g \ne e_G$ as $G$ has a trivial center, so
we can view Aut$(G)$ as a group extending $G$.  Also the extension
Aut$(G)$ is a group with trivial center, so we can continue defining
$G^{<\alpha>}$ increasing with $\alpha$ for every ordinal $\alpha$; 
let $\tau_G$ be when we stop,
i.e., the first $\alpha$ such that $G^{< \alpha +1>} = G^{<\alpha>}$ 
(or $\alpha = \infty$ but see below) hence
$\beta > \alpha \Rightarrow G^{<\beta>} = G^{<\alpha>}$, (see Definition
\scite{0.1}).  How large can $\tau_G$ be?

Weilant \cite{Wel39} proves that for finite $G,\tau_G$ is finite.
Thomas \cite{Th85} celebrated work proves for infinite $G$ that 
$\tau_G \le (2^{|G|})^+$, in fact as noted by Felgner and Thomas
$\tau_G < (2^{|G|})^+$.  Thomas shows also that $\tau_\kappa \ge
\kappa^+$.  Later he (\cite{Th98}) showed that if $\kappa =
\kappa^{< \kappa},2^\kappa = \kappa^+$ (hence $\tau_\kappa <
\kappa^{++}$ in $\bold V$) and $\lambda \ge \kappa^{++}$ and we force by
$\Bbb P$, the forcing of
adding $\lambda$ Cohen subsets to $\kappa$, then in 
$\bold V^{\Bbb P}$ we still have $\tau_\kappa < \kappa^{++}$ though
$2^\kappa$ is $\ge \lambda$ (and $\bold V,\bold V^{\Bbb P}$ has the
same cardinals).  \nl

Just Shelah and Thomas \cite{JShT:654} prove that when $\kappa =
\kappa^{< \kappa} < \lambda$, in some forcing
extension (by a specially constructed $\kappa$-complete $\kappa^+$-c.c. forcing
notion) we have $\tau_\kappa \ge \lambda$, so consistently
$\tau_\kappa > 2^\kappa > \kappa^+$ for some $\kappa$.  
An important lemma there which
we shall use (see \scite{0.5} below) 
is that if $G$ is the automorphism group of a structure of
cardinality $\kappa,H \subseteq G,|H| \le \kappa$ then $\tau'_{G,H}$, 
the normalizer
length of $H$ in $G$ (see Definition \scite{0.2}(2)), is $<
\tau_\kappa$.  Concerning groups with center Hamkins show that $\tau_G
<$ the first strongly inaccessible cardinal $> |G|$.  On the subject
see the forthcoming book of Thomas.
\bn
We shall show, e.g.
\proclaim{\stag{0.0} Theorem}   If $\kappa$ is 
strong limit singular of uncountable cofinality \ub{then}
$\tau_\kappa > 2^\kappa$. 
\endproclaim
\bn
It would have been nice if the lower bound for $\tau_\kappa,\kappa^+$
would (consistently) be the correct one, 
but Theorem \scite{0.0} shows that this is
not so.
Note that Theorem \scite{0.0} shows that provably in ZFC, in general 
the upper bound $(2^\kappa)^+$ cannot be improved.
See Conclusion \scite{p.5} for proof of the theorem,
quoting results from
pcf theory.  We thank Simon Thomas, the referee and Itay Kaplan for
many valuable complaints detecting serious problems in earlier versions.

The program, described in a simplified way, is that for each 
so called ``$\kappa$-parameter $\bold p$" which
includes a partial order $I$ we define a group $G_{\bold p}$
and a two element subgroup $H_{\bold p}$ such that $\langle
\text{nor}^\alpha_{G_{\bold p}}(H_{\bold p}):\alpha \le { \text{\rm rk\/}}
^{< \infty}_I\rangle$  ``reflect" rk$^{< \infty}_I = { \text{\rm
rk\/}}_{\bold p}$, the natural rank on $I$ (see Definition 
\scite{m.1}), so in particular $\tau'_{G_{\bold p},H_{\bold p}} = 
\text{ rk}^{< \infty}_{\bold p}$.  (Actually in the end we shall get
only $H$ of cardinality $\le \kappa$).

We use an inverse system ${\frak s} = \langle {\bold p}_u,\pi_{u,v}:u \le_J v
\rangle$ of $\kappa$-parameters $\pi_{u,v}$ maps $I_{\bold p_v}$ to
$I_{\bold p_u}$; however, in general 
the $\pi_{u,v}$'s do not preserve order (but do preserve in
some weak global sense) where $J$ is an $\aleph_1$-directed
partial order.  Now for each $u \in J$, we can define the group
$G_{{\bold p}_u}$; and we can take inverse limit in two ways.
\bn
\ub{Way 1}:  The inverse limit $\bold p_{\frak s}$ (with $\pi_{u,{\frak s}}$
for $u \in J$ of ${\frak s}$) is a $\kappa$-parameter and so the 
group $G_{{\bold p}_{\frak s}}$ is well defined.
\mn
\ub{Way 2}:  The inverse system $\langle G_{{\bold p}_u},\hat \pi_{u,v}:u
\le_J v \rangle$, of groups were $\hat \pi_{u,v}$ is the (partial)
homomorphism from $G_{{\bold p}_v}$ to $G_{{\bold p}_u}$ induced 
by $\pi_{u,v}$, has an inverse limit $G_{\frak s}$. \nl
Now
\mr
\item "{$(A)$}"  concerning $G_{I_{\frak s}}$ we normally have good
control over rk$({\bold p}_{\frak s})$ hence on the normalizer length of
$H_{{\bold p}_{\frak s}}$ inside $G_{{\bold p}_{\frak s}}$
\sn
\item "{$(B)$}"  $G_{\frak s}$ is (more exactly can be
represented as good enough) inverse limit of groups of cardinality $\le
\kappa$ hence is isomorphic to Aug$({\frak A})$ for some structure of
cardinality $\le \kappa$
\sn
\item "{$(C)$}"  in the good case $G_{{\bold p}_{\frak s}} = G_{\frak s}$ so
we are done (by \scite{0.5}).
\ermn
In \S3 we work to get the main result. 
\nl
There are obvious possible improvement of the results here, say trying
to prove $\delta_\kappa \le \tau_\kappa$ (see Definition \scite{0.4}) 
for every $\kappa$.  But more importantly, a 
natural conjecture, at least for me was $\tau_\kappa =
\delta_\kappa$ because all the results so far on $\tau_\kappa$ has
parallel for $\delta_\kappa$ (though not inversely).  
In particular it seems reasonable that for $\kappa =
\aleph_0$ the lower bound was right, i.e., $\tau_\kappa = \omega_1$.
[We shall try to return to those problems in a sequel \cite{Sh:F579}.] 
\bigskip

\definition{\stag{0.1} Definition}  1) For a group $G$ with trivial
center, define the group $G^{<\alpha>}$ with trivial center for an
ordinal $\alpha$, increasing continuous with $\alpha$ such that
$G^{<0>} = G$ and $G^{<\alpha +1>}$ is the group of automorphisms of
$G^{<\alpha>}$ identifying $g \in G^{<\alpha>}$ with the inner
automorphisms it defines.  We may stipulate $G^{<-1>} = \{e_G\}$.
\nl
[We know that $G^{<\alpha>}$ is a group with trivial center increasing
continuous with $\alpha$ and for some $\alpha <
(2^{|G|+\aleph_0})^+$ we have 
$\beta > \alpha \Rightarrow G^{<\beta>} = G^{<\alpha>}$.]
\nl
2) The automorphism tower height of the group $G$ is 
$\tau_G = \tau^{\text{atw}}_G = 
\text{ Min}\{\alpha:G^{<\alpha>} = G^{<\alpha +1>}\}$;
clearly $\beta \ge \alpha \ge \tau_G \Rightarrow G^{<\beta>} =
G^{<\alpha>}$, atw stands for automorphism tower. \nl
3) Let $\tau_\kappa = \tau^{\text{atw}}_\kappa$ 
be the least ordinal $\tau$ such that $\tau(G) <
\tau$ for every group $G$ of cardinality $\le \kappa$; we call it 
the group tower ordinal of $\kappa$.
\enddefinition
\bn
Now we define normalizer (group theorist write $N_G(H)$, but probably
for others nor$_G(H)$ will be clearer, at least this is so for the author).
\definition{\stag{0.2} Definition}  1) Let $H$ be a subgroup of $G$.

We define nor$^\alpha_G(H)$, a subgroup of $G$, by induction on the
ordinal $\alpha$, increasing continuous with $\alpha$.  We may add
nor$^{-1}_G(H) = \{e_G\}$.
\mn
\ub{Case 1}: $\alpha = 0$.

$\text{nor}^0_G(H) = H$.
\mn
\ub{Case 2}: $\alpha = \beta +1$.

$\text{nor}^\alpha_G(H) = \text{ nor}_G(\text{nor}^\beta_G(H))$, see below.
\mn
\ub{Case 3}:  $\alpha$ a limit ordinal

$$
\text{nor}^\alpha_G(H) = \cup\{\text{nor}^\beta_G(H):\beta < \alpha\}
$$
\mn
where

$$
\align
\text{nor}_G(H) = \{g \in G:&g \text{ normalize } H, \text{ i.e. }
g N g^{-1} = N, \text{ equivalently} \\
  &(\forall x \in H)[g x g^{-1} \in H \and g^{-1} x g \in H]\}.
\endalign
$$
\mn
2) Let $\tau'_{G,H}= \tau^{\text{nlg}}_{G,H}$, the normalizer 
length of $H$ in $G$, be
Min$\{\alpha:\text{nor}^\alpha_G(H) = \text{ nor}^{\alpha +1}_G(H)\}$;
so $\beta \ge \alpha \ge \tau'_{G,H} = \text{ nor}^\beta_G(H) = \text{
nor}^\alpha_G(H)$; nlg stands for normalizer length.  \nl
3) Let $\tau'_\kappa = \tau^{\text{nlg}}_\kappa$ 
be the least ordinal $\tau$ such that $\tau >
\tau'_{G,H}$ whenever $G = \text{ Aut}({\frak A})$
for some structure ${\frak A}$ on $\kappa$ and $H \subseteq G$ is a
subgroup satisfying $|H| \le \kappa$.
\nl
4) $\tau''_\kappa = \tau^{\text{nlf}}_\kappa$ is the least
ordinal $\tau$ such that $\tau > \tau^{\text{nlf}}_{G,H}$ 
wherever $G = \text{ Aut}({\frak A}),
{\frak A}$ a structure of cardinality $\le \kappa,H$ a
subgroup of $G$ of cardinality $\le \kappa$ and nor$^{\infty}_G(H) =
\cup\{\text{nor}^\alpha_G(H):\alpha$ an ordinal$\} = G$.
\enddefinition
\bigskip

\definition{\stag{0.3} Definition}  We say that $G$ is a 
$\kappa$-automorphism group
\ub{if} $G$ is the automorphism group of some structure of cardinality
$\le \kappa$.
\enddefinition
\bigskip

\definition{\stag{0.4} Definition}  Let $\delta_\kappa =
\delta(\kappa)$ be the first ordinal
$\alpha$ such that there is no sentence $\psi \in \Bbb L_{\kappa^+,\omega}$
satisfying:
\mr
\item "{$(a)$}"  $\psi \vdash$ ``$<$ is a linear order"
\sn
\item "{$(b)$}"   for every $\beta < \alpha$ there is a model $M$ of
$\psi$ such that $(|M|,<^M)$ has order type $\ge \beta$
\sn
\item "{$(c)$}"  for every model $M$ of $\psi,(|M|,<^M)$ is a well
ordering.
\ermn
See on this, e.g.  \cite[VII,\S5]{Sh:c}.
\enddefinition
\bn
Our proof of better lower bounds rely on the following result 
from \cite{JShT:654}.
\proclaim{\stag{0.5} Lemma}  $\tau'_\kappa \le \tau_\kappa$.
\endproclaim
\bn
\margintag{0.6}\ub{\stag{0.6} Question}:  1) Is it 
consistent that for some $\kappa,\tau'_\kappa < \tau_\kappa$?  Is this
provable in ZFC?  Is the negation consistent? 
\nl
2) Similarly for the inequalities $\delta_\kappa < \tau'_\kappa$, (and
$\delta_\kappa < \tau'_\kappa < \tau_\kappa$).  
\nl
See on those in \cite{Sh:F579}.
\bigskip

\demo{\stag{0.7} Observation}  For every $\kappa \ge \aleph_0$ we have
$\tau^{\text{atw}}_\kappa \ge \tau^{\text{nlg}}_\kappa \ge
\tau^{\text{nlf}}_\kappa$. 
\enddemo
\bigskip

\demo{Proof}  By \scite{0.5} and checking the definitions of
$\tau^{\text{nlg}}_\kappa,\tau^{\text{nlf}}_\kappa$.  In fact we
mostly work on proving that in \scite{0.1}, $\tau^{\text{nlf}}_\kappa
> 2^\kappa$.
\enddemo
\bn
\ub{Notation}:  For a group $G$ and $A \subseteq G$ let $\langle
A\rangle_G$ be the subgroup of $G$ generated by $A$.
\bn
\centerline{$* \qquad * \qquad *$}
\bn
\ub{Explanation of the proof}:

We would like to derive the desired group from a partial order $I$
representing the ordinal desired as $\tau_{G,H}$ in some way and the
tower of normalizers of an appropriate subgroup will reflect.  It
seems natural to say that if $t \in I$ represent the ordinal $\alpha$
then the $s <_I t$ will represent ordinals $< \alpha$ so we use the
depth in $I$

$$
\text{dp}_I(t) = \cup\{\text{dp}_I(s)+1:s <_I t\}.
$$
\mn
For each $t \in I$ we will like to have a generator $g_t$ of the group
(really denoted by $g_{(<t>,<>)}$) take care of the normalizer tower
not sloping at $\alpha = \text{ dp}_I(t)$ say $g_t$ will be in the
$(\alpha +1)$-th normalizer but not in the $\alpha$-th normalizer.
But we need a witness for $g_t$ not being in earlier normalizer
$(\beta +1)$-th normalizer $\beta < \alpha$.

Now $\beta$ is represented by some $s <_I t$, so we have witness
$g_{(<(t,x),(<>)>)},g_{<(t,x),(1)>}$, the first in the beginning, the
second in the $(\beta +1)$-th normalizer not in the $\beta$-th
normalizer.  So we have a long normalizer tower of the subgroup
$G^{<0}_I$, the one generated by $\{g_{(\bar t,\eta)}:\eta(\ell) =0$
for some $\ell < \ell g(\eta)\}$.

However $G^{<0}_I$ is too big.  So we use a semi-direct product $K_I =
G_I * L_I$, where $L_I$ is an abelian group with every element of order
two, generated by $\{\bold h_{g G^{<0}_I}:g \in G^{<0}_I\}$ with $g_1
\bold h_{g G^{<0}_I} = \bold h_{(g_1 g)G^{<0}_I}$ and show that the
normalizer wins of the subgroup $H_I = \{e,\bold h_{G^{<0}_I}\}$ of
$K_I$ has the same height.

But we have to make $K_I$ a $\kappa$-automorphism group.  We only
almost have it: (and has too) we will represent it as aut$(M)/N$ for
some structure $M$ of cardinality $\le \kappa$ and normal subgroup $N$
of it of cardinality $\le \kappa$; this suffices.

\relax From where will $M$ come from?  We will represent $I$ as a universe
limit of some kind of ${\frak t} = \langle I_u,\pi_{u,v}:u \le_J
v\rangle$ where $I_u$ is a partial order of cardinality $\le
\kappa,\pi_{u,v}$ a mapping from $I_v$ to $I_u$ (commuting).  It
seemed a priori natural to have $\pi_{u,v}$ is order preserving but it
seemingly does not work out.  It seemed a priori natural to prove that
whenever ${\frak t}$ is as above there is a universe limit, etc.  We
find it more transparent to treat axiomatically: the limit is given
inside, i.e. as ${\frak s}$ which is ${\frak t} +$ a limit $v^*$; and
$J^{\frak t} = J^{\frak s} \backslash \{v^*\}$ is directed.

Also we demand that $J^{\frak t}$ is $\aleph_1$-directed (otherwise in
the limit we have words come.
\nl
We shall derive the structure $M$ from ${\frak t}$ so its automorphism
comes from members of $K_{I_u}(u \in J^{\frak t})$.  Well, not exactly
but for formal terms for it, to enable us to project to $u'
\le_{J[{\frak t}]} u$; as recall that $\pi_{u,v}$ does not necessarily
preserve order.  To make things smooth we demand that if $J^{\frak t}$
is a linear order (say cf$(\kappa))$ when as in the main case,
$\kappa$ is singular strong limit of uncountable cofinality.

More specifically, if $s,t \in I$ then for every large enough $u \in
J^{\frak t},s <_{I_{v^*}} t \Leftrightarrow \pi_{u,v^*}(s) <_{I_u}
\pi_{u,v}(t)$; note the order of the quantifiers.  Then we define a
structure $M$ derived from ${\frak t}$.  So the automorphism group of
$M$ is the inverse limit of groups which comes from the formal
definitions of elements of $K_{I_u}$'s.  Each depend on finitely many
generators, which in different $u$'s give different reduced forms.

Now they are defined from some $\bar t \in {}^k(I_u)$ using
``$I_{v^*}$ is the inverse limit..." the ``important" $t_u$'s, those
which really affect, well form an inverse system (\wilog \, the length
$k$ is constant on an end segment here we use ``$J^{\frak t}$ is
$\aleph_1$-directed) so for those $\ell$'s \, $\langle t_{u,\ell}:u
\in J^{\frak t}\rangle$ has limit $t_{v^*,\ell}$ say for $\ell <
k_*$. 

So $\langle t_{u^*,\ell}:\ell < k_*\rangle$ has the same quantifier
type in $I_u$ whenever $u_* \le u \le v^*$ for some $u_* < v^*$.  The
other $t$'s still has influence, so it is enough to find for them a
pseudo limit: $t_{v^*,\ell}$ such that they will have the same affect
on how the ``important" $t_{u,\ell}$ are used (this is the essential
limit).

All this gives an approximation to aut$(M) \cong K_{I_{v^*}}$.  They
almost mean that we divide by the subgroup of the automorphism of $M$
which are id$_{K_u}$ after $u \in J^{\frak t}$ large enough.  This is
a normal subgroup of cardinality $\le \kappa$ so we are done except
constructing such systems.  

\newpage

\head {\S1 The groups} \endhead  \resetall \sectno=1
 \spuriousreset
\bn
\ub{Discussion}:  Our aim is for a partial order $I$
to define a group $G = G_I$ and a
subgroup $H=H_I$ such that the normalizer length of $H$ inside $G$
reflects the depth of the well founded part of $I$.  Eventually we
would like to use $I$ of large depth such that $|H_I| \le \kappa$ and
the normalizer length of $H$ inside $G_I$ is $> \kappa$, even equal
to the depth of $I$.

For clarity we first define an approximation, in particular, $H$
appears only in \S2.
How do we define the group $G = G_I$ from the partial order $I$?
For each $t \in I$ we would like to have an element
associated with it (it is $g_{(<t>,<>)}$) such that it will ``enter" 
$\text{nor}^\alpha_G(H)$ exactly for $\alpha = \text{ rk}_I(t)+1$.  We
intend that among the generators of the group
commuting is the normal case, and we need witnesses
that $g_{(<t>,<>)} \notin 
\text{ nor}^{\beta +1}_G(H)$ wherever $\beta < \alpha = \text{
rk}_I(t),\beta >0$.  It is natural that if rk$_I(t_1) = \beta$ and
$t_1 <_I t_0 =: t$ then we use $t_1$ to represent $\beta$, 
as witness; more specifically, we construct the group such that conjugation by
$g_{(<t>,<>)}$ interchange $g_{(<t_0,t_1>,<0>)}$ and
$g_{(<t_0,s_0>,<1>)}$ and one of them, say $g_{(<t_0,t_1>,<0>)}$
belongs to $\text{nor}^{\beta +1}_G(H) \backslash \text{
nor}^\beta_G(H)$ whereas the other one, $g_{(<t_0,0>,<1>)}$, belongs
to nor$^1_G(H)$.  Iterating we get the elements 
$x \in X_I$ defined below.
  
In an earlier version, to ``start the induction", some additional
generators $g_{(\alpha,\ell)}(\alpha \in Z^I,\ell <2)$ 
were used to generate $H$ and not using all of them had helped to
make nor$^1_{G_I}(H_I)$ having the desired value.  However, we have to
decide for each $g_{(\bar t,\nu)}$ for $(\bar t,\nu)$ as above, for which
$g_{(\alpha,\ell)}(\alpha \in Z^I,\ell < 2)$ does conjugation
by $g_{(\bar t,\nu)}$ maps $g_{(\alpha,\ell)}$ to itself and for which it
does not.  For this we chose subsets $A_{(\bar t,\nu)}
\subseteq Z^I$ to code our decisions when 
$(\bar t,\nu)$ is as above and well defined, and
make the conjugation with the generators intended to generate
nor$^1_G(H)$ appropriately.

Now we do it by adding to $G$ an element $g_*$ of
order 2 getting $K_I$, commuting with $g \in G$ iff $g$ is intended
to be in the low level (e.g. $g_{(\bar t,\eta)},t_n \in I$ is minimal,
see notation below).

We could have in this section considered only a partial  order $I$, and
the groups $G_I$ (and later $K_I$) derived from it.  But as anyhow we shall
use it in the context of $\kappa$-p.o.w.i.s., we do it in this frame
(of course if $J^{\frak s} = \{u\}$, then ${\frak s}$ is essentially
just $I_u$).
\bn
Note that for our main result it suffices to deal with the case rk$(I)
< \infty$. 
\definition{\stag{m.1} Definition}  Let $I$ be a partial order (so $\ne
\emptyset$).
\nl
1)  rk$_I:I \rightarrow \text{ Ord} \cup \{\infty\}$ is defined by
rk$_I(t) \ge \alpha$ iff $(\forall \beta < \alpha)(\exists s <_I t)
[\text{rk}_I(s) \ge \beta]$. \nl
2) rk$^{< \infty}_I(t)$ is defined as rk$_I(t)$ if rk$_I(t) < \infty$ and is
defined as $\cup\{\text{rk}_I(s)+1:s$ satisfies $s <_I t$ and rk$_I(s) <
\infty\}$ in general.
\nl
3) Let rk$(I) = \cup \{\text{rk}_I(t) +1:t \in I\}$ stipulating
$\alpha < \infty = \infty +1$. \nl
4) rk$^{< \infty}_I = \text{ rk}^{< \infty}(I) =
\cup\{\text{rk}^{< \infty}_I(t) +1:t \in I\}$. \nl
5) Let $I_{[\alpha]} = \{t \in I:\text{rk}(t)=\alpha\}$.
\nl
6) $I$ is non-trivial \ub{when} 
$\{s:s \le_I t$ and $\text{rk}_I(s) \ge \beta\}$ is
infinite for every $t \in I$ satisfying rk$^{< \infty}_I(t) > \beta$
(used in the proof of \scite{m.11}(1); it is equivalent to demand
``rk$_I(s)=\beta$").
\nl
7) $I$ is explicitly non-trivial if each $E_I$-equivalence class is
infinite where $E_I = \{(t_1,t_2):t_2 \in I,t_2 \in I$ and
$(\forall s \in I)(s <_I t_1 \equiv s <_I t_2)\}$.
\enddefinition
\bigskip

\definition{\stag{m.2} Definition}  1) ${\frak s}$ is a
$\kappa$-p.o.w.i.s. (partial order weak inverse system) \ub{when}:
\mr
\item "{$(a)$}"  ${\frak s} = (J,\bar I,\bar \pi)$ so $J =
J^{\frak s} = J[{\frak s}],\bar I = \bar I^{\frak s},\bar \pi = \bar
\pi^{\frak s}$
\sn
\item "{$(b)$}"  $J$ is a directed partial order of cardinality $\le
\kappa$
\sn
\item "{$(c)$}"  $\bar I = \langle I_u:u \in J \rangle = \langle
I^{\frak s}_u:u \in J\rangle$ and we may write $I[u]$ or $I^{\frak s}[u]$
\sn
\item "{$(d)$}"  $I_u = I^{\frak s}_u$ is a partial 
order of cardinality $\le \kappa$
\sn
\item "{$(e)$}"  $\bar \pi = \langle \pi_{u,v}:u \le_J v \rangle$
\sn
\item "{$(f)$}"  $\pi_{u,v}$ is a partial mapping from $I_v$ into
$I_u$ (no preservation of order is required!)
\sn
\item "{$(g)$}" if $u \le_J v \le_J w$ then $\pi_{u,w} = \pi_{u,v}
\circ \pi_{v,w}$.
\ermn
2) ${\frak s}$ is a p.o.w.i.s. mean $\kappa$-p.o.w.i.s. for some
$\kappa$.
\nl
3) For $u \in J$ let $X_u = X^{\frak s}_u$ be the set of $x$ such that
for some $n < \omega$:
\mr
\item "{$(a)$}"  $x = (\bar t,\eta) = (\bar t^x,\eta^x)$
\sn
\item "{$(b)$}"  $\eta^x$ is a function from $\{0,\dotsc,n-1\}$ to $\{0,1\}$
\sn
\item "{$(c)$}"   $\bar t =
\langle t_\ell:\ell \le n \rangle = \langle t^x_\ell:\ell \le n
\rangle$ where $t_\ell \in I^{\frak s}_u$ is 
$<_{I^{\frak s}_u}$-decreasing, i.e., $t_n
<_{I^{\frak s}_u} t_{n-1} <_{I^{\frak s}_u} \ldots <_{I^{\frak s}_u} t_0$.
\ermn
3A) In fact for every partial order $I$ we define $X_I$ similarly, so
$X^{\frak s}_u = X_{I^{\frak s}[u]}$.
\nl
4) In part (3) for $x \in X^{\frak s}_u$ let $n(x) = 
\ell g(\bar t^x)-1$ and $t^x = t(x) := t^x_{n(x)}$.
\nl
5) For $x \in X^{\frak s}_u$ and $n \le n(x)$ let $y = x
\upharpoonleft n \in X^{\frak s}_u$ be defined by:

$$
\bar t^y := \bar t^x \restriction (n+1) = \langle t^x_0,\dotsc,t^x_n \rangle
$$

$$
\eta^y = \eta^x \restriction n(y) =: \eta^x \restriction \{0,\dotsc,n-1\}.
$$
\mn
6) We define rk$^1_u = \text{ rk}^{1,{\frak s}}_u$ and rk$^2_u =
\text{ rk}^{2,{\frak s}}_u$ as follows:
\mr
\item "{$(a)$}"    let rk$^1_u:X_u \rightarrow \text{ Ord } 
\cup \{\infty\}$ be defined by $x \in X_u \Rightarrow 
\text{ rk}^{1,{\frak s}}_u(x)  = \text{ rk}^1_u(x) = \text{ rk}_{I[u]}(t^x)$
\sn
\item "{$(b)$}"  let rk$^2_u:X_u \rightarrow \{-1\} \cup \text{ Ord} \cup 
\{\infty\}$ 
{\roster
\itemitem{ $(\alpha)$ }   if $x \in X_u$ and $\{\eta^x(\ell):\ell < n(x)\}
\subseteq \{1\}$ (e.g., $n(x)=0$) then let rk$^2_u(x) = \text{
rk}^{2,{\frak s}}_u(x) = \text{ rk}_{I[u]}(t(x))$
\sn
\itemitem{ $(\beta)$ }    if $x \in X_u$ and $\{\eta^x(\ell):\ell < n(x)\}
\nsubseteq \{1\}$ then let 
rk$^{2,{\frak s}}_u(x) = -1$ (yes, -1).
\endroster}
\ermn
7) We say that ${\frak s}$ is nice \ub{when} every $I^{\frak s}_u$ is
  non-trivial and $\pi_{u,w}$ is a function from $I_v$ into $I_u$,
  i.e., the domain of $\pi^{\frak s}_{u,v}$ is $I_v$.
\nl
8) $X^{< \alpha}_u :=
\{x \in X^{\frak s}_u:\text{rk}^2_u(x) < \alpha\}$ and 
$X^{\le \alpha}_u := \{x \in X^{\frak s}_u:\text{ rk}^2_u(x) 
\le \alpha\}$.  Note that $X^{\le \alpha}_u = X^{< \alpha +1}_u$ when
$\alpha < \infty$.  Of course, we may write $X^{< \alpha,{\frak
s}}_u,X^{\le \alpha,{\frak s}}_u$ and note that $X^{<0}_u = \{x \in X^{\frak
s}_u:0 \in \text{ Rang}(\eta^x)\}$.
\enddefinition
\bigskip

\definition{\stag{m.4} Definition}  Assume ${\frak s}$ is a
$\kappa$-p.o.w.i.s. and $u \in J^{\frak s}$.
\nl
1) Let $G_u = G^{\frak s}_u = G_u[{\frak s}]$ be the group 
generated by $\{g_x:x \in X^{\frak s}_u\}$ freely except the 
equations in $\Gamma_u = \Gamma^{\frak s}_u$ where $\Gamma_u$ consists of
\mr
\item "{$(a)$}"  $g^{-1}_x = g_x$, that is $g_x$ has order $2$, for
each $x \in X_u$
\sn
\item "{$(b)$}"  $g_{y_1} g_{y_2} = g_{y_2} g_{y_1}$ when $y_1,y_2 \in
X_u$ and $n(y_1) = n(y_2)$
\sn
\item "{$(c)$}"  $g_x g_{y_1} g^{-1}_x = g_{y_2}$ when
$\circledast^{u,{\frak s}}_{x,y_1,y_2}$, see below.
\ermn
1A) Let $\circledast_{x,y} = \circledast^u_{x,y} =
\circledast^{u,{\frak s}}_{x,y}$ means that
$\circledast_{x,y_1,y_2}$ for some $y_1,y_2$ such that $y \in
\{y_1,y_2\}$, see below.
\nl
1B) Let $\circledast_{x,y_1,y_2} = \circledast^u_{x,y_1,y_2} = 
\circledast^{u,{\frak s}}_{x,y_1,y_2}$ means that:
\mr
\item "{$(a)$}"  $x,y_1,y_2 \in X_u$
\sn
\item "{$(b)$}"  $n(x) < n(y_1) = n(y_2)$
\sn
\item "{$(c)$}"  $y_1 \upharpoonleft n(x) = y_2 \upharpoonleft n(x)$
\sn
\item "{$(d)$}"  $\bar t^{y_1} = \bar t^{y_2}$
\sn
\item "{$(e)$}"  for $\ell < n(y_1)$ we have:
$\eta^{y_1}(\ell) \ne \eta^{y_2}(\ell)$ \ub{iff} $\ell = n(x)
\wedge x = y_1 \upharpoonleft n(x)$.
\ermn
2)  Let $G^{< \alpha}_u = G^{< \alpha,{\frak s}}_u$ 
be defined similarly to $G^{\frak s}_u$
except that it is generated only by $\{g_x:x \in X^{< \alpha}_u\}$,
freely except the equations from $\Gamma^{< \alpha}_u = 
\Gamma^{< \alpha,{\frak s}}_u$, 
where $\Gamma^{< \alpha}_u $ is the set of equations from $\Gamma_u$ among 
$\{g_x:x \in X^{< \alpha}_u\}$.

Similarly $G^{\le \alpha}_u,\Gamma^{\le \alpha}_u$; note 
that $G^{\le \alpha}_u = 
G^{< \alpha +1}_u,\Gamma^{\le \alpha}_u = 
\Gamma^{< \alpha +1}_u$ if $\alpha < \infty$. 
\nl
3) For $X \subseteq X_u$ let $G_{u,X} = G^{\frak s}_{u,X}$ be 
the group generated by $\{g_y:y \in X\}$ freely except the equations in 
$\Gamma_{u,X} = \Gamma^{\frak s}_{u,X}$ which is the set of equations from 
$\Gamma_u$ mentioning only generators among $\{g_y:y \in X\}$.
\enddefinition
\bigskip

\demo{\stag{m.5} Observation}  1) The sequence
$\langle X^{< \alpha}_u:\alpha \le \text{ rk}(I^{\frak s}_u)
\rangle$ is $\subseteq$-increasing continuous.
\nl
2) If $\ell \in \{1,2\}$ and 
$x,y \in X_u$ are such that $x \ne y=x \upharpoonleft n$ and $\ell
\in \{1,2\}$ then rk$^\ell_{\bold p}(y) \ge \text{ rk}^\ell_{\bold
p}(x)$ and if equality holds then rk$^1_u(x) = \infty =
 \text{\rm rk}^1_u(y)$ \ub{or} both are $-1$ and $\ell=2$. 
\nl
3) If a partial order $I$ is explicitly non-trivial \ub{then} $I$ is
non-trivial. 
\enddemo
\bigskip

\demo{Proof}  Check.
\enddemo
\bigskip

\demo{\stag{m.6} Observation}  For a $\kappa$-p.o.w.i.s. ${\frak s}$.
\nl
1) $\circledast^{u,{\frak s}}_{x,y}$ holds \ub{iff}:
\mr
\item "{$(\alpha)$}"  $x,y \in X_u$ and 
\sn
\item "{$(\beta)$}"  $n(y) \ge n(x) +1$.
\ermn
2) If $x \in X^{\frak s}_u$ then 
$\{(y_1,y_2):\circledast^{u,{\frak s}}_{x,y_1,y_2}$ holds$\}$ is a
permutation of order two of 
$Y_{> n(x)} =: \{y \in X^{\frak s}_n:n(y) > n(x)\}$.
\nl
3) Moreover, the permutation in part (2) maps each $Y_{n+1} \backslash
Y_n$ onto itself when 
$n \in [n(x),\omega)$ and so it maps $\Gamma_{Y_{>n}}$ onto itself
 when $n(*) \le n <  \omega$.
\nl
4) If $\circledast^{u,{\frak s}}_{x,y_1,y_2}$ then 
$y_1 \upharpoonleft n(x) = y_2 \upharpoonleft n(x)$ and $n(x) <
n(y_1) = n(y_2)$.
\nl
5) $\circledast^{u,{\frak s}}_{x,y_1,y_2}$ iff 
$\circledast^{u,{\frak s}}_{x,y_2,y_1}$.
\nl
6) For $x,y \in X^{\frak s}_u$, in the group 
$G^{\frak s}_u$ the elements $g_x,g_y$ commute except when $x
\ne y \wedge (x=y \upharpoonleft n(x) \vee y = x \upharpoonleft
n(y))$.  In this case, if $n(x) < n(y)$ there is $y' \ne y$ such that
\mr
\item "{$\circledast_{x,y,y'}$}"   so $n(y') = n(y)$ and $\eta^y(\ell)
= \eta^{y'}(\ell) \Leftrightarrow \ell \in n(x)$.
\endroster  
\enddemo
\bigskip

\demo{Proof}  Straight (details on (2),(3) see the proof of
   \scite{m.11}).
\hfill$\square_{\scite{m.6}}$
\enddemo
\bn
We first sort out how elements in $G^{\frak s}_u$ and various subgroups
can be (uniquely) represented as products of the generators.
\proclaim{\stag{m.7} Claim}  Assume that ${\frak s}$ is a
$\kappa$-p.o.w.i.s., $u \in J^*$ and $<^*$ is any linear order of
$X_u$ such that 
\mr
\item "{$\boxdot$}"   if $x \in X_u,y \in X_u$ and $n(x) > n(y)$ 
\ub{then} $x <^* y$.
\ermn
1) Any member of $G_u$ is equal to a product of the form $g_{x_1}
\ldots g_{x_m}$ where $x_\ell <^* x_{\ell +1}$ for $\ell
=1,\dotsc,m - 1$.  Moreover, this representation is unique.
\nl
2) Similarly for $G^{\le \alpha}_u,G^{< \alpha}_u$
(using $X^{\le \alpha}_u,X^{< \alpha}_u$ respectively
instead $X_u$) hence $G^{\le \alpha}_u,G^{< \alpha}_u$ are subgroups
of $G_u$. 
\nl
3) In part (1) we can replace $G_u$ and $X_u$ by 
$G = G_{u,X}$ and $X$ respectively when $X \subseteq X_u$ is such that 
$[\{x,y_1,y_2\} \subseteq X_u \wedge \circledast^{u,{\frak s}}_{x,y_1,y_2}
\wedge \{x,y_1\} \subseteq X \Rightarrow y_2 \in X]$.  Hence $G_{u,X}$
is equal to $\langle \{g_x:x \in X\}\rangle_{G_u}$.
\nl
4) If $g = g_{y_1} \ldots g_{y_m}$ where $y_1,\dotsc,y_m \in X_u$
and $g = g_{x_1} \ldots g_{x_n} \in G_u$ and $x_1 <^* \ldots <^* x_n$
\ub{then} $n \le m$. \nl
5) $\langle G^{< \alpha}_u:\alpha \le { \text{\rm rk\/}}
(I^{\frak s}_u),\alpha$ an ordinal$\rangle$ is an 
increasing continuous sequence of groups with last element $G^{< \infty}_u$.
\nl
6) $\{g G^{< 0}_u:g \in G_u\}$ is a partition of $G_u$ (to left
cosets of $G_u$ over $G^{< 0}_u$).
\nl
7) If $<^1,<^2$ are two linear orders of $X_u$ as in $\boxdot$ above and
$G_u \models ``g_{x_1} \ldots g_{x_k} = g_{y_1} \ldots g_{y_m}"$ and $x_1 <^1
\ldots <^1 x_k$ and $y_1 <^2 \ldots <^2 y_m$ (or just $x_1 < \char 94
\ldots \char 94 x_k,n(y_1) \ge n(y_2) \ge \ldots n(y_n)$ and $\langle
y_\ell:\ell=1,m\rangle$ is with no repetitions), \ub{then}:
\mr
\item "{$(\alpha)$}"  $k=m$
\sn
\item "{$(\beta)$}"  for every $i$ we have $\{\ell:n(x_\ell)=i\} =
\{\ell:n(y_\ell) = i\}$ and this set is a convex subset of $\{1,\dotsc,m\}$.
\ermn
(So the only difference is permuting $g_{x_{\ell(2)}},g_{x_{\ell(1)}}$
when $n(x_{\ell(1)}) = n(x_{\ell(2)})$.
\nl
8) If $I \subseteq I_u$ and $X = X_I$ then $G_{u,X} \cap G^{< 0}_u$ is
the subgroup of $G_{u,X}$ generated by $\{g_x:x \in
X,\text{\rm Rang}(\eta^x) \nsubseteq 1\}$, i.e., the (naturally defined)
$G^{< 0}_I$. 
\nl
9) If $I_\ell \subseteq I^{\frak s}_u$ for $\ell=1,2,3$ (so
$\le_{I_\ell} = \le_I \restriction I_\ell$) and $I_1 \cap I_2 =
I_3$ \ub{then} $G_{I_1} \cap G_{I_2} = G_{I_3}$ and $G^{< 0}_{I_1}
\cap G^{<0}_{I_2} = G^{<0}_{I_3}$.
\endproclaim
\bigskip

\demo{Proof}  1),2),3)  Recall that each generator has order two.  We can 
use standard
combinatorial group theory (the rewriting process but below we do not
assume knowledge of it); the point is that
in the rewriting the number of generators in the word do not
increase (so no need of $<^*$ being a well ordering). 
\nl
We now give a full self-contained proof, for part of (2) we consider $G = 
G^{< \alpha}_u,X = X^{< \alpha}_u \subseteq X_u,\Gamma = \Gamma^{< \alpha}_u$ 
for $\alpha$ an ordinal or infinity and for part (1) and 
the rest of part (2) consider $G = G^{\le \beta}_u,X = X^{\le \beta}_u 
\subseteq X_u,\Gamma = \Gamma^{\le \beta}_u$ for 
$\beta$ an ordinal or infinity (recall that $G_u,X_u$
is the case $\beta = \infty$).  Now in parts (1),(2) for the set $X$, 
the condition from part (3) holds by \scite{m.5}(2).
\nl
[Why?  So assume $\circledast^u_{x,y_1,y_2}$ and e.g. $x,y_1 \in
X^{\le \alpha}_u$ and we should prove that $y_2 \in X^{\le
\alpha}_u$.  If $y_1 = y_2$ this is trivial so assume $y_1 \ne y_2$,
hence necessarily $y_1 \upharpoonleft n(x) = x = y_2 \upharpoonleft
n(x)$ and $n(x) < n(y_1) = n(y_2)$ and $\bar t^{y_1} = \bar t^{y_2}$
and $\eta^{y_1}(\ell) = \eta^{y_2}(\ell) \Leftrightarrow \ell \ne n(x)$.
If $\eta^x$ is not constantly one then also $\eta^{y_1}$ is not
constantly one hence $y_2 \in X^{< 0}_u$ so fine.  If
$\eta^x$ is constantly one then $\alpha \ge \text{ rk}^1_u(t^x) >
\text{ rk}^1_u(t^{y_1}) = \text{ rk}^1_u(t^{y_2}) \ge \text{
rk}^2_u(t^{y_2})$ hence $y_2 \in X^{\le \alpha}_u$ so fine.]
\nl
So it is enough to prove part (3).  Now recall that 
$G = G_{u,X}$ and 
\mr
\item "{$\circledast_1$}"  every member of $G$ can be written as a
product $g_{x_1} \ldots g_{x_n}$ for some $n < \omega,x_\ell \in X$
\nl
[Why?  As the set $\{g_x:x \in X\}$ generates $G$.]
\sn
\item "{$\circledast_2$}"  if in $g = g_{x_1} \ldots g_{x_n}$ we have
$x_\ell = x_{\ell +1}$ \ub{then} we can omit both
\nl
[Why?  As $g_x g_x = e_G$ for every $x \in X$ by clause (a) of
Definition \scite{m.4}(1)]
\sn
\item "{$\circledast_3$}"  if $1 \le \ell < n$ and 
$g = g_{x_1} \ldots g_{x_n}$ and we have
$x_{\ell +1} <^* x_\ell$ and $m \in \{1,\dotsc,n\}
\backslash \{\ell,\ell +1\} \Rightarrow y_m = x_m$ 
\ub{then} we can find $y_\ell,y_{\ell+1} \in X$
such that $g = g_{y_1} \ldots g_{y_n}$ and $y_\ell <^* y_{\ell +1}$
and, in fact, $y_{\ell +1} = x_\ell$.
\ermn
[Why does $\circledast_3$ hold?  By Definition \scite{m.4}(1) and
Observation \scite{m.6}(6) one of the following cases occurs.
\bn
\ub{Case 1}:  $g_{x_\ell},g_{x_{\ell +1}}$ commutes. 

Let $y_\ell = x_{\ell +1},y_{\ell +1} = x_\ell$.
\mn
\ub{Case 2}:  Not Case 1 but  
$\circledast^{u,{\frak s}}_{x_{\ell +1},x_\ell}$, see 
Definition \scite{m.4}(1A).

By clause (b) of Definition \scite{m.4}(1) we have $n(x_{\ell +1})
< n(x_\ell)$. 
So by $\boxdot$ of the assumption of the present claim we have
$x_\ell <^* x_{\ell +1}$, contradiction.
\mn
\ub{Case 3}:  Not case 1 but 
$\circledast^{u,{\frak s}}_{x_\ell,x_{\ell +1}}$, see
Definition \scite{m.4}(1B).

By \scite{m.6}(6) there is $y_\ell \in X$ such that
$n(y_\ell) = n(x_{\ell +1}) > n(x_\ell),\bar t^{y_\ell} = 
\bar t^{x_{\ell +1}}$ and $i < n(x_{\ell +1}) \Rightarrow 
(\eta^{y_\ell}(i) = \eta^{x_{\ell +1}}(i)) \equiv(i
\ne n(x_\ell))$.

Let $y_{\ell +1} = x_\ell$, clearly $y_{\ell +1},y_\ell \in X$.  By
Definition \scite{m.4}(1), we have 
$g_{x_\ell} g_{x_{\ell +1}} g^{-1}_{x_\ell} = g_{y_\ell}$ hence
$g_{x_\ell} g_{x_{\ell +1}} = g_{y_\ell} g_{x_\ell} = g_{y_\ell}
g_{y_{\ell +1}}$ and clearly $n(y_{\ell +1}) = n(x_\ell) < n(y_\ell)$ hence
$y_\ell <^* x_\ell = y_{\ell +1}$, so we are
done.  \nl
The three cases exhaust all possibilities hence $\circledast_3$ is proved.]
\mr
\item "{$\circledast_4$}"  every $g \in G$ can be represented as
$g_{x_1} \ldots g_{x_n}$ with $x_1 <^* x_2 <^* \ldots <^* x_n$.
\ermn
[Why?  Without loss of generality $g$ is not the unit of $G$.
By $\circledast_1$ we can find $x_1,\dotsc,x_n \in X_1$ such that $g =
g_{x_1} \ldots g_{x_n}$ and $n \ge 1$.  Choose such a representation satisfying
\mr
\item "{$\otimes$}"   $(a) \quad$ with minimal $n$ and 
\sn
\item "{${{}}$}"  $(b) \quad$ for this $n$, 
with minimal $m \in \{1,\dotsc,n+1\}$ such
that $x_m <^* \ldots <^* x_n$ 
\nl

\hskip25pt and $1 \le \ell < m \le n \Rightarrow \dsize
\bigwedge^{m-1}_{\ell=1} x_\ell <^* x_m$, and
\sn
\item "{${{}}$}"   $(c) \quad$ for 
this pair $(n,m)$ if $m > 2$ then with maximal $\ell$ where $\ell \in$
\nl

\hskip25pt $\{1,\dotsc,m-1\}$ satisfies $x_\ell$ is $<^*$-maximal among
$\{x_1,\dotsc, x_{m-1}\}$

\hskip25pt that is $k \in \{1,\dotsc,m-1\} \Rightarrow x_k \le^* x_\ell$.
\ermn
Easily there is such a sequence $(x_1,\dotsc,x_n)$, noting that $m=n+1$
is O.K. for $(b)$ and there is $\ell$ as in $\otimes(c)$.

By $\circledast_2$ and clause (a) of $\otimes$ 
we have $x_\ell \ne x_{\ell +1}$ when $\ell$ from $\otimes(c)$ is 
well defined, i.e., if $m>2$). \nl
Now $m=2$ is impossible (as then $m=1$ can serve), if $m=1$ we are
done, and if $m>2$ then $\ell$ is well defined and 
$\ell =m-1$ is impossible (as then $m-1$ can serve
instead $m$).  Lastly by $\circledast_3$ applied to this $\ell$, we could
have improved $\ell$ to $\ell +1$, contradiction.]
\mr
\item "{$\circledast_5$}"  the representation in $\circledast_4$ is
unique.
\ermn
[Why does $\circledast_5$ hold?  Assume 
toward contradiction that $g_{x'_1} \ldots g_{x'_{n_1}} =
g_{y'_1} \ldots g_{y'_{n_2}}$ where $x'_1 <^* \ldots <^* x'_{n_1}$ and 
$y'_1 <^* \ldots <^* y'_{n_2}$ and $(x'_1,\dotsc,x'_{n_1}) \ne 
(y'_1,\dotsc,y'_{n_2})$.  Without loss of generality among all
such examples, $(n_1 + n_2 + 1)^2 + n_1$ is minimal.

Recall $Y_n =: \{x \in X:n(x)=n\}$.

So $\langle Y_n:n < \omega \rangle$ is a partition of $X$.

For $k \le m < \omega$ let $X^{<k,m>} = 
\bigcup\{Y_\ell:k \le \ell < m\}$ and let $G^{<k,m>}$ be the group
generated  by $\{g_x:x \in X^{<k,m>}\}$ freely except the equations in
$\Gamma^{<k,n>}$, i.e., the equations from $\Gamma_{u,X^{<k,m>}}$, 
i.e., the equations from Definition \scite{m.4}(4)
mentioning only its generators, $\{y_x:x \in X^{<k,m>}\}$.  
Now clearly if $\circledast^{u,{\frak s}}_{x,y_1,y_2}$, see Definition
\scite{m.4}(1B) then $n < \omega \Rightarrow [y_2 \in Y_n \equiv y_2
\in Y_n]$ so the set $X^{<k,m>} \subseteq X$ satisfies the requirement
in part (3) of \scite{m.7} which we are proving; so what we have
proved for $X$ holds for $X^{<k,m>}$.  In particular
$\circledast_1 - \circledast_4$
above gives that for every $g \in G^{<k,m>}$ there are $n$ and
$x_1 <^* \ldots <^* x_n$ from $X^{<k,m>}$ such that $G^{<n,m>} \models
``g = g_{x_1} \ldots g_{x_n}"$.  Also it is enough to prove the uniqueness
for $G^{<k,m>}$ (for every $k \le m < \omega$), i.e., we can assume
$x'_1,\dotsc,x'_{n_1},y'_1,\dotsc,y'_{n_2} \in X^{<k,m>}$
as if the equality holds (though $\langle x'_1,\dotsc,x'_{n_1}\rangle
\ne \langle y'_1,\dotsc,y'_{n_2}\rangle$), finitely many equations of
$\Gamma_{u,X}$ implies the undesirable equation
and for some $k \le  m < \omega $ they are all from $\Gamma^{<k,m>}$
and $\{x'_1,\dotsc,x'_{n_1},y'_1,\dotsc,y'_{n_2}\} \subseteq X$, 
hence already in $G^{\langle k,m \rangle}$ we get this undesirable equation.

Now for $k <  m < \omega$ and $ x \in Y_k $ let 
$ \pi ^{k,m}_x $ be the following permutation of 
$X^{\langle k+1,m \rangle}$: it maps $y_1 \in 
X^{\langle k+1,m \rangle}$ to $y_2$ if 
$\circledast^{u,{\frak s}}_{x,y_1,y_2}$.

It is easy to check that
\mr
\item "{$\boxdot_1$}"   For $k,m,x$ as above, 
{\roster
\itemitem{ $(i)$ }  $\pi^{k,m}_x$
is a permutation of $ X^{\langle k+1,m \rangle}$ which maps 
$\Gamma^{\langle k+1,m \rangle}$ onto itself
\sn
\itemitem{ $(ii)$ }    $\pi^{k,m}_x$ induce an automorphism
$\hat{\pi}^{k,m}_x $  of  $G^{\langle k,m \rangle}$: the one
mapping $g_{y_1}$ to $ g_{y_2} $ when  $\pi^{k,m}_x (y_1)=y_2 $ 
\sn
\itemitem{ $(iii)$ }  the automorphisms $\hat{\pi}^{k,m}_x$
of $G^{\langle k,m \rangle}$ for $x \in Y_k$ pairwise commute 
\sn
\itemitem{ $(iv)$ }   the automorphism $\hat{\pi}^{k,m}_x$ 
of $G^{\langle k,m \rangle}$ is of order two. 
\endroster}
\ermn
We prove this revised formulation of the uniqueness by induction on
$m-k$.
\nl
Note that
\mr
\item "{$(*)$}"  if $x \in Y_k,y \in Y_\ell$ and $x <^* y$ then $\ell \le k$.
\ermn
If $m-k=0$, 
then $G^{<k,m>}$ is the trivial group so the uniqueness is trivial.

Also the case $ k=m-1 $ is trivial too as in this 
case $G^{\langle k,m \rangle}$ 
is actually a vector space over $ {\Bbb Z}/2{\Bbb Z}$
with basis $\{g_x:x \in Y_k\}$, well in additive notation
so the uniqueness is clear.

So assume that $m- k \ge 2$, now we need
\mr
\item "{$\boxdot^2_{k,m}$}"    if $x'_1,\dots,x'_{n_1},y'_1,
\dots, y'_{n_2 }$  from 
$X^{\langle k,m \rangle}$ are as above in $G^{<k,m>}$ then 
$\langle x'_1,\dots,x'_{n_1} \rangle = \langle y'_1,\dots,y'_{n_2}\rangle$.
\ermn
We can prove the induction step by \scite{m.8} below.
\nl
So 1),2),3) holds.
\nl
4)  Included in the proof of $\circledast_4$ inside the proof of parts
(1),(2),(3).
\nl
5) For $\alpha < \beta \le \infty$,  clearly $X^{< \alpha}_u
\subseteq X^{< \beta}_u$ and $\Gamma^{< \alpha}_u
\subseteq \Gamma^{< \beta}_u$ hence there is a homomorphism
from $G^{< \alpha}_u$ into $G^{< \beta}_u$. This
homomorphism is one-to-one (because of the uniqueness clause in
part (2)) hence the homomorphism is the identity.  So the sequence is
$\subseteq$-increasing, the $\subset$ follows by 
part (1), the uniqueness we have rk$_I(t) = \alpha < \infty \Rightarrow
g_{(\langle t \rangle,<>)} \in G^{< \alpha +1}_u \backslash G^{< \alpha}_u$.  
\nl
6),7),8),9)   Easy.   \hfill$\square_{\scite{m.7}}$ 
\enddemo
\bigskip

\demo{\stag{m.8} Observation}  Assume that
\mr
\item "{$(a)$}"  $G$ is a group
\sn
\item "{$(b)$}"  $f_t$ is an automorphism of $G$ for $t \in J$
\sn
\item "{$(c)$}"  $f_t,f_s \in \text{ Aut}(G)$ commute for any $s,t \in
J$.
\ermn
\ub{Then} there are $K$ and $\langle g_t:t \in J \rangle$ such that
\mr
\item "{$(\alpha)$}"  $K$ is a group
\sn
\item "{$(\beta)$}"  $G$ is a normal subgroup of $K$
\sn
\item "{$(\gamma)$}"  $H$ is generated by $G \cup \{g_t:t \in J\}$
\sn
\item "{$(\delta)$}"  if $a \in G$ and $t \in G$ then $g_t a g^{-1}_t
= f_t(a)$
\sn
\item "{$(\varepsilon)$}"  if $<_*$ is a linear order of $J$ then every
member of $K$ has a one and only one representation as $x
g^{b_1}_{t_1} g^{b_2}_{t_2} \ldots g^{b_n}_{t_n}$ where $x \in G,n <
\omega,t_1 <_* \ldots <_* t_n$ are from $J$ and $b_1,\dotsc,b_n \in
\Bbb Z \backslash \{0\}$.
\endroster
\enddemo
\bigskip

\demo{Proof}  A case of semi-direct product, see below.  (It is also a case of 
repeated HNN extensions).  \hfill$\square_{\scite{m.8}}$
\enddemo
\bigskip

\definition{\stag{m.9} Definition/Claim}  1) Assume $G_1,G_2$ are groups and
$\pi$ is a homomorphism from $G_1$ into Aut$(G_2)$, we define the sem-direct
product $G = G_1 *_\pi G_2$ as follows:
\mr
\item "{$(a)$}"   the set of elements is $G_1 \times G_2 =
\{(g_1,g_2):g_1 \in G_1,g_2 \in G_2\}$
\sn
\item "{$(b)$}"  the product operation is $(g_1,g_2) * (h_1,h_2) =
(g_1 h_1,g^{\pi(h_1)}_2 h_2)$ where
{\roster
\itemitem{ $(\alpha)$ }  $g^{\pi(h_1)}_2$ is the image of $g_2$ by the
automorphism $\pi(h_1)$ of $G_2$
\sn  
\itemitem{$ (\beta)$ }  $g_1 h_1$ is a $G_1$-product
\sn  
\itemitem{ $(\gamma)$ }  $g^{\pi(h_1)}_2 h_2$ is a $G_2$-product.
\endroster}
\ermn
2)
\mr
\item "{$(a)$}"  such group $G$ exists
\sn
\item "{$(b)$}"  in $G$ every member has one and only one
representation as $g'_1 g'_2$ where $g'_1 \in G_1 \times
\{e_{G_2}\},g'_2 \in \{e_{G_1}\} \times G_2$
\sn
\item "{$(c)$}"  the mapping $g_1 \mapsto (g_1,e)$ embeds $G_1$ into
$G$
\sn
\item "{$(d)$}"  the mapping $g_2 \mapsto (e,g_2)$ embeds $G_2$ into
$G$
\sn
\item "{$(e)$}"  so up to renaming, for each $h_1 \in G_1$ conjugating by
it (i.e. $g \mapsto h^{-1}_1 g h_1$) inside $G$ acts on $G_2$ as 
the automorphism $\pi(h_1)$ of $G_1$.
\ermn
3) If $H_1,H_2$ is a subgroup of $G_1,G_2$ respectively, and $g_1 \in
H_1 \Rightarrow \pi(g_1)$ maps $H_2$ onto itself and $\pi':H_1
\rightarrow \text{ Aut}(H_2)$ is $\pi'(x) = \pi(x) \restriction H_2$
\ub{then} $\{(h_1,h_2):h_1 \in H_1,h_2 \in H_2\}$ is a subgroup of
$G_1 *_\pi G_2$ and is in fact $H_1 *_{\pi'} H_2$; we denote $\pi'$ by
$\pi_{[H_1/H_2]}$.
\nl
4) If the pairs $(H^a_1,H^a_2)$ and $(H^b_1,H^b_2)$ are as in part (3)
 and $H^c_1 := H^a_1 \cap H^b_1,H^c_2 := H^a_2 \cap H^b_2$ \ub{then}
the pair $(H^c_1,H^c_2)$ is as in part (3) and $(H^a_1
*_{\pi[H^a_1,H^a_2]} H^a_2) \cap (H^b_1 *_{[H^b_1,H^b_2]} H^b_2) =
(H^c_1 *_{\pi[H^c_1,H^c_2]} H^c_2)$.
\enddefinition
\bigskip

\demo{Proof}  Known and straight.  \hfill$\square_{\scite{m.9}}$
\enddemo
\bigskip

\proclaim{\stag{m.11} Claim}  Let ${\frak s}$ be a $\kappa$-p.o.w.i.s., 
$u \in J^{\frak s}$ and
$I_u = I^{\frak s}_u$ be non-trivial.
\nl
1)  If $0 \le \alpha < \infty$ \ub{then} 
the normalizer of $G^{< \alpha}_u$ in $G_u$ is 
$G^{< \alpha +1}_u$. \nl
2) If $\alpha = { \text{\rm rk\/}}(I_u)$ \ub{then} the normalizer 
of $G^{< \alpha}_u$ in $G_u$ is 
$G^{< \infty}_u = G^{< \alpha}_u$.
\endproclaim
\bigskip

\demo{Proof}  1) First
\mr
\item "{$(*)_1$}"   if $x \in X_u$ and rk$^2_u(x) = \alpha$ then 
conjugation by $g_x$ in $G_u$ maps $\{g_y:y \in X^{< \alpha}_u\} = \{g_y:
y \in X_u$ and rk$^2_u(y) < \alpha\}$  onto itself.
\ermn
[Why?  As $g_x = g^{-1}_x$ it is enough to prove for 
every $y \in X^{< \alpha}_u$
that: $g_x g_y g^{-1}_x \in X^{< \alpha}_u$.  Now 
for each such $y$, one of the following cases occurs.
\bn
\ub{Case $(i)$}:  $g_x,g_y$ commutes so $g_x g_y g^{-1}_x = g_y \in
X^{< \alpha}_u$.

In this case the desired conclusion holds trivially.
\mn
\ub{Case $(ii)$}:  $n(y) \le n(x)$ and not case (i).

As case (i) does not occur, necessarily $n(y) < n(x)$ and $y=x
\upharpoonleft n(y)$ by \scite{m.6}(6).  Also it
follows that $t^x_{n(x)} <_{I_u[{\frak s}]} t^y_{n(y)}$, i.e., $t(x)
<_{I_u[{\frak s}]} t(y)$ but rk$^2_u(x) = \alpha$ hence rk$^2(y)
\in \{-1\} \cup [\alpha +1,\infty]$.  However we are assuming 
$y \in X^{< \alpha}_u$ hence necessarily $y \in X^{< 0}_u$, 
so $\langle \eta^y(\ell):\ell < n(y)\rangle$ is not constantly 1
hence $\langle \eta^x(\ell):\ell < n(x)\rangle$ is not constantly 1
hence rk$^2_u(x)=0$, contradiction.
\mn
\ub{Case $(iii)$}:  $n(y) > n(x)$ and not case (i).

As in case (ii) by \scite{m.6}(6) we have $x = y \upharpoonleft n(x)$.

Clearly $t(y) = t^y_{n(y)} <_{I_u[{\frak s}]} t^y_{n(x)} = t^x_{n(x)}
= t(x)$ so as rk$^2_u(x) \ge 0$ necessarily
rk$^1_u(x) = \text{ rk}^2_u(t(x)) = \alpha \in [0,\infty)$
hence rk$_{I_u}(t^y) < \text{ rk}_{I_u}(t^x) = \alpha$ and so
rk$^2_u(y) \le \text{ rk}^1_u(t^y) < \alpha$.

Let $y_1=y$ and by \scite{m.6}(1),(6) and Definition \scite{m.4}(1A)
there is $y_2$ such that $\circledast^{u,{\frak
s}}_{x,y_1,y_2}$ hence $G_u \models g_x g_y g^{-1}_x = g_{y_2}$ and
$\bar t^y = \bar t^{y_1} = \bar t^{y_2}$, so rk$^2_u(y_2) \le
\text{ rk}^1_u(y_2) = \text{ rk}^1_u(t^{y_2}) = 
\text{ rk}^1_u(t^{y_1}) < \alpha$ hence $y_2 \in X^{< \alpha}_u$ 
and so $g_{y_2} \in G^{< \alpha}_u$ so we are done.
\nl
So $(*)_1$ holds.]
\mn
Now by $(*)_1$ it follows that $g_x$ normalize $G^{< \alpha}_u$  
for every member $g_x$ of $\{g_x:\text{rk}^2_u(x) = \alpha\}$, hence clearly  
nor$_{G_u}(G^{< \alpha}_u) \supseteq (G^{< \alpha}_u) \cup 
\{g_x:\text{rk}^2_u(x) =  \alpha$ and $x \in X_u\}$
but the latter generates $G^{< \alpha +1}_u$ hence 
\mr
\item "{$(*)_2$}"  nor$_{G_u}
(G^{< \alpha}_u) \supseteq  G^{< \alpha +1}_u$.
\ermn
Second assume $g \in G_u \backslash G^{< \alpha +1}_u$, 
let $<^*$ be a linear ordering of $X_u$ as in $\boxdot$ of \scite{m.7}; so 
we can find $k < \omega$ and $x_1 <^* \ldots 
<^* x_k$ from $X_u$ such that $g = g_{x_1} g_{x_2} \ldots g_{x_k}$ and
so it suffices to prove by induction on $k$ that if $g = g_{x_1}
\ldots g_{x_k} \in G_u \backslash G^{< \alpha +1}_u$ then $g \notin
\text{ nor}_{G_u}(G^{< \alpha})$.  By \scite{m.7}(2),(4) \wilog \, $x_1
<^* \ldots <^* x_k$.  
As $g \notin G^{< \alpha +1}_u$ necessarily not all the $x_m$'s 
are from $X^{< \alpha +1}_u$ hence for some
$m,g_{x_m} \notin G^{< \alpha +1}_u$.
\mr
\item "{$(*)_3$}"  \wilog \, $x_1,x_k \notin G^{< \alpha +1}_u$.
\nl  
[Why?  So assume $x_k \in G^{< \alpha +1}_u$ hence
{\roster
\itemitem{ $(a)$ }  $x_k \in \text{\rm nor}_{G_u}(G^{< \alpha}_u)$ (as we
have already proved $G^{<(\alpha +1)}_u \subseteq 
\text{\rm nor}_{G_u}(G^{< \alpha}_u))$
\sn 
\itemitem{ $(b)$ }   nor$_{G_u}(G^{< \alpha}_u)$ is a subgroup of $G_u$
hence
\sn
\itemitem{ $(c)$ }    $g = g_{x_1} \ldots g_{x_{k-1}} g_{x_k} 
\in \text{ nor}_{G_u}(G^{< \alpha}_u)$ iff 
$g_{x_1} \ldots g_{x_{k-1}} \in \text{ nor}_{G_u}(G^{< \alpha}_u)$. 
\endroster}
By the induction hypothesis on $k$ we get are done.
Similarly if $g_{x_1} \in G^{< \alpha +1}_u$ then derive 
$g \in \text{ nor}_{G_u}(G^{< \alpha}_u)$ 
iff $g_{x_2} \ldots g_{x_k} \in \text{ nor}_{G_u}(G^{< \alpha}_u)$ to finish.]
\ermn
Now we can find $t^* \in I_u$ such that
\mr 
\item "{$(*)_4$}"   $(a) \quad t^* <_{I_u} t(x_1)$
\sn
\item "{${{}}$}"    $(b) \quad$ rk$_{I_u}(t^*) \ge \alpha$
\sn
\item "{${{}}$}"    $(c) \quad t^* \notin \{t_\ell(x):x \in
\{x_1,\dotsc,x_k\}$ and $\ell \in \{0,\dotsc,n(x)\}\}$.
\ermn
[Why?  As we assume that ${\frak s}$ is nice which implies that each
$I_u$ is non-trivial, see Definition \scite{m.1}(6) and Definition
\scite{m.2}(7).] 

Let $m(*)$ be maximal such that $1 \le m(*) \le k$ and
$(\exists i)(x_{m(*)} = x_1 \upharpoonleft i)$.
\nl
Now we choose $y \in X^{\frak s}_u$ as follows:
\mr
\item "{$(*)_5$}"  $(a) \quad \bar t^y = \bar t^{x_{m(*)}} \char 94
\langle t^* \rangle$
\sn
\item "{${{}}$}"  $(b) \quad \eta^y \restriction n(x_{m(*)}) = \eta^{x_{m(*)}}$
\sn
\item "{${{}}$}"  $(c) \quad \eta^y (n(x_{m(*)})) = 0$.
\ermn
Note that
\mr
\item "{$(*)_6$}"  $y \in X^{< 0}_u$ and $n(y) = n(x_{m(*)}) + 1$ and
\sn
\item "{$(*)_7$}"   $n(x_1) \ge \ldots \ge n(x_{m(*)}) 
\ge n(x_{m(*)+1}) \ge \ldots \ge n(x_k)$.
\ermn
We now try to define $\langle y_\ell:\ell=1,\dotsc,k+1 \rangle$ by
induction on $\ell$ as follows : 
\mr
\item "{$(*)_8$}"  $y_1 = y$ and $G_u \models g^{-1}_{x_\ell}
g_{y_\ell} g_{x_\ell} = g_{y_{\ell +1}}$ if well defined.
\ermn
So
\mr
\item "{$(*)_9$}"  $y_\ell = y$ for $\ell = 1,\dotsc,m(*)$ and is well defined.
\nl
[Why?  We prove it by induction on $\ell$.  For $\ell=1$ this is
given.  So assume that this holds for $\ell$ and we shall prove it for
$\ell+1$ when $\ell + 1 \le m(*)$.  Now $\neg(\bar t^y \triangleleft
\bar t^{x_\ell})$ by the choice of $t^*$ (and $y$) and hence $\neg(y = x_\ell
\upharpoonleft n(y) \wedge n(y) < n(x_\ell))$ and we also have
$\neg(x_\ell = y \upharpoonleft n(x_\ell) \wedge n(x_\ell) < n(y))$
as otherwise $x_\ell = x_{m(*)} \upharpoonleft n(x_\ell)$ 
but $n(x_\ell) \ge n(x_{m(*)})$ as $x_\ell
<^* x_{m(*)}$ hence
$x_\ell = x_{m(*)}$, but $\ell \ne m$ hence $x_\ell \ne x_{m(*)}$,
contradiction. 
Together by \scite{m.6}(6) the elements $g_y,g_{x_\ell}$ commute so as
by the induction hypothesis $y_\ell = y$ it follows that $y_{\ell +1}
= y$ so we are done.]
\ermn
Now:
\mr
\item "{$(*)_{10}$}"  $y_{m(*)+1}$ is well defined and 
satisfies $(*)_5(a),(b)$ and $(*)_5(c)$ when we replace 0 by 1.
\nl
[Why?  By the definition of $G_u$ in \scite{m.4}(1),(1B).]
\sn
\item "{$(*)_{11}$}"  $y_{m(*)+1} \notin G^{< \alpha}_u$.
\nl
[Why?  By $(*)_3,x_1 \notin G^{< \alpha +1}_u$ hence $\eta^{x_1}$ is
constantly 1; but $x_{m(*)} = x_1 \upharpoonleft n(x_{m(*)})$ 
hence $\eta^{x_{m(*)}}$ is constantly
one.  Now $\eta^{y_{m(*)+1}} = \eta^{x_{m(*)}} \char 94 \langle 1 \rangle$ by
$(*)_{10}$ hence $\eta^{y_{m(*)+1}}$ is constantly one.  
So rk$^2_u(y_{m(*)+1}) = \text{ rk}_{I[u]}(t^{y_{m(*)+1}}) 
= \text{ rk}_u(t^*) \ge \alpha$ so we are done.] 
\sn
\item "{$(*)_{12}$}"  if $\ell \in \{m(*)+1,\dotsc,k+1\}$ then $y_\ell =
y_{m(*)+1}$ and $y_\ell$ is well defined.
\nl
[Why?  We prove this by induction on $\ell$.  For $\ell = m(*)+1$ this
is trivial by $(*)_{10}$.  
For $\ell + 1 \in \{m(*)+2,\dotsc,k+1\}$, it is enough to
prove that $y_{m(*)+1},x_\ell$ commute.  Now $\neg(\bar t^{y_{m(*)+1}}
\triangleleft \bar t^{x_\ell})$ because $\ell g(\bar t^{y_{m(*)+1}}) =
\ell g(\bar t^y) = \ell g(\bar t^{x_{m(*)}}) + 1 \ge \ell
g(t^{x_\ell})+1 > \ell g(\bar t^{x_\ell})$ hence
$\neg \bigl( y_{m(*)+1} = x_\ell \upharpoonleft n(y_{m(*)+1}) \wedge
n(y_{m(*)+1}) < n(x_\ell) \bigr)$; also $\neg \bigl( x_\ell = y_{m(*)+1}
\upharpoonleft n(x_\ell) \wedge n(x_\ell) < n(y_{m(*)+1}) \bigr)$ 
as otherwise this contradicts the choice of $m(*)$.  So
by \scite{m.6}(6) they commute indeed.]
\sn
\item "{$(*)_{13}$}"  $g^{-1} g_y g = g_{y_{k+1}}$.
\nl
[Why?  We can prove by induction on $\ell = 1,\dotsc,k+1$ that
\nl
$(g_1 \ldots g_{\ell -1})^{-1} g_y (g_1 \ldots g_{\ell -1}) =
g_{y_\ell}$, by the definition of the $y_\ell$'s, i.e., by $(*)_8$ and
they are well defined by $(*)_9 + (*)_{10} + (*)_{12}$.]
\sn
\item "{$(*)_{14}$}"  $g^{-1} g_y g = g_{m(*)+1}$.
\nl
[Why?  By $(*)_{12}$ and $(*)_{13}$.]
\sn
\item "{$(*)_{15}$}"  $g^{-1} g_y g \notin G^{< \alpha}_u$.
\nl
[Why?  By $(*)_{14} + (*)_{11}$.]
\ermn
So by $(*)_6$ we have $g_y \in G^{< 0}_u \subseteq G^{< \alpha}_u$
and by $(*)_{15}$ we have $g^{-1} g_y g \notin G^{< \alpha}_u$ hence 
$g$ does not normalize $G^{< \alpha}_u$, so we have carried the
induction on $k$.  As $g$ was any
member of $G_u \backslash G^{<(\alpha +1)}_u$ we get nor$_{G_u}(G^{<
\alpha}_u) \subseteq G^{<(\alpha +1)}_u$.

Together with $(*)_2$ we are done.  
\nl
2) Follows.  \hfill$\square_{\scite{m.11}}$
\enddemo
\newpage

\head {\S2 Easier group} \endhead  \resetall \sectno=2
 \spuriousreset
\bigskip

The $G^{\frak s}_u$'s from \S1 has long towers of normalizers but the ``base",
$G^{< 0,{\frak s}}_u$ is in general of large cardinality.  Hence we
replace below $G^{\frak s}_u$ by $K^{\frak s}_u$ and 
$G^{< 0,{\frak s}}_u$ by $H^{\frak s}_u$. 

\definition{\stag{n.1} Definition}  Let ${\frak s}$ be a 
$\kappa$-p.o.w.i.s.
\nl
1) For $u \in J^{\frak s}$:
\mr
\item "{$(a)$}"  recall \scite{m.7}(6): ${\Cal A}_u = {\Cal A}^{\frak s}_u
= \{g G^{< 0}_u:g \in G_u\}$
is a partition of $G$ (to left cosets of $G^{< 0}_u$ inside $G_u$);
\sn
\item "{$(b)$}"  we define for every $f \in G_u$ a permutation
$\partial_f$ of ${\Cal A}_u$ defined by $\partial_f(g_1 G^{<0}_u) = (f g_1)
G^{<0}_u$, we may write it also as $f(g_1 G)$
\sn
\item "{$(c)$}"  let $L_u = L^{\frak s}_u$ be the group 
generated by $\{h_{\bold a}:\bold a \in {\Cal A}_u\}$ freely except
$h_{\bold a} h_{\bold b} = h_{\bold b} h_{\bold a}$ and $h^{-1}_{\bold a} =
h_{\bold a}$ for ${\bold a},{\bold b} \in {\Cal A}_u$; for $g \in G_u$
let $h_g = h_{g G^{<0}_u}$
\sn
\item "{$(d)$}"  let ${\bold h}_u = {\bold h}^{\frak s}_u$ be the
homomorphism from $G_u$ into the automorphism group of $L_u$ such
that $f \in G_u \wedge {\bold a} \in {\Cal A}_u \Rightarrow ({\bold
h}_u(f))(h_{\bold a}) = h_{f{\bold a}}$
\sn
\item "{$(e)$}"  let $K_u = K^{\frak s}_u$ be $G_u *_{\bold h_u} L_u$,
the twisted product of
$G_u,L_u$ with respect to the homomorphism 
${\bold h}_u$, see \scite{m.9}, and we
identify $G_u$ with $G_u \times \{e_{L_u}\}$ and $L_u$ with
$\{e_{G_u}\} \times L_u$
\sn
\item "{$(f)$}"  let $H_u = \{(e_{G_u},h_{e G^{\le
0}_u}),(e_{G_u},e_{L_u})\}$, a subgroup of $K_u$ and let $h_* :=
h_{e_{G_u}} = h_{e_{G_u} G^{\le 0}_u} \in L_u$, i.e. the pair
$(e_{G_u},g_{e G^{\le 0}_u})$, this is the unique member of $H_u$
which is not the unit.
\ermn
2) For $\alpha \le \infty$ let $K^{< \alpha}_u = 
K^{< \alpha,{\frak s}}_u$ be the subgroup $\{(g,h):g \in G^{<
\alpha}_u$ and $h \in L_u\}$ of $K_u$.  Similarly $K^{\le \alpha}_u =
K^{\le \alpha,{\frak s}}_u$. 
\nl
3) For $u \in J^{\frak s}$ let
\mr
\item "{$(a)$}"   $D_u = D^{\frak s}_u = \{(v,g):v
\le_{J[{\frak s}]} u$ and $g \in K^{\frak s}_v\}$ 
\sn
\item "{$(b)$}"   $Z^0_u = \{(\bar t,\eta):\bar t = \langle t_\ell:
\ell \le n\rangle,n < \omega,t_\ell \in I$ for each $\ell < n$ 
and $\eta \in {}^n 2\}$ and let $z = (\bar t^z,\eta^z) = (\langle
t^z_\ell:\ell \le n \rangle,\eta^z)$ and $n(z) = n$ for $z \in Z^0_u$;
this is compatible with Definition \scite{m.2}(4); note that here
$\bar t$ is not necessarily decreasing
\sn
\item "{$(c)$}"  $Z^1_u := \{\langle x_\ell:\ell < k \rangle:k <
\omega$, each $x_\ell$ is from $Z^0_u\}$ and let $z = (\langle
x^z_\ell:\ell < k(z)\rangle)$ if $z \in Z^1_u$
\sn
\item "{$(d)$}"  $Z_u := Z^0_u \cup Z^1_u$
\sn
\item "{$(e)$}"   for $z \in Z_u$ we define his$(z)$, a finite subset
of $I_u$ by
{\roster
\itemitem{ $(\alpha)$ }  if $z = (\langle t_\ell:\ell \le
n\rangle,\eta) \in Z^0_u$ then his$(z) = \{t_\ell:\ell \le n\}$
\sn
\itemitem{ $(\beta)$ }  if $z \in Z^1_u$ say $z = 
\langle (\langle t^k_\ell:\ell \le
\ell_k\rangle,\eta^k):k<k^*\rangle \in Z^1_u$ then his$(z) 
= \{t^k_\ell:k < k^*$ and $\ell \le \ell_k\}$
\endroster}
\item "{$(f)$}"  for $z = Z_u$ let $n(z) = \Sigma\{\ell_k:k < k^*\}$ if $z =
\langle (\langle t^k_\ell:\ell \le \ell_k \rangle,\eta^k):k
<k^*\rangle \in Z^1_u$ and $n(z)$ is already defined if $z \in Z^0_u$
in clause (b).
\endroster
\enddefinition
\bigskip

\demo{\stag{n.1.11} Observation}  In Definition \scite{n.1}.
\nl
1) For $u \in J^{\frak s},K_u$ is well defined and 
$G_u,L_u$ are subgroups of  $K_u$ (after the identification).
\nl
2) For $I \subseteq I^{\frak s}_u$ let $L^{\frak s}_{u,I}$ be the
subgroup of $L^{\frak s}_u$ be generated by $\{h_{g G^{<0}_u}:g \in
G^{\frak s}_{u,X_I}\}$.  If $I_1,I_2 \subseteq I^{\frak s}_u$ then
$L^{\frak s}_{u,I_1} \cap L^{\frak s}_{u,I_2} = L^{\frak s}_{u,I_1
\cap I_2}$.
\nl
3) For $I \subseteq I^{\frak s}_u$ let $K^{\frak s}_{u,I}$ be the
subgroup of $K^{\frak s}_u$ generated by $G^{\frak s}_{u,X_I} \cup
L^{\frak s}_{u,I}$.  Then
\mr
\item "{$(a)$}"   $G^{\frak s}_{u,X_I}$ normalized 
$L^{\frak s}_{u,I}$ inside $K^{\frak s}_u$
\sn
\item "{$(b)$}"  $K^{\frak s}_{u,I}$ is $G^{\frak s}_{u,X_I} *_\pi
L^{\frak s}_{u,I}$ for the natual $\pi$.
\ermn
Also
\mr
\item "{$(c)$}"   if 
$I_1,I_2 \subseteq I^{\frak s}_u$ then $K^{\frak s}_{u,I_1}
\cap K^{\frak s}_{u,I_2} = K^{\frak s}_{u,I_1 \cap I_2}$.
\endroster
\enddemo
\bigskip

\demo{Proof}  Easy (recall \scite{m.7}(8),(9), \scite{m.9}(2),(3)).
\enddemo
\bigskip

\definition{\stag{n.1.15} Definition}  1) If $I$ is a partial order
\ub{then} ${}^k I$ is the set of $\bar t = \langle t_\ell:\ell < k
\rangle$ where $t_\ell \in I$.
\nl
2) If $\bar t \in {}^k I$ then tp$_{\text{qf}}(\bar t,\emptyset,I) =
\{(\iota,\ell_1,\ell_2):\iota =0$ and $I \models t_{\ell_1} <
t_{\ell_2}$ or $\iota=1$ and $t_{\ell_1} = t_{\ell_2}$ or $\iota
=2$ and $I \models t_{\ell_1} > t_{\ell_2}$ and $\iota = 3$ if none of
the previous cases$\}$.
\nl
2A) Let ${\Cal S}^k = \{\text{tp}_{\text{qf}}(\bar t,
\emptyset,I):\bar t \in {}^k I$ and $I$ is a partial order$\}$.
\nl
3) We say $\bar t \in {}^k I$ realizes $p \in {\Cal S}^k$ when $p =
\text{ tp}_{\text{qf}}(\bar t,\emptyset,I)$.
\nl
4) If $k_1 < k_2$ and $p_2 \in {\Cal S}^{k_2}$ then $p_1:= p_2
\restriction k_1$ is the unique $p_1 \in {\Cal S}^{k_1}$ such that
if $p_2 = \text{ tp}_{\text{qf}}(\bar t,\emptyset,I)$ then $p_1 =
\text{ tp}_{\text{qf}}(\bar t \restriction k_1,\emptyset,I)$.
\enddefinition
\bigskip

\remark{Remark}  Below each member of
$\Lambda^0_k,\Lambda^1_k,\Lambda^2_k$ will be a description of an
element of $G^{\frak s}_u,{\Cal A}^{\frak s}_u,K^{\frak s}_u$
respectively from a $k$-tuple of members of $I^{\frak s}_u$.  Of
course, a member of $Z^{\frak s}_u$ is a description of a generator of
$K^{\frak s}_u$. 
\endremark
\bigskip

\definition{\stag{n.1.17} Definition}  1) For $k < \omega$ 
let $\Lambda^0_k = \cup\{\Lambda^0_{k,p}:p \in {\Cal
S}^k\}$ where for $p \in {\Cal S}^k$ we let
$\Lambda^0_{k,p}$ be the set of sequences of the form
$\langle (\bar \ell_j,\eta_j):j < j(*)\rangle$ such that:
\mr
\item "{$(a)$}"  for each $j$ for some $n=n(\bar \ell_j,\eta_j)$ we
have $\bar \ell_j = \langle \ell_{j,i}:i \le n(\bar \ell_j,\eta_j)\rangle$
is a sequence of numbers $<k$ of length $n+1$ such that $p =
\text{ tp}_{\text{qf}}(\bar t,\emptyset,I) \Rightarrow \langle
t_{\ell_{j,i}}:i \le n(\ell_j,\eta_j)\rangle$ is decreasing
\sn
\item "{$(b)$}"  for each $j,\eta_j \in {}^n 2$ where $n = n(\bar
\ell_j,\eta_j)$.
\ermn
2) For any p.o.w.i.s.
${\frak s},u \in J^{\frak s},\bar t \in {}^k(I_u)$ and $\rho =
\langle(\bar \ell_j,\eta_j):j<j(*)\rangle \in \Lambda^0_k$, let
$g^u_{\bar t,\rho} = g^{u,{\frak s}}_{\bar t,\rho} =
(\ldots g_{(\bar t^j,\eta_j)} \ldots)_{j<j(*)}$, the product in $G_u
\subseteq K_u$ (so if $j(*)=0$ it is $e_{G_u} = e_{K_u}$) where
\mr
\item "{$(a)$}"  $\bar t^j = \text{ seq}_{\rho,j}(\bar t) := \langle
t_{\ell_{j,i}}:i \le n(\ell_j,\eta_j)\rangle$
\sn
\item "{$(b)$}"  if $\bar t^j$ is decreasing (in $I_u$) then $g_{(\bar
t^j,\eta_j)} \in G_u \subseteq K_u$ is already well defined, if not then
$g_{(\bar t^j,\eta_j)} = e_{K_u}$.
\ermn
2A) For a p.o.w.i.s. ${\frak s},u \in J^{\frak s},t \in {}^k(I^{\frak
s}_u)$ and $\rho = \langle(\bar \ell_j,\eta_j):j < j(*)\rangle \in
\Lambda^0_k$ let $z^u_{\bar t,\rho} = z^{u,{\frak s}}_{\bar t,\rho}$
be the following member of $Z^{1,{\frak s}}_u$: it is $\langle x_{\bar
t,\rho,j}:j < j(*)\rangle$ where $x_{\bar t,\rho,j} = x_{\bar t,(\bar
\ell_j,\eta_j)} = (\langle t_{\ell_{j,i}}:i \le n(\bar
\ell_j,\eta_j)\rangle,\eta_j)$. 
For $p \in {\Cal S}^k$ and $\rho = \langle(\bar \ell_j,\eta_j):j <
j(*)\rangle \in \Lambda^0_{k,p}$ let supp$(\rho) =
\cup\{\text{Rang}(\bar \ell_j):j<j(*)\}$ and if $\bar t \in
{}^k(I^{\frak s}_u)$ let sup$(\bar t,\rho) = \{t_\ell:\ell \in \text{
supp}(\rho)\}$.
\nl
2C) We say $\rho \in \Lambda^0_k$ is $p$-reduced when: $p \in {\Cal
S}^k$ and for every p.o.w.i.s. ${\frak s},u \in J^{\frak s}$ and $t
\in {}^k(I^{\frak s}_u)$ realizing $p$ (in $I^{\frak s}_u$), for no
$\rho' \in \Lambda^0_k$ do we have supp$(\rho') \subset \text{ supp}(\rho)$ and
$g^{u,{\frak s}}_{\bar t,\rho'} = g^{u,{\frak s}}_{\bar t,\rho'}$.
\nl
2D) We say that $\rho \in \Lambda^0_k$ is explicitly $p$-reduced when
the sequence is with no repetitions and $\langle
n(\bar \ell_j,\eta_j):j < j(*)\rangle$ is non-increasing (the
length can be zero).
\nl
3) For $k < \omega$ let $\Lambda^1_k =
\cup\{\Lambda^1_{k,p}:p \in {\Cal S}^k\}$ where for $p \in {\Cal S}^k$
we let $\Lambda^1_{k,p}$ be 
the set of $\rho = \langle (\bar \ell_j,\eta_j):j<j(*)\rangle
\in \Lambda^0_{k,p}$ such that: for every ${\frak s}$ and $u \in 
J^{\frak s}$ if 
$\bar t \in {}^k(I^{\frak s}_u)$ realizes $p$ \ub{then} there is no
$\rho' \in \Lambda^0_{k,p}$ with supp$(\rho') \subset \text{\rm
supp}(\rho)$ and satisfying
$g^{u,{\frak s}}_{\bar t,\rho} G^{< 0}_u = g_{\bar t,\rho'} 
G^{< 0}_u$.
\nl
4) For $k < \omega$ and $p \in {\Cal S}^k$ let $\Lambda^2_{k,p}$
be the set of finite sequences $\varrho$ of 
length $\ge 1$ such that $\varrho(0)
\in \Lambda^0_{k,p}$ and $0 < i < \ell g(\varrho) \Rightarrow \ell
g(\varrho(i)) > 0 \wedge \varrho(i) \in \Lambda^1_{k,p}$.
Let $\Lambda^2_k = \cup\{\Lambda^2_{k,p}:p \in {\Cal S}^k\}$.
\nl
5) For any ${\frak s}$, if 
$u \in J^{\frak s},\bar t \in {}^k(I_u)$ and $\varrho = \langle
\rho_i:i <i(*)\rangle \in \Lambda^2_k$ then $g_{\bar t,\varrho} \in
K_u$ (recalling $i(*) \ge 1$) is $g_{\bar t,\rho_0} 
h_{g_{\bar t,\rho_1}} h_{g_{\bar t,\rho_2}},
\dotsc,h_{g_{\bar t,\rho_{i(*)-1}}}$ (product in $K_u$) where $g_{\bar
t,\rho_\ell}$ is from Part (2), recalling that $h_g = h_{g G^{<0}_u}$
is from clause (c) of Definition \scite{n.1}(2).
\nl
5A) For any p.o.w.i.s. ${\frak s},u \in J^{\frak s},\bar t \in
{}^k(I^{\frak s}_u)$ and $\varrho = \langle \rho_i:u < i(*)\rangle
\in \Lambda^2_k$, let $z^u_{\bar t,\varrho} = z^{u,{\frak s}}_{\bar t,\varrho}
\in Z^{2,{\frak s}}_u$ be $\langle z^u_{\bar t,\rho_i}:i < i(*)\rangle$.
\nl
5B) For $p \in {\Cal S}^k$ and $\varrho \in \Lambda^2_{k,p}$ let 
supp$(\varrho) = \cup\{\text{supp}(\varrho(i):i < \ell g(\varrho)\}$.
\nl
5C) We say $\varrho \in \Lambda^2_{k,p}$ is $p$-reduced when for every
p.o.w.i.s. ${\frak s},u \in J^{\frak s}$ and $\bar t \in {}^k(I^{\frak
s}_u)$ realizing $p$, for no $\varrho' \in \Lambda^2_{k,p}$ do we have
(in $K^{\frak s}_u$) $g^{u,{\frak s}}_{\bar t,\varrho'} = 
g^{u,{\frak s}}_{\bar t,\varrho}$ and supp$(\varrho') 
\subset \text{ supp}(\varrho)$.
\enddefinition
\bigskip

\definition{\stag{n.1.18} Definition}   1) For 
$\rho_1,\rho_2 \in \Lambda^0_{k,p}$ we say $\rho_1 
{\Cal E}^0_{k,p} \rho_2$ or $\rho_1,\rho_2$ are 0-p-equivalent when:
for every p.o.w.i.s. ${\frak s}$ and $u \in J^{\frak s}$ and $\bar t
\in {}^k(I^{\frak s}_u)$ realizing $p$ the elements $g^{u,{\frak
s}}_{\bar t,\rho_1},g^{u,{\frak s}}_{\bar t,\rho_2}$ of $G^{\frak
s}_u$ are equal.
\nl
2) For $\rho_1,\rho_2 \in \Lambda^1_{k,p}$ we say $\rho_1 
{\Cal E}^2_{k,p} \rho_2$ or $\rho_1,\rho_2$ are 1-$p$-equivalent when: for
every p.o.w.i.s. ${\frak s}$ and $u \in J^{\frak s}$ and $\bar t \in {}^k(I_u)$
realizing $p$ we have $g^{u,{\frak s}}_{\bar t,\rho_1} G^{< 0}_u =
g^{u,{\frak s}}_{\bar t,\rho_2} G^{< 0}_u$.
\nl
3) For $\varrho_1,\varrho_2 \in \Lambda^2_{k,p}$ we say that
$\varrho_1 {\Cal E}_{2,p} \varrho_2$ or $\varrho_1,\varrho_2$ are 
2-$p$-equivalent, \ub{when}: for every p.o.w.i.s. ${\frak s}$ and $u
\in J^{\frak s}$ and $\bar t \in {}^k(I_u)$ realizing $p$ the element
$g^{u,{\frak s}}_{\bar t,\rho_1}$ and $g^{u,{\frak s}}_{\bar
t,\rho_2}$ of $K^{\frak s}_u$ are equal. 
\enddefinition
\bigskip

\proclaim{\stag{n.1.19} Claim}  1) In Definition \scite{n.1.17} parts
(2C),(3),(5B) saying ``for every p.o.w.i.s.
${\frak s},u \in J^{\frak s}$ and $\bar t \in {}^k(I_u)$
realizing $p$" it is equivalent to saying ``for some ...".
\nl
2) In Definition \scite{n.1.17}, $E^\iota_{k,p}$ is an
equivalence relation on $\Lambda^\iota_{k,p}$ for $\iota = 0,1,2$.
Every $E^\iota_{k,p}$-equivalence class contains a reduced member and
for $\iota=0$ even an explicitly reduced one.  Explicitly reduced
implies reduced.
\nl
3) For every p.o.w.i.s. ${\frak s}$, \ub{if} $u \in J^{\frak s}$ and
$\bar t \in {}^k(I^{\frak s}_u)$ realizes $p \in {\Cal S}^k$ then
\mr
\item "{$(a)$}"  for $\rho_1,\rho_2 \in \Lambda^0_{k,p}$ we have
{\roster
\itemitem{ $(\alpha)$ }  $g^{u,{\frak s}}_{\bar t,\rho_1} = 
g^{u,{\frak s}}_{\bar t,\rho_2}$ iff $\rho_1 {\Cal E}^0_{k,p} \rho_2$
\sn
\itemitem{ $(\beta)$ }  $\{\ell^{\rho_1}_{j,i}:j < \ell g(\rho_1),i \le 
n(\bar \ell^{\rho_1}_j,\eta^{\rho_1}_j)\} =
\{\ell^{\rho_2}_{j,i}:j < \ell g(\rho_1),i \le 
n(\bar \ell^{\rho_1}_j,\eta^{\rho_2}_j)\}$
\sn
\itemitem{ $(\gamma)$ }  if $\rho_1,\rho_2$ are explicitly
$p$-reduced, \ub{then} they are $\rho_1 {\Cal E}^0_{k,p} \rho_2$ \ub{iff}
letting $\rho_i = 
\langle(\bar \ell^i_j,\eta^i_j):j<j_i\rangle$ for $i=1,2$ we have
\sn
\itemitem{ ${{}}$ }  $\qquad (a) \quad j_1 = j_2$
\sn
\itemitem{ ${{}}$ }  $\qquad (b) \quad$ for some 
permutation $\pi$ of $\{0,\dotsc,j_1-1\}$ we
have 
\nl

\hskip50pt $(\bar \ell^2_j,\eta^2_j) = (\bar \ell^1_{\pi(j)},\eta^2_{\pi(j)})$
(so actually only the domain of ${\Cal E}_{0,p}$
\nl

\hskip50pt  depends on $p$).
\endroster}
\item "{$(b)$}"  for $\rho_1,\rho_2 \in \Lambda^1_{k,p}$ we have
{\roster
\itemitem{ $(\alpha)$ }  $g^{u,{\frak s}}_{\bar t,\rho_1} G^{< 0}_u =
g^{u,{\frak s}}_{\bar t,\rho_2} G^{< 0}_u$ iff $\rho_1 
{\Cal E}^1_{k,p} \rho_2$
\endroster}
\ermn
4) For every p.o.w.i.s. ${\frak s}$ if $u \in J^{\frak s}$ and 
$\bar \ell \in {}^{\bar k}(I^{\frak s}_u)$ realizes $p \in {\Cal S}^k$ then
\mr
\item "{$(c)$}"  for $\varrho_1,\varrho_2 \in \Lambda^2_{k,p}$ we have
{\roster
\itemitem{ $(\alpha)$ }  $g^{u,{\frak s}}_{\bar t,\varrho_1} 
= g^{u,{\frak s}}_{\bar t,\varrho_2}$ iff 
$\varrho_1 {\Cal E}^2_{k,p} \varrho_2$
\sn
\itemitem{ $(\beta)$ }  if $\varrho_1 {\Cal E}^2_{k,p} \varrho_2$ and
$\varrho_1,\varrho_2$ are $p$-reduced then
{\rm supp}$(\varrho_1) = \,\text{\rm supp}(\varrho_2)$.
\endroster}
\endroster
\endproclaim
\bigskip

\demo{Proof}  Straight, (recalling \scite{m.7}(7) and note that (3)
elaborate (1)).  \hfill$\square_{\scite{n.1.19}}$ 
\enddemo
\bigskip

\proclaim{\stag{n.1.23} Claim}  Assume $k < \omega,p \in {\Cal
S}^k,{\frak s}$ is a p.o.w.i.s., $u \in J^{\frak t}$ and $\bar
t_1,\bar t_2 \in {}^k I$ satisfies $p = \text{\rm tp}(\bar
t_\ell,\emptyset,I^{\frak s}_u)$ for $\ell = 1,2$.
\nl
1) If $\rho \in \Lambda^0_{k,p}$ and $\rho$ is $p$-reduced and
$g_{\bar t_1,\rho} = g_{\bar t_2,\rho} \in G^{\frak s}_u$, \ub{then} $\bar t_2
\restriction \text{ \rm supp}(\rho)$ is a permutation of $\bar t_1
\restriction \text{\rm supp}(\rho)$. 
\nl
2) If $\rho \in \Lambda^1_{k,p}$ is $p$-reduced and
$g^{u,{\frak s}}_{\bar t_1,\rho} G^{<0}_u = g^{u,{\frak s}}_{\bar
t_2,\rho} G^{< 0}_u$ \ub{then} $\bar t_1 \restriction$ 
{\rm supp}$(\rho)$ is a permutation of $\bar t_2 \restriction$ 
{\rm supp}$(\rho)$.
\nl
3) If $\varrho \in \Lambda^2_{k,p}$ is $p$-reduced and $g^{u,{\frak s}}_{\bar
t_1,\varrho} = g^{u,{\frak s}}_{\bar t_2,\varrho}$ so both are well
defined \ub{then} similarly
$\bar t_1 \restriction \text{\rm supp}(\varrho)$ is a permutation of
$\bar t_2 \restriction \text{\rm supp}(\varrho)$ and both are with no
repetition. 
\nl
4) For every $\varrho_1 \in \Lambda^2_{k,p}$ there is a $p$-reduced
$\varrho_2$ such that for every p.o.w.i.s., $u \in J^{\frak s}$ and
$\bar t \in {}^k(I^{\frak s}_u)$ realizing $p$ we have 
$g^{u,{\frak s}}_{\bar t,\varrho_1} = g^{u,{\frak s}}_{\bar t,\varrho_2}$.
(Similarly for $\Lambda^0_{k,p},\Lambda^1_{k,p}$).
\endproclaim
\bigskip

\demo{Proof}  Straight.
\enddemo
\bigskip

\definition{\stag{nx.0} Definition}  Let ${\frak s}$ be a
$\kappa$-p.o.w.i.s.
\nl
1) For $u \le_{J[{\frak s}]} v$ let $\hat \pi^0_{u,v}$ be the following
partial mapping from $Z^{0,{\frak s}}_v$ to $Z^{0,{\frak s}}_u$,
recalling Definition \scite{n.1}(3)(b):

$x \in \text{ Dom}(\hat \pi^0_{u,v})$ iff $x \in Z^{0,{\frak s}}_v$ and
$\pi_{u,v}(t^x_\ell)$ is well defined for $\ell \le n(x)$ and then $\hat
\pi_{u,v}(x) = (\langle \pi_{u,v}(t^x_\ell):\ell \le n(x)\rangle,\eta^x)$.
\nl
2) For $u \le_{J[{\frak s}]} v$ let $\hat \pi^1_{u,v} =
\hat\pi^{1,{\frak s}}_{u,v}$ be the
following partial mapping $Z^1_v$ to $Z^1_u$: if $z \in Z^1_u$ so $z =
\langle (\bar t^k,\eta^k):k < k^*\rangle$ 
and $\bar t^k = \langle t^k_\ell:\ell
\le \ell_k\rangle,t^k_\ell \in I_v$ for $k<k^*,\ell \le \ell_k$ \ub{then} $\hat
\pi^1_{u,v}(z) = \langle(\langle\pi_{u,v}(t^k_\ell):
\ell < \ell_k \rangle,\eta^k):k<k^*\rangle$ when
each $\pi_{u,v}(t^k_\ell)$ is well defined.
\nl
3) Let $u \le_{J[{\frak s}]} v$ let $\hat \pi_{u,v}$ be 
$\hat \pi^0_{u,v} \cup \hat\pi^1_{u,v}$.  
\nl
4) For $u \in J^{\frak s}$ and $z \in Z_u$ let $\partial_{u,z}$ be the 
following permutation of $D_u = D^{\frak s}_u$ where $D_u$ is from
Definition \scite{n.1}(3)(a).
\bn
For each $(v,g) \in D_u$ we define $\partial_{u,z}((v,g))$ as follows:
\bn
\ub{Case 1}:  $z \in \text{ Dom}(\hat \pi^0_{v,u}) \subseteq Z^0_u$ 
and $\hat \pi_{v,u}(z) \in
X^{\frak s}_v$, i.e., $\langle \hat \pi_{v,u}(t^z_\ell):\ell\le n(*)
\rangle$ is $\le_{I_u}$-decreasing.

Then let $\partial_{u,z}((v,g)) = (v,g_{\hat \pi_{v,u}(z)} g)$ noting
$g_{\pi_{v,u}(z)} \in G_v \subseteq K_v$.
\bn
\ub{Case 2}:  $z \in \text{ Dom}(\hat\pi^1_{v,u}) \subseteq Z^1_u$ so $z =
\langle x_\ell:\ell < k\rangle$ and $x_\ell \in \text{
Dom}(\hat\pi^0_{v,u})$ for $\ell < k$ and let 
$x'_\ell := \hat \pi^0_{v,u}(x_\ell) \in X^{\frak s}_v$ for $\ell < k$.

Then let $\partial_{u,z}((v,g)) = (v,g')$ when $g' \in K_v$ is defined
by as $h_{g_{x'_0} \ldots g_{x'_{k-1}}} g$, product in $K_u$ 
noting $g_{x'_\ell} \in G_v \subseteq K_v$ for $\ell < k$.
\bn
\ub{Case 3}:  Neither case 1 nor case 2.

Then let $\partial_{u,x}((v,g)) = (v,g)$.
\enddefinition
\bigskip

\demo{\stag{nx.2} Observation}  In Definitions \scite{n.1}, \scite{nx.0}:
\nl
1) If $u \le_{J[{\frak s}]} v$ \ub{then} $\hat \pi_{u,v}$ is a partial
mapping from $Z_v$ to $Z_u$.
\nl
2) In part (1), $\hat \pi_{u,v}$ maps $Z^0_v,Z^1_v$ to $Z^0_u,Z^1_u$
respectively, that is it maps $Z^\ell_v \cap \text{ Dom}(\hat \pi_{u,v})$
into $Z^\ell_u$ for $\ell = 0,1$.
\nl
3) If $u \le_{J[{\frak s}]} v$ and ${\frak s}$ is nice or just
Dom$(\pi_{u,v}) = I_v$ \ub{then} Dom$(\hat \pi_{u,v}) = Z_v$.
\nl
4) nor$_{K_u}(H_u)$ is $K^{< 0}_u$ where $H_u$ is from Definition
\scite{n.1}(1)(f).
\nl
5) nor$^{1 +\alpha}_{K_u}(H_u)$ is $K^{< \alpha}_u$ for $\alpha \ge 0$
if ${\frak s}$ is non-trivial.
\enddemo
\bigskip

\demo{Proof}  1),2),3) Check.
\nl
4) As $H_u$ has two elements $e_{K_u}$ and $h_*$ clearly an
element of $K_u$ normalize $H_u$ iff it commutes with $g_*$.  Now
when does $(g,h) \in G_u *_{\bold h_u} L_u$ commute with $g_* =
(e_{G_u},h_{e_{G_u}} G^{<0}_u)$?  Note that

$$
(g,h)(e_{G_u},h_{e_{G_u}G^{\le 0}_u}) = (g,h+h_{e_{G_u}G^{< 0}_u})
$$

$$
(e_{G_u},h_{e_{G_u} G^{< 0}_u})(g,h) = (g,((\bold h_u(g))
(h_{e_{G_u} G^{< 0}_u}) + h).
$$
\mn
As $L_u$ is commutative, ``$(g,h)$ commute in $K_u$" iff in $L_u$

$$
(\bold h_u(g))(h_{e_{G_u} G^{< 0}_u}) = h_{e_{G_u} G^{< 0}_n}.
$$
\mn
By the definition of $\bold h_u \in \text{ Hom}(G_u,\text{Aut}(L_u))$
in \scite{n.1}(1)(d),(e)
this means

$$
(g e_{G_u}) G^{< 0}_u = e_{G_u} G^{< 0}_u.
$$
\mn
i.e.

$$
g \in G^{< 0}_u.
$$
\mn
We can sum that: $(g,h) \in G_u *_{\bold h_u} L_u$ belongs to
nor$_{K_u}(H_u)$ iff $(g,h)$ commutes with $h_*$ iff $g \in G^{< 0}_u$
iff $(g,h) \in K^{< 0}_u$, as required.
\nl
5) Let $\bold f_u:K_u \rightarrow G_u$ be 
defined by $\bold f_u((g,h))=g$.  Clearly
\mr
\item "{$(*)_1$}"  $\bold f_u$ is a homomorphism from $K_u$ onto $G_u$ and
for every ordinal $\alpha \ge 0$, it maps $K^{< \alpha}_u$ onto $G^{<
\alpha}_u$ so $\bold f_u(K^{< \alpha}_u)= G^{< \alpha}_u$ and moreover
$\bold f^{-1}(G^{< \alpha}_u) = K^{< \alpha}_u$
(see the definition of $K^{< \alpha}_u$ in \scite{n.1}(2)). 
\ermn
Also
\mr
\item "{$(*)_2$}"  Ker$(\bold f_u) = \{e_{G_u}\} \times L_u 
\subseteq K^{< 0}_u$.
\ermn
Now we prove by induction on the ordinal $\alpha \ge 0$ that
nor$^{1 +\alpha}_{K_u}(H_u) = K^{< \alpha}_u$.  For $\alpha = 0$ 
this holds by part (4).  For $\alpha$ limit this holds as both
$\langle\text{nor}^\beta_{K_u}(H_u):\beta \le \alpha \rangle$ and
$\langle K^{< \beta}_u:\beta \le \alpha\rangle$ are increasing
continuous.

Lastly, for $\alpha = \beta +1 > 0$ we have for any $f \in K_u$

$$
\align
f \in \text{ nor}^{1 +\alpha}_{K_u}(H_\beta) &\Leftrightarrow f \in \text{
nor}_{K_u}(\text{nor}^{1 +\beta}_{K_u}(H_\beta)) \\
  &\Leftrightarrow f \in \text{ nor}_{K_u}(\bold f^{-1}_u(G^{<
\beta}_u)) \\
  &\Leftrightarrow f(\bold f^{-1}_u(G^{< \beta}_u))f^{-1} = \bold
f^{-1}_u(G^{< \beta}_u) \\
  &\Leftrightarrow \bold f_u(f) G^{< \beta}_u \bold f_u(f)^{-1} = G^{<
\beta}_u \\
  &\Leftrightarrow \bold f_u(f) \in \text{ nor}_{G_u}(G^{<
  \beta}_u) \\
  &\Leftrightarrow \bold f_u(f) \in G^{< \alpha}_u \Leftrightarrow f \in
  K^{< \alpha}_u.
\endalign
$$
\mn
[Why?  The first $\Leftrightarrow$ by the definition of nor$^{\beta
+1}_{K_u}(-)$, the second $\Leftrightarrow$ by the induction
hypothesis, the third $\Leftrightarrow$ by the definition of
nor$_{K_u}(-)$, the fourth $\Leftrightarrow$ by $(*)_1$, the fifth
$\Leftrightarrow$ by the definition of nor$_{G_u}(-)$, the sixth
$\Leftrightarrow$ by \scite{m.11}(1), 
the seventh $\Leftrightarrow$ by $(*)_1$.]
${{}}$ \hfill$\square_{\scite{nx.2}}$
\enddemo
\bigskip

\demo{\stag{nx.3} Observation}  Let ${\frak s}$ be a p.o.w.i.s.
\nl
1) For $u \in J^{\frak s}$ and $x \in
Z^{\frak s}_u$ we have: $\partial_{u,x}$ is a
well defined function and is a permutation of $D^{\frak s}_u$.
\nl
2) If $u \le_{J[{\frak s}]} v$ \ub{then} $D^{\frak s}_u \subseteq
D^{\frak s}_v$.
\nl
3)  If $u \le_{J[{\frak s}]} v$ and $y \in Z^{\frak s}_v$ and
$x = \hat \pi_{u,v}(y)$ then $\partial_{u,x} = \partial_{v,y} 
\restriction D_u$.
\nl
4) If ${\frak s}$ is nice and $u \in J^{\frak s}$ and $z \in 
Z^{\frak s}_u$ \ub{then} in the definition \scite{nx.0}(4) of
 $\partial_{u,z}$ Case 3 never occurs.
\enddemo
\bigskip

\demo{Proof}  Straight.
\enddemo
\bigskip

\definition{\stag{nx.0.7} Definition}  Let ${\frak s}$ be
$\kappa$-p.o.w.i.s.
\nl
1) Let $\bold S^k = \{\bold q:\bold q$ is a function with domain
${\Cal S}^k$ and for $q \in {\Cal S}^k,\bold q(p) \in
\Lambda^2_{k,p}\}$, on $\Lambda^2_{k,p}$, see Definition
\scite{n.1.17}(4) above. 
\nl
2) We say that $\bold q \in \bold S^k$ is disjoint when $\langle
\text{supp}(\bold q(p)):p \in {\Cal S}^k\rangle$ is a sequence of
pairwise disjoint sets.  We say that $\bold q$ is reduced when $\bold
q(p)$ is $p$-reduced for every $p \in {\Cal S}^k$.
\nl
3) Let $Z^2_u = Z^{2,{\frak s}}_u$ be $\cup\{Z^{2,k}_u:k <
\omega\}$, where $Z^{2,k}_u = Z^{2,k,{\frak s}}_u$ is the set of
pairs $(\bar t,\bold q)$ where for some $k < \omega,\bar t \in
{}^k(I^{\frak s}_u)$ and $\bold q \in \bold S^k$.
\nl
4) For $z = (\bar t,\bold q) \in Z^2_u$ let $\partial_{u,z} =
\partial^{\frak s}_{u,z}$ be the following permutation of $D_u$: if $v
\le_{J[{\frak t}]} u$ and $(v,g) \in \{v\} \times K_v$ then
$\partial^{\frak s}_{u,z}((v,g)) = (v,g' g)$ where $g' = g^{v,{\frak
s}}_{\pi_{v,u}(\bar t),\bold q(p)}$ where $p = \text{
tp}_{\text{qf}}(\pi_{v,u}(\bar t),\emptyset,I^{\frak s}_v)$, and, of
course, $\pi_{v,u}(\langle t_\ell:\ell < k\rangle) = \langle
\pi_{v,u}(t_\ell):\ell < k\rangle$.
\nl
5) For $(\bar t,\bold q) \in Z^2_u$ let $g_{\bar t,\bold q} =
g^u_{\bar t,\bold q} = g^{u,{\frak s}}_{\bar t,\bold q} = 
g_{\bar t,\bold q(p)}$ when $p = \text{ tp}_{\text{qf}}(\bar t,\emptyset,I_u)$.
Let $g^{v,{\frak s}}_{\bar t,\bold q} = g^{v,{\frak s}}_{\bar t,\bold q} 
= g^v_{\pi_{v,u}(\bar t),\bold q}$ when $v \le_{J[{\frak s}]} u$.
\enddefinition
\bigskip

\remark{\stag{nx.0.9} Remark}  We can add $\{\partial^{\frak s}_{u,z}:
z \in Z^{2,{\frak s}}_u\}$ to the 
generators of $F^{\frak s}_u$ defined in \scite{nx.4} below.
\endremark
\bigskip

\demo{\stag{nx.0.11} Observation}  In Definition \scite{nx.0.7}(2),
$\partial^{\frak s}_{u,z}$ is a well defined permutation of
$D^{\frak s}_u$.
\enddemo
\bigskip

\demo{Proof}  Easy.
\enddemo
\bigskip

\definition{\stag{nx.4} Definition}  Let ${\frak s}$ be a p.o.w.i.s.
\nl
1) Let $F_u = F^{\frak s}_u$ be
the subgroup of the group of permutations of $D^{\frak s}_u$ generated
by $\{\partial_{u,z}:z \in Z^{\frak s}_u\}$.
\nl
2) For a p.o.w.i.s. ${\frak s}$ let $M_{\frak s}$ be the following model:
\sn
\ub{set of elements}:  $\{(u,g):u \in J^{\frak s}$ and $g \in K^{\frak s}_u\}
\cup \{(1,u,f):u \in J^{\frak s}$ and $f \in F^{\frak s}_u\}$.
\sn
\ub{relations}:  $P^{M_{\frak s}}_{1,u}$, a unary relation, 
is $\{(u,g):g \in K_u\}$ for $u \in J^{\frak s}$,

$P^{M_{\frak s}}_{2,u}$, a unary relation is $\{(1,u,f):f \in F_u\}$ for $u \in
J^{\frak s}$

$R^{M_{\frak s}}_{u,v,h}$, a binary relation, is
$\{((v,g),(1,u,f)):f \in F_u,g \in K_v$ and $f((v,h)) = (v,g)\}$ 
for $u \in J^{\frak s}$ and $v \le_{J[{\frak s}]} u$ and $h \in K_v$.
\enddefinition
\bigskip

\demo{\stag{nx.4.1} Observation}  If ${\frak s}$ is a
$\kappa$-p.o.w.i.s. and $v \le_{J[{\frak s}]} u$ and $f \in F_u$
\ub{then} $f$ maps $\{\nu\} \times K_v = P^{M_{\frak s}}_{1,v}$ onto itself.
\enddemo
\bigskip

\remark{Remark}  If $\pi \in F^{\frak s}_u$ and $v \le_{I^{\frak
s}_u[{\frak s}]} u$ then $\pi \restriction (\{v\} \times K_v)$ comes
directly from $K^{\frak s}_v$, but the relation between the $\langle
\pi \restriction (\{v\} \times K_v):v \le_{I_u[{\frak s}]} u\rangle$ are
less clear.
\endremark
\bigskip

\proclaim{\stag{nx.5} Claim}  Let ${\frak s}$ is a p.o.w.i.s.
\nl
1) $\varkappa$ is an automorphism of $M_{\frak s}$
\ub{iff}:
\mr
\item "{$\circledast$}"  $(a) \quad \varkappa$ is a function with domain
$M_{\frak s}$
\sn
\item "{${{}}$}"  $(b) \quad$ for every $u \in J^{\frak s}$ we have:
{\roster
\itemitem{ ${{}}$ }  $(\alpha) \quad \varkappa 
\restriction D_u \in F^{\frak s}_u$ for every $u \in J^{\frak s}$
\sn
\itemitem{ ${{}}$ }  $(\beta) \quad$ letting 
$f_u = \varkappa \restriction D_u$ we have $(1,u,f) \in P^{M_{\frak s}}_{2,u} 
\Rightarrow \varkappa((1,u,f))$
\nl

\hskip30pt $= (1,u,f_u f)$ where $f_u f$ is the 
product in $F_u$.
\endroster}
\ermn
2) If $f_u \in F_u$ for $u \in J^{\frak s}$ and 
$f_u \subseteq f_v$ for $u \le_{J[{\frak s}]}
v$ \ub{then} there is one and only one automorphism $\varkappa$ of
$M_{\frak s}$ such that $u \in J^{\frak s} \Rightarrow f_u \subseteq
\varkappa$. 
\endproclaim
\bigskip

\demo{Proof}  First assume that $\bar f = \langle f_u:u \in J^{\frak
s}\rangle$ is as in part (2).  We define $\varkappa_{\bar f}$, a
function with domain $M_{\frak s}$ by:
\mr
\item "{$\circledast_1$}"  $(a) \quad$ if $a = (u,g) \in 
P^{M_{\frak s}}_{1,u}$ and $u \in
J^{\frak s}$ then $\varkappa_{\bar f}(a) = f_u(a)$
\sn
\item "{${{}}$}"  $(b) \quad$ if $a = (1,u,f) \in P^{M_{\frak s}}_{2,u}$ then
$\varkappa_{\bar f}(a) = (1,u,f_u f)$.
\ermn
So
\mr
\item "{$\circledast_2$}"  $(a) \quad \varkappa_{\bar f}$ is a well
defined function
\sn
\item "{${{}}$}"  $(b) \quad \varkappa_{\bar f}$ is one to one
\sn
\item "{${{}}$}"  $(c) \quad \varkappa_{\bar f}$ is onto $M_{\frak s}$
\sn
\item "{${{}}$}"  $(d) \quad \varkappa_{\bar f}$ maps 
$P^{M_{\frak s}}_{1,u}$ onto $P^{M_{\frak s}}_{1,u}$ and 
$P^{M_{\frak s}}_{2,u}$ onto $P^{M_{\frak s}}_{2,u}$ for $u \in J^{\frak s}$
\sn
\item "{${{}}$}"  $(e) \quad$ also $\bar f' = \langle f^{-1}_u:u \in
J^{\frak s}\rangle$ satisfies the condition of part (2) and
\nl

\hskip25pt $\varkappa_{{\bar f}'}$ is the inverse of $\varkappa_{\bar f}$
\sn
\item "{${{}}$}"  $(f) \quad \varkappa_{\bar f}$ maps 
$R^{M_{\frak s}}_{u,v,h}$ onto itself.
\ermn
[Why?  The only non-trivial one is clause (f) and in it by clause (e)
it is enough to prove that $\varkappa_{\bar f}$ maps 
$R^{M_{\frak s}}_{u,v,h}$ into $R^{M_{\frak s}}_{u,v,h}$.  So 
assume $v \le_{J[{\frak s}]} u,h \in K_v$ and 
$((v,g),(1,u,f)) \in R^{M_{\frak s}}_{u,v,h}$ hence $f \in
F_u,g \in K_v$ and $f((v,h)) = (v,g)$.  So $\varkappa_{\bar
f}((v,g)) = f_v((v,g))$ and $\varkappa_{\bar f}(1,u,f) = (1,u,f_u f)$ and we
would like to show that $(f_v((v,g)),(1,u,f_u f)) 
\in R^{M_{\frak s}}_{u,v,h}$.
\nl
This means that $(f_u f)((v,h)) = f_v((v,g))$.  We know that
$f((v,h)) = (v,g)$ hence $(f_u f)((v,h)) = f_u(f((v,h))) = 
f_u((v,g))$ so we have to show that $f_u((v,g)) =
f_v((v,g))$.  But $v \le_{J[{\frak s}]} u$ hence (by the assumption
on $\bar f$) we have $f_u \subseteq f_v$ hence $f_u((v,g)) = f_v((v,g))$
so we are done.]

So we have shown that
\mr
\item "{$\circledast_3$}"   if $\bar f = \langle f_u:u \in J^{\frak
s}\rangle$ is as in part (2) then $\varkappa_{\bar f}$ is an
automorphism of $M_{\frak s}$.
\ermn
Next
\mr
\item "{$\circledast_4$}"  if $\varkappa \in \text{ Aut}(M_{\frak s})$
and $\varkappa \restriction P^{M_{\frak s}}_{1,u}$ is the 
identity for each $u \in
J^{\frak s}$ then $\varkappa = \text{ id}_{M_{\frak s}}$.
\ermn
[Why?  By the $R^{M_{\frak s}}_{u,v,h}$'s and $F^{\frak s}_u$ being a group
of permutations of $D_u$.]
\mr
\item "{$\circledast_5$}"  the mapping $\varkappa \mapsto \langle
\varkappa \restriction P^{M_{\frak s}}_{1,u}:u \in J^{\frak s}\rangle$ is a
homomorphism from Aut$(M_{\frak s})$ into $\{\varkappa_{\bar f}:\bar f$ as
above$\}$ with coordinatewise product, with kernel
$\{\varkappa \in \text{ Aut}(M_{\frak s}):\varkappa \restriction
P^{M_{\frak s}}_{1,u} = \text{ id}_{P^{M_{\frak s}}_{1,u}}$ for every 
$u \in J^{\frak s}\}$.
\ermn
[Why?  Easy.]
\mr
\item "{$\circledast_6$}"  the mapping above is onto.
\ermn
[Why?  Given $\varkappa \in \text{ Aut}(M_{\frak s})$, let $f_u =
\varkappa \restriction P^{M_{\frak s}}_{1,u}$.  Clearly $f_u \in F_u$ and $u
\le_{J[{\frak s}]} v \Rightarrow f_u \subseteq f_v$ so $\bar f =
\langle f_u:u \in J^{\frak s}\rangle$ is as above so by
$\circledast_3$ we know $\varkappa_{\bar f}$ is an 
automorphism of $M_{\frak s}$ and $\varkappa_{\bar f} 
\varkappa^{-1}$ is an automorphism of $M_{\frak s}$ which is the
identity on each $P^{M_{\frak s}}_{1,u}$ hence by 
$\circledast_4$ is id$_{M_{\frak s}}$.  
So $\varkappa = \varkappa_{\bar f}$, is as required.]
\mr
\item "{$\circledast_7$}"  the mapping above is one to one.
\ermn
[Why?  Easy by $\circledast_4$.]

Together both parts should be clear.  \hfill$\square_{\scite{nx.5}}$
\enddemo
\bigskip

\definition{\stag{nx.15} Definition}  1) We say that $\bold q_1,\bold
q_2 \in \bold S^k$ are ${\Cal S}$-equivalence where ${\Cal S}
\subseteq {\Cal S}^k$ when $p \in {\Cal S} \Rightarrow \bold q_1(p)
{\Cal E}_{2,p} \bold q_2(p)$.
\nl
2) Omitting ${\Cal S}$ means ${\Cal S} = {\Cal S}^k$.
\enddefinition
\bigskip

\proclaim{\stag{nx.13} Claim}  1) If $u \in J^{\frak s}$ and $f \in
F^{\frak s}_u$ \ub{then} for some $k$ and $\bar t = \langle \bar
t_\ell:\ell < k\rangle \in {}^k(I^{\frak s}_u)$ and $\bold q \in \bold S^k$ we
have:
\mr
\item "{$(*)$}"   $f = \partial_{u,(\bar t,\bold q)}$ (so 
if $v \le_{J[{\frak s}]} u$ then
$f \restriction (\{v\} \times K^{\frak s}_v$) is moving by multiplication by
$g_{(\pi_{v,u}(\bar{\bold t}),\bold q}$, e.g. $g \in K_v \Rightarrow
f((v,g)) = (v,g_{\pi_{v,u}(\bar{\bold t}),\bold q})$.
\ermn
2) $\{\partial_{u,(\bar t,\bold q)}:\bar t \in {}^k(I^{\frak s}_u)$
and $\bold q \in \bold S^k$ for some $k\}$ is a group of
permutations of $D^{\frak s}_u$ which include $F^{\frak s}_u$. 
\nl
3) For every $\bold q \in \bold S^k$ there is a reduced $\bold q' \in
\bold S^k$ which is $\bold S^k$-equivalent to it (see Definition
\scite{nx.0.7}(2). 
\endproclaim
\bigskip

\remark{\stag{nx.14} Remark}  1) We can be somewhat more restrictive.
\endremark
\bigskip

\demo{Proof}  We use freely Definition \scite{nx.0.7}.
Recall that $F^{\frak s}_u$ is the group of permutations
of $D^{\frak s}_u$ generated by $\{\partial_{u,z}:
z \in Z^{\frak s}_u\}$.   Hence it is enough to prove that 
$f \in F^{\frak s}_u$ satisfies the conclusion
of the claim in the following cases.
\bn
\ub{Case 0}:  $f$ is the identity.

It is enough to let $k=0$ so ${\Cal S}^k$ is a singleton $\{p\}$ and
$\bold q(p)$ is the sequence $<<>>$, i.e. we use in Definition
\scite{n.1.17}(1) the case $j(*)=0$, i.e. \scite{n.1}(3) for $k=0$.
\bn
\ub{Case 1}:  $f = \partial_{u,z}$ where $z \in Z^0_u$.

So $z = X_{I_u}$ let $k = n(z)+1,\bar t = \bar t^z$.  
We define $\bold q$ as follows:
\mr
\item "{$(a)$}"  if $q \in {\Cal S}^k$ ``says" that 
$\bar t = \langle t_\ell:\ell \le n(z)\rangle$ is decreasing then 
$g_{\bar t,\bold q}$ is $g_z$
\sn
\item "{$(b)$}"  if not then $g_{\bar t^z,\bold q} = e_{K_u}$.
\endroster
\bn
\ub{Case 2}:  $f = \partial_{u,z}$ where $z \in Z^1_u$.

Also clear.
\bn
\ub{Case 3}:  $f = f_1 f_2$ (product in $F^{\frak s}_u$) where
$f_1,f_2 \in F^{\frak s}_u$ satisfies the conclusion of the claim.

Just combine the definitions.
\bn
\ub{Case 4}:  $f = f^{-1}$ where $f \in F^{\frak s}_u$ satisfies the
conclusion of the claim.

Easy, too.  \hfill$\square_{\scite{nx.13}}$
\enddemo
\bigskip

\remark{\stag{nx.16} Remark}  If $q \in {\Cal S}^k$ 
and $\bold q_1,\bold q_2 \in {\bold S}^k$
and $v \le_{J[{\frak s}]} u,\bar{\bold t} = {}^k(I_u)$ 
and $q = \text{ tp}_{\text{qf}}
(\pi^{\frak s}_{v,u}(\bar t),\emptyset,I_v)$ and
$\bold q_1(q),\bold q_2(q)$ are not ${\Cal E}_{2,q}$-equivalent,
\ub{then} $g_{\bar{\bold t},\bold q_1} \ne g_{\bar{\bold t},\bold q_2}$.
\endremark
\bigskip

\demo{Proof}  This is by \scite{n.1.17}(4).
\enddemo
\newpage

\head {\S3 The main result} \endhead  \resetall \sectno=3
 \spuriousreset
\bn
We can prove that every $\kappa$-parameter has a limit, but for our
application it is more transparent to consider $\kappa$-parameber
${\frak s}$ which is the $\kappa$-parameter ${\frak t}$ + its limit.
\bigskip

\definition{\stag{p.1} Definition}  We say that ${\frak s}$ is the
limit of ${\frak t}$ as witnessed by $v_*$ when (both are
p.o.w.i.s. and)
\mr
\item "{$(a)$}"  $J^{\frak t} \subseteq J^{\frak s}$ and
$J^{\frak s} = J^{\frak t} \cup \{v_*\},v_* \notin
J^{\frak t}$ and $u \in J^{\frak s} \Rightarrow u \le_{J[{\frak s}]} v_*$
\sn
\item "{$(b)$}"  $I^{\frak s}_u = I^{\frak t}_u$ and
$\pi^{\frak s}_{u,v} =
\pi^{\frak t}_{u,v}$ when $u \le_{J[{\frak s}]} v <_{J[{\frak s}]}
v_*$
\sn
\item "{$(c)$}"   if $t \in I^{\frak s}_{v_*}$ then for some $u = u_t
\in J^{\frak s}$ we have $t \in \text{ Dom}(\pi^{\frak s}_{u_t,v_*})$,
moreover (if ${\frak s}$ is nice this follows)
$J^{\frak s} \models ``u_t \le v < v_*" 
\Rightarrow t \in \text{ Dom}(\pi^{\frak s}_{v,v_*})$
\sn
\item "{$(d)$}"  if $s,t \in I^{\frak s}_{v^*}$
\ub{then} for some $u = u_{s,t} \in J^{\frak t}$ for every $v$ satisfying
$u \le_{J[{\frak s}]} v <_{J[{\frak s}]} v_*$ we have 
$I^{\frak s}_{v_*} \models ``s  < t" \Leftrightarrow
\pi^{\frak s}_{v,v_*}(s) <_{I^{\frak s}_v} 
\pi^{\frak s}_{v,v_*}(t)$
\sn
\item "{$(e)$}"  if $\langle t_u:u \in J^{\frak t}_{\ge w}\rangle$ is
a sequence satisfying $w \in J,J_{\ge w} = \{u:w \le u \in J\};t_u \in 
I^{\frak s}_u$ and $w \le u_1 \le u_2 \in J^{\frak t} \Rightarrow 
\pi_{u_1,u_2}(t_{u_2}) =
 t_{u_1}$, \ub{then} there is a unique $t \in I^{\frak s}_{v_*}$
 such that $u \in J^{\frak t}_{\ge w} \Rightarrow \pi_{u,v_*}(t) = t_u$.
\endroster
\enddefinition
\bigskip

\definition{\stag{p.1.7} Definition}  We say that ${\frak s}$ is an
 existential limit of ${\frak t}$ \ub{when}: clauses (a)-(e) of Definition
 \scite{p.1} holds and
\mr
\item "{$(f)$}"  assume that
{\roster
\itemitem{ $(\alpha)$ }  $u_* \in J^{\frak t}$
\sn
\itemitem{ $(\beta)$ }   $k_1,k_2 < \omega$ and $k = k_1 + k_2$
\sn
\itemitem{ $(\gamma)$ }   $E$ is an equivalence relation on ${\Cal
S}^k$
\sn
\itemitem{ $(\delta)$ }  $\bar e = \langle e_u:u \in J^{\frak t}_{\ge
u_*}\rangle$, where $e_u$ is an $E$-equivalence class
\sn
\itemitem{ $(\varepsilon)$ }   $\bar t \in {}^{k_1}(I^{\frak s}_{v_*})$
\sn
\itemitem{ $(\zeta)$ }   for every $v \in J^{\frak t}_{\ge u_*}$
there is $\bar s_v \in {}^{k_2}(I^{\frak t}_{w(v)})$ such that:
\nl

$\quad$ if $u_* \le_{J[{\frak t}]} u \le_{J[{\frak t}]} v$ then $e_u$ is the
$E$-equivalence class of 
\nl

$\quad$ tp$_{\text{qf}}(\bar t^u \char 94 \bar s^{u,v},
\emptyset,I^{\frak t}_u)$ where $\bar t^u = \pi^{\frak s}_{u,v_*}(\bar t)$ and
$\bar s^{u,v} = \pi^{\frak t}_{u,v}(\bar s_v)$.
\endroster}
\ermn
\ub{Then} there is $\bar s \in {}^{k^*}(I^{\frak s}_{v_*})$ such that
for every $u \in J^{\frak t}$ large enough tp$(\pi^{\frak s}_{u,v_*}
(\bar t \char 94 \bar s),\emptyset,I^{\frak t}_u)$ belongs
to $e_u$ (and is constantly $p_*$ for some $p_* \in {\Cal S}^k$).
\enddefinition
\bigskip

\remark{\stag{p.1A} Remark}  We may say ``${\frak s}$ is semi-limit of
${\frak t}$" when in clause (d) we replace $\Leftrightarrow$ by
$\Rightarrow$.  We may consider using this weaker version and/or omit
linearity in our main theorem, but the present version suffices.
\endremark
\bigskip

\proclaim{\stag{p.2} Main Claim}  $K^{\frak s}_{v_*}$ is an almost
 $\kappa$-automorphism group (see below) \ub{when}:
\mr
\item "{$\boxtimes$}"  $(a) \quad {\frak s},{\frak t}$ are both p.o.w.i.s
\sn
\item "{${{}}$}"  $(b) \quad {\frak s}$ is an existential 
limit of ${\frak t}$ as  witnessed by $v_*$ 
\sn
\item "{${{}}$}"  $(c) \quad J^{\frak t}$ is $\aleph_1$-directed and
is linear (i.e., for every $u,v \in J^{\frak t}$ we have 
\nl

\hskip25pt $u \le_{J[{\frak t}]} v$ or $v \le_{J[{\frak t}]} v)$
\sn
\item "{${{}}$}"   $(d) \quad {\frak t}$ is a $\kappa$-p.o.w.i.s (so
$\kappa \ge |J^{\frak t}|$ and $\kappa \ge
|I^{\frak t}_u|$ for $u \in J^{\frak t}$)
\sn
\item "{${{}}$}"   $(e) \quad {\frak t}$ is non-trivial (see
Definition \scite{m.1}(6).
\endroster
\endproclaim
\bigskip

\remark{Remark}  Not much harm in adding ${\frak t}$ is 
nice (see Definition \scite{m.2}(7)) so for 
$u \le_{J[{\frak t}]} v$ the functions 
$\pi^{\frak t}_{u,v},\hat \pi^{\frak t}_{u,v}$ has full
domain, see Definition \scite{nx.0}(1),(2),(3) and Claim \scite{nx.2}(3)).
\endremark
\bigskip

\definition{\stag{p.1.1} Definition}  $G$ is an almost
$\kappa$-automorphism group when: there is a $\kappa$-automorphism
group $G^+$ and a normal subgroup $G^-$ of $G^+$ of cardinality $\le
\kappa$ such that $G$ is isomorphic to $G^+/G^-$, i.e., there is a
homomorphism from $G^+$ onto $G$ with kernel $G^-$.
\enddefinition
\bn
Before proving \scite{p.2} we explain: why
being almost $\kappa$-automorphism group help us in proving our
intended result?
\nl
Recalling \scite{0.6}:

\proclaim{\stag{p.1.2} Claim}  For any ordinal $\alpha$, if there is an almost
$\kappa$-automorphism group $G$ with a subgroup $H$ of cardinality $\le
\kappa$ such that $\tau'_{G,H} = \alpha$ [such that
{\rm nor}$^\alpha_G(H)=G \wedge (\forall \beta < \alpha)(\text{\rm
nor}^\beta_G(H) \ne G)$] \ub{then} there is a $\kappa$-automorphism group
$G'$ with a subgroup $H'$ of cardinality $\le \kappa$ such that
$\tau'_{G',H'} = \alpha$ [such that {\rm nor}$^\alpha_G(H') = G \wedge
(\forall \beta < \alpha)(\text{\rm nor}^\beta_G(H) \ne G)$]. 
\endproclaim
\bigskip

\demo{Proof}  Easy.  

Let $G^+,G^-$ be as in Definition \scite{p.1.1} and $h$ be a
homomorphism from $G^+$ onto $G$ with kernel $G^-$ and let $H^+ = \{x
\in G^+:h(x) \in H\}$.

So it is easy to check each of the following statements (similar to
\scite{nx.2}(5)):
\mr
\item "{$\circledast$}"  $(a) \quad H^+$ is a subgroup of $G^+$
\sn
\item "{${{}}$}"  $(b) \quad |H^+| \le |H| \times |G^-| \le \kappa 
\kappa = \kappa$
\sn
\item "{${{}}$}"  $(c) \quad G^+$ is a $\kappa$-automorphism group
\sn
\item "{${{}}$}"  $(d) \quad$ nor$^\beta_{G^+}(H^+) = \{x \in
G^+:h(x) \in \text{\rm nor}^\beta_G(H)\}$ for every $\beta \le \infty$
\sn
\item "{${{}}$}"  $(e) \quad \tau_{G,H} = \tau_{G^+,H^+}$
\sn
\item "{${{}}$}"  $(f) \quad$ nor$^\beta_G(H)=G$ then
nor$^\beta_{G^+}(H^+) = G^+$ for every $\beta \le \infty$.
\ermn
Together $(G^+,H^+)$ exemplifies the desired conclusion.
\hfill$\square_{\scite{p.1.2}}$ 
\enddemo
\bigskip

\demo{Proof of \scite{p.2}}  Let $G^+$ be the 
automorphism group of $M_{\frak t}$ and let $G^-$ be
the following subgroup of $G^+$

$$
\align
\{\varkappa \in G^+:&\text{ for some } u \in J^{\frak t} \text{ we have} \\
  &u \le_J v \wedge g \in K_v \Rightarrow \varkappa((v,g)) = (v,g)\}.
\endalign
$$
\mn
Easily
\mr
\item "{$\circledast_1$}"  $G^-$ is a subgroup of $G^+$
\nl
[Why?  As $J^{\frak t}$ is directed]
\sn
\item "{$\circledast_2$}"  for every $\varkappa \in G^+$ we can find $\bar
f^\varkappa = \langle f^\varkappa_u:u \in J^{\frak t}\rangle$ such that
{\roster
\itemitem{ $(a)$ }  $f^\varkappa_u \in F^{\frak t}_u$
\sn
\itemitem{ $(b)$ }  $\varkappa \restriction D^{\frak t}_u = f_u$
\sn
\itemitem{ $(c)$ }  $\varkappa \restriction P^{M_{\frak t}}_{2,u}$ is
$(1,u,f) \mapsto (1,u,f_u f)$.
\nl
[Why?  By Claim \scite{nx.5}.]
\endroster}
\item "{$\circledast_3$}"  $G^-$ has cardinality $\le \kappa$.
\nl
[Why?  As $|J^{\frak t}| \le \kappa$, it 
suffices to prove that for each $u \in J^{\frak t}$, the
subgroup $G^-_u := \{\varkappa \in G^+:\varkappa \restriction 
P^{M_{\frak t}}_{1,v}$ is the identity when $u \le_{J[{\frak s}]} v\}$
has cardinality $\le
\kappa$, but this has the same number of elements as $F^{\frak s}_u$
because $\varkappa \mapsto \varkappa \restriction D_u$ is a 
one-to-one function from
$G^-_u$ onto $F^{\frak s}_u$ and ${\frak t}$ is linear.  
As $|F^{\frak s}_u| \le \aleph_0 +
|Z_u| = \aleph_0 + |I_u| \le \kappa$ we are done.]
\sn
\item "{$\circledast_4$}"  $G^-$ is a normal subgroup of $G^+$.
\nl
[Why?  By its definition,  more elaborately
{\roster
\itemitem{ $(a)$ }   each $G^-_u$ is a normal subgroup of $G^+$.
\nl
[Why?  As all members of Aut$(M_{\frak s})$ maps each $\{v\} \times
K_v$ onto itself so $G^-_u$ is even an definable subgroup]
\sn
\itemitem{ $(b)$ }  $u \le_{J[{\frak t}]} v 
\Rightarrow G^-_u \subseteq G^-_v$.
\nl
[Why?  Check the definitions.]
\sn
\itemitem{ $(c)$ }   $G^- = \cup\{G^-_u:u \in J\}$.
\nl
[Why?  Trivially.]
\endroster}
\ermn
Together we are done proving $\circledast_4$.]
\mr
\item "{$\circledast_5$}"  For 
$x \in Z^{\frak s}_{v^*}$ let $\varkappa_x$ be the
following automorphism of $M_{\frak t}$, it is defined
as in $\circledast_2$ by $\langle
f^x_u:u \in J^{\frak t}\rangle$ where $f_u = 
\partial^{\frak t}_{u,\hat\pi_{u,v_*}(x)}$ is from Definition
\scite{nx.0}(4)
\sn
\item "{$\circledast_6$}"  for every $x \in Z^{\frak s}_{v^*},
\varkappa_x$ is a well defined automorphism of $M_{\frak t}$.
\nl
[Why?  Look at the definitions and \scite{nx.5}.]
\ermn
The main point is
\mr
\item "{$\circledast_7$}"  $G^+$ is generated by $\{\varkappa_x:x \in
Z^{\frak s}_{v^*}\} \cup G^-$.
\ermn
Why?  Clearly the set is a set of elements of $G^+$.  So assume $\varkappa \in
G^+$ and let $\bar f^\varkappa = \langle f^\varkappa_u:u \in J^{\frak t}
\rangle$ be as in $\circledast_2$, they are fixed for awhile.

By \scite{nx.13} for each $u \in J^{\frak t}$ there are $k=k^u$ and
$\bar t = \bar t^u \in {}^{k^u}(I^{\frak s}_u)$ and $\bold q = \bold
q^u \in \bold S^{k^u}$ such that (the ``disjoint" as we can replace
$\bar t$ by $\bar t \char 94 \bar t$ or even $\bar t \char 94 \bar t
\char 94 \ldots \char 94 \bar t$ with $|S^{k^u}|$ copies note that we
can demand that $\bold q$ is reduced by \scite{nx.13}(3)):
\mr
\item "{$\boxdot_1$}"   $f^{\varkappa}_u = 
\partial_{u,(\bar t^u,\bold q)}$, i.e., if
$v \le_{J[{\frak t}]} u$ then $f \restriction (\{v\} \cap K^{\frak
t}_v)$ is a multiplication from the left (of the $K^{\frak t}_v$-coordinate)
 by $g_{\pi^{\frak
t}_{v,u}(\bar t),\bold q^u}$ and $\bold q^u$ is reduced and
disjoint, see Definition \scite{nx.0.7}(2),(5).
\ermn
The choices are not necessarily unique, in particular
\mr
\item "{$\boxdot_2$}"  if $u^1 \le_{J[{\frak t}]} u^2$ then
$(k^{u^2},\pi_{u^1,u^2}(\bar t^{u^2}),\bold q^{u^2})$ can serve as
$(k^{u^1},\bar t^{u^1},\bold q^{u^1})$.
\ermn
Also
\mr
\item "{$\boxdot_3$}"  the set of possible $(k^u,\bold q^u)$ is
countable.  
\ermn
As $J^{\frak t}$ is $\aleph_1$-directed
\mr
\item "{$\boxdot_4$}"  for some pair $(k^*,\bold q^*)$ the set $\{u
\in J^{\frak t}:k^u = k^*$ and $\bold q^u = \bold q^*\}$ is cofinal in
$J^{\frak t}$.
\ermn
Together, \wilog \, for some $k^*,\bold q$
\mr
\item "{$\boxdot_5$}"  $k^u = k^*$ and $\bold q^u = \bold q$ for
every $u \in J^{\frak t}$.
\ermn
Let $E$ be an ultrafilter on $J^{\frak t}$ such that $u \in J^{\frak
t} \Rightarrow \{v:u \le_{J[{\frak t}]} u\}\in E$, exists as $J^{\frak
t}$ is directed.  For each $u \in J^{\frak t}$ there are
$A_u,p_u,w(u)$ such that
\mr
\item "{$\boxdot_6$}"  $(a) \quad A_u \in E$ and
\sn
\item "{${{}}$}"  $(b) \quad p_u \in {\Cal S}^{k^*}$
\sn
\item "{${{}}$}"  $(c) \quad$ if $v \in A_u$ then $u \le_{J[{\frak
t}]} v$ and $p_u = \text{ tp}(\pi_{u,v}(\bar t^v),\emptyset,I_u)$
\sn
\item "{${{}}$}"  $(d) \quad w(u) \in A_u$.
\ermn
For $p \in {\Cal S}^{k^*}$ let
\mr
\item "{$\boxdot_7$}"  $(a) \quad Y_p = \{u \in J^{\frak t}:p_u =
p\}$
\sn
\item "{${{}}$}"  $(b) \quad \bar s^{u,v} = \pi^{\frak t}_{u,v}(\bar
t^v) \restriction \text{ supp}(\bold q(p_u))$ for $u \in J^{\frak t},v
\in A_u$
\sn
\item "{${{}}$}"  $(c) \quad \bar s^u = \bar s^{u,w(u)}$.
\ermn
So
\mr
\item "{$\boxdot_8$}"  $\langle Y_p:p \in S^{k^*}\rangle$ is a
partition of $J^{\frak t}$.
\ermn
Fix $p \in {\Cal S}^k$ for awhile so for each $u \in Y_p$ and $v \in
A_u$ by $\boxdot_1$, $\varkappa \restriction (\{u\} \times K_u)$ is
multiplication from the left by $g^{u,{\frak s}}_{\pi^{\frak t}_{u,v}(\bar
t^v),\bold q}$ (it was $\bold q^v$ but we have already agreed that
$\bold q^v = \bold q$).  But $p = \text{ tp}_{\text{qf}}
(\pi^{\frak t}_{u,v}(\bar
t^v),\emptyset,J_u)$ as $u \in Y_p,v \in A_u$ and so by Definition
\scite{nx.0.7}(4) we know that
$g^{u,{\frak s}}_{\pi^{\frak t}_{u,v}(\bar t^v),\bold q}$ 
is $g^{u,{\frak s}}_{\pi^{\frak t}_{u,v}(\bar t^v),\bold q(p)}$.

Now $\bold q(p) \in \Lambda^2_{k^*}$ so $\bold q(p) = \langle
\rho^p_0,\rho^p_1,\dotsc,\rho^p_{i(p)-1}\rangle$ and recall

$$
g_{\pi^{\frak t}_{u,v}(\bar t^v),\bold q(p)} \text{ is } 
g_{\bar t,\rho^p_0} h_{g_{\bar t,\rho^p_1} G^{<0}_u} \ldots;
$$
\mn
so it depends only on $\bar t \restriction \text{ supp}(\bold q(p))$ only.

Now consider any two members $v_1,v_2$ of $A_u$ (so they are above $u$)
comparing the two expressions for $\varkappa \restriction (\{u\}
\times K_u)$ one coming from $v^1$ the second from $v^2$ we conclude
that $g_{\pi^{\frak t}_{u,v^s_2}(\bar t^{v^*_1}),\bold q(p)} =
g_{\pi^{\frak t}_{u,v_2}(\bar t^{v_2}),\bold q(p)}$.  As $\bold q$ is
reduced also $\bold q(p)$ is $p$-reduced hence by \scite{n.1.23}(3) 
we conclude that
\mr
\item "{$\boxdot_9$}"   if ($p \in {\Cal S}^{k^*},u \in Y_p \subseteq
J^{\frak t}$ and) $v_1,v_2 \in A_u$ then $\pi^{\frak t}_{u,v_2}(\bar
t^{v_1}) \restriction$ supp$(\bold q(p))$ is a permutation of
$\pi^{\frak t}_{u,v_2}(\bar t^{v_2}) \restriction \text{ supp}(\bold
q(p))$
\nl
this means
\sn
\item "{$\boxdot_{10}$}"  if $p \in {\Cal S}^{k^*},u \in J^{\frak t}$
and $v_1,v_2 \in A_u$ then $\bar s^{u,v_1}$ is a permutation of $\bar
s^{u,v_2}$.
\ermn
Hence for each $u \in J^{\frak t}$
\mr 
\item "{$\boxdot_{11}$}"  if $v \in A_u$ then $\bar s^{u,v}$ is a
permutation of $\bar s^u = \bar s^{u,w(u)}$.
\ermn
As there are only finitely many permuations of $\bar s^{u,v_u}$, there
are $\bar s^u,A'_u$ such that
\mr
\item "{$\boxdot_{12}$}"  for $u \in J^{\frak t}$:
{\roster
\itemitem{ $(a)$ }  $A'_u \in E$
\sn
\itemitem{ $(b)$ }  $A'_u \subseteq A_u$
\sn
\itemitem{ $(c)$ }  $\bar s^u = \bar s^{u,v}$ for every $v \in A'_u$.
\endroster}
\ermn
Now
\mr  
\item "{$\boxdot_{13}$}"  if $p \in {\Cal S}^k$ and $u_1 \le_{J[{\frak
t}]} u_2$ are from $Y_p$ then $\pi^{\frak t}_{u_1,u_2}(\bar s^{u_2}) =
\bar s^{u_2}$.
\ermn
[Why?  As $E$ is an ultrafilter on $J^{\frak t}$ and
$A'_{u_1},A'_{u_2} \in E$ we can find $v \in A'_{u_1} \cap A'_{u_2}$.
So for $\ell=1,2$ we have $\bar s^{u_\ell} = \pi^{\frak
t}_{u_\ell,v}({\frak t}^v) \restriction \text{ supp}(\bold q(p)) =
\pi^{\frak t}_{u_\ell,v}(\bar t^v \restriction \text{ supp}(\bold
q(p))$.

As $\pi^{\frak t}_{u_1,v} = \pi^{\frak t}_{u_1,u_0} \circ \pi^{\frak
t}_{u_2,v}$ we conclude $\bar s^{u_1} = \pi^{\frak t}_{u_1,u_2}(\bar
s^{u_2})$ is as required.]

Let ${\Cal S}' = \{p \in {\Cal S}^{k^*}:Y_p$ is an unbound subset of
$J^{\frak t}\}$, so for some $u_* \in J^{\frak t}$ we have
\mr
\item "{$\boxdot_{14}$}"   $J^{\frak t}_{\ge u_*} \subseteq
\cup\{Y_p:p \in {\Cal S}'\}$.
\ermn
Also \wilog 
\mr
\item "{$\boxdot_{15}$}"  $k^* = k^*_1 + k^*_2$ and
$\{0,\dotsc,k^*_1-1\} = \cup\{\text{supp}(\bold q(p):p \in {\Cal S}'\}$
\sn
\item "{$\boxdot_{16}$}"  for $p \in {\Cal S}'$ and $\ell \in 
\text{ supp}(\bold q(p))$, so $s^u_\ell$ is well defined for $u \in
Y_p$, there is a unique $t \in J_{\frak s}$ such that:
$$
u \in Y_p \Rightarrow \pi^{\frak s}_{u,v_*}(t) = s^u_\ell.
$$
\ermn
[Why?  By clause (d) of Definition \scite{n.1}.]

Next we can find $\bar t$ such
that
\mr
\item "{$\boxdot_{17}$}"  $(a) \quad \bar t = \langle t_\ell:\ell <
k^*_1\rangle$
\sn
\item "{${{}}$}"  $(b) \quad$ if $p \in {\Cal S}'$ and $\ell \in 
\text{ supp}(\bold q(p))$ then $t_\ell \in I^{\frak s}_{v^*}$ 
is as in $\boxdot_{16}$.
\ermn
[Why?  For $i \in \cup\{\text{supp}(\bold q(p)):p \in {\Cal S}'\}$ use
$\boxdot_{16}$, as $\bold q$ is disjoint (see Definition
\scite{nx.0.7}(2)) there is no case of ``double definition".]

By clause (d) of Definition \scite{p.1}, possibly increasing $u_*$
\mr
\item "{$\boxdot_{18}$}"  $p^* = \text{ tp}(\pi^{\frak
s}_{u,v_*}(\bar t),\emptyset,I^{\frak s}_u)$ for every $u \in
J^{\frak t}_{\ge u_*}$
\sn
\item "{$\boxdot_{19}$}"  let ${\Cal E}$ be the following equivalence
relation on ${\Cal S}^{k^*},p_1 {\Cal E} p_2 \Leftrightarrow \bold
q(p_1) {\Cal E}^1_{k^*_1,p \restriction k^*_1} \bold q(p_2)$; note
they are actually from ${\Cal S}^{k^*_1}$ and so ``${\Cal
E}^1_{k^*_1,p \restriction k^*_1}$-equivalent" is meaningful, see
Definition \scite{n.1.15}(4) 
\sn
\item "{$\boxdot_{20}$}"  let $\bar e = \langle e_u:u \in 
J^{\frak t}_{\ge u_*}\rangle$ be defined by $ e_u = p_u/E$
\sn
\item "{$\boxdot_{21}$}"  $E,\bar t,\bar e,\langle
\pi^{\frak t}_{u,w(u)}(\bar t^{w(u)}):u \in J^{\frak t}_{\ge
u_*}\rangle$ satisfies the demands $(f)(\alpha)-(\zeta)$ from
Definition \scite{p.1.7}.
\ermn
[Why?  Check.]

Recall $p^* = \text{ tp}(\bar t,\emptyset,I^{\frak s}_{v_*})$ here so
let $\bar s \in {}^{(k^*_2)}(I^{\frak s}_{v_*})$ be as guaranteed to
exist by Definition \scite{p.1.7}.  
Let $\bar t^{v^*} := \bar t \char 94 \bar s$.  So
possibly increasing $u_* \in J^{\frak t}$ for some $p^*$ we have
\mr
\item "{$\boxdot_{22}$}"   if $u \in J^{\frak t}_{\ge u_*}$
then $p^* = \text{ tp}(\pi^{\frak s}_{u,v_*}(\bar t \char 94 \bar
s),\emptyset,I^{\frak s}_u) = \text{ tp}(\bar t \char 94 \bar
s,\emptyset,I^{\frak s}_{v_*})$.
\ermn
Let
\mr
\item "{$\boxdot_{23}$}"  $(a) \quad \varrho^* = \bold q(p^*)$ so
$\varrho^* \in \Lambda^2_{k^*_1,p^*}$ and let $\varrho^* = \langle
\rho_\ell:\ell < \ell(*)\rangle$
\sn
\item "{${{}}$}"  $(b) \quad \bar t_u = \pi^{\frak s}_{u,v_*}(\bar t)$
for $u \in J^{\frak t}$
\sn
\item "{${{}}$}"  $(c) \quad$ let $z_u = z^{u,s}_{\bar t_u,\varrho}
\in Z^{1,{\frak s}}_u$ (see Definition \scite{n.1.17}(5A))
\sn
\item "{${{}}$}"  $(d) \quad$ let $f_u = \partial^{\frak s}_{u,z_u}
\in F^{\frak s}_u$; (this is not the same as $f^\varkappa_u$!).
\ermn
Now
\mr
\item "{$\boxdot_{24}$}"  for $u_1 \le_{J[{\frak t}]} u_2$ we have
$f_{u_1} \subseteq f_{u_2}$.
\ermn
[Why?  Check.]
\mr 
\item "{$\boxdot_{25}$}"   $\varkappa_{\bar f}$ is a finite product of
members of $\{\varkappa_x:x \in Z^{\frak s}_{v_*}\}$.
\ermn
[Why?  Recall $\varkappa_x$ for $x \in Z^{\frak s}_{v^*}$ is from
$\circledast_5$.  Now use $\boxdot_{23}$.]
\sn
Lastly
\mr
\item "{$\boxdot_{26}$}"  $(\varkappa^{-1}_{\bar f})\varkappa \in G^+ 
= \text{ Aut}(M_{\frak t})$ is the identity on $P^{M_{\frak t}}_u$ whenever $u
\in J^{\frak t}_{\ge u_*}$.
\ermn
[Why?  By $\boxdot_{24}$ and our choices.]
\mr
\item "{$\boxdot_{25}$}"  $(\varkappa_{\bar f}) \in (G^-_{u_*}
\subseteq) G^-$.
\ermn
[Why?  By $\boxdot_{25}$ and the definition of $(G_{u^*}$ and) $G^-$.]
\mr 
\item "{$\boxdot_{28}$}"  $\varkappa$ is the product (in $G^+$) of 
$\varkappa_{\bar f} \in G^-$ and $(\varkappa^{-1}_f)\varkappa \in
\langle \{\varkappa_x:x \in Z^{\frak s}_{v_*}\}\rangle$.
\ermn
[Why?  $\boxdot_{25} + \boxdot_{27}$ this is clear.]

As $\varkappa$ was any a member of $G^+$ we are done proving $\circledast_7$.
\mr
\item "{$\circledast_8$}"  there is a homomorphism $\bold h$ from
$K^{\frak s}_{v_*}$ onto $G^+/G^-$ which maps $g_x$ to $\varkappa_x
G^-$ for $x \in Z^{\frak s}_{v_*}$.
\ermn
[Why?  By $\circledast_7$ there is at most one such homomorphism and
if it exists it is onto.

So it is enough to show that for any group term, $\sigma$ if 
$K^{\frak s}_{v_*}$
satisfies $K_{v_*} \models ``\sigma(g_{x_1},\dotsc,g_{x_{k-1}})=e"$ then
$\sigma(\varkappa_{x_0},\dotsc,\varkappa_{x_{k-1}}) \in G^-$.  Let
$\langle t_\ell:\ell < \ell^*\rangle$ list
$\cup\{\text{his}(x_\ell):\ell < k\} \subseteq I^{\frak s}_{v_*}$ 
and let $u_* \in J^{\frak t}$ be
such that: if $u_* \le_{J[{\frak t}]} u$ and $\ell(1),\ell(2) < \ell^*$
we have $I^{\frak s}_{v_*} \models t_{\ell(1)} <_I t_{\ell(2)}$ iff
$I^{\frak t}_u \models \pi_{u,v_*}(t_{\ell(1)}) <
\pi_{u,v^*}(t_{\ell(2)})$ and similarly for equality, see clause (d)
of Definition \scite{p.1}.

Let $t_{u,\ell} = \pi_{u,v_*}(t_\ell),x_{u,\ell} = \hat
\pi_{u,v_*}(x_\ell)$.  By the definition of $G^-$ it is enough to show
that: if $u_* \le_{J[{\frak t}]} u$ then $K_u \models
``\sigma(g_{x_{u,0}},\dotsc,g_{x_{u,k_1}}) = e_{K_u}"$.  By the
analysis in \scite{m.7} and \S2 (i.e., twisted product) this should be clear.]
\mr
\item "{$\circledast_9$}"  $\varkappa^*$ is one to one.
\ermn
[Why?  By part of the analysis as for $\circledast_7$.]

By $\circledast_8 + \circledast_9$ we are done.
\hfill$\square_{\scite{p.1.2}}$ 
\enddemo
\bigskip

\proclaim{\stag{p.3} Theorem}   Assume
\mr
\item "{$(a)$}"  $\aleph_0 < { \text{\rm cf\/}}(\theta) \le \theta \le
\kappa$
\sn
\item "{$(b)$}"  ${\Cal F}_\alpha \subseteq {}^\alpha \kappa$ for $\alpha <
\theta$ has cardinality $\le \kappa$ (also ${\Cal F}_\alpha \subseteq
{}^\alpha \beta$ for some $\beta < \kappa^+$ is O.K.)
\sn
\item "{$(c)$}"  ${\Cal F} = \{f \in{}^\theta \kappa:f
\restriction \alpha \in {\Cal F}_\alpha$ for every $\alpha < \theta\}$
\sn
\item "{$(d)$}"  $\gamma = { \text{\rm rk\/}}({\Cal F},
<_{J^{\text{bd}}_\theta})$, necessarily $< \infty$ so $<
(\kappa^\theta)^+$
\sn
\item "{$(e)$}"  for $f_1,f_2 \in {\Cal F}$, \ub{then} $f_1 
<_{J^{\text{bd}}_\theta} f_2$ or $f_2 <_{J^{\text{bd}}_\theta} f_1$
or $f_2 =_{J^{\text{bd}}_\theta} f_1$; follows from (f)
\sn
\item "{$(f)$}"  for stationarily many $\delta < \theta$ we have: if
$f_1,f_2 \in {\Cal F}_\delta$, \ub{then} for some $\alpha < \delta$ we
have $\beta \in (\alpha,\delta) \Rightarrow (f_1(\beta) < f_2(\beta)
\equiv f_1(\alpha) < f_2(\alpha))$.
\ermn
\ub{Then} $\tau^{\text{atw}}_\kappa \ge \tau^{\text{nlg}}_\kappa 
\ge \tau^{\text{nlf}}_\kappa > \gamma$
(on $\tau^{\text{nlf}}_\kappa$ see Definition \scite{0.2}(4)).
\endproclaim
\bigskip

\proclaim{\stag{p.3A} Theorem}  We can in Theorem \scite{p.3} weaken
clause (f) to
\mr
\item "{$(f)'$}"  $(\alpha) \quad S \subseteq \theta$ is a stationary set
consisting of limit ordinals
\sn
\item "{${{}}$}"  $(\beta) \quad D$ is a normal filter on $\theta$
\sn
\item "{${{}}$}"  $(\gamma) \quad S \in D$
\sn
\item "{${{}}$}"  $(\delta) \quad \bar J = \langle J_\delta:\delta \in
S\rangle$
\sn
\item "{${{}}$}"  $(\varepsilon) \quad J_\delta$ is an ideal on
$\delta$ extending $J^{\text{bd}}_\delta$ for $\delta \in S$
\sn
\item "{${{}}$}"  $(\zeta) \quad$ if $S' \subseteq S,S' \in D^+$ and
$w_\delta \in J_\delta$ for $\delta \in S'$ then 
\nl

\hskip25pt $\cup\{\delta \backslash w_\delta:\delta \in S'\}$
contains an end segment of $\theta$
\sn
\item "{${{}}$}"  $(\eta) \quad$ if $\delta \in S$ and $f_1,f_2 \in
{\Cal F}$ \ub{then} $f_1 \restriction \delta <_{J_\delta} 
f_2 \restriction \delta$ or 
\nl

\hskip25pt $f_2 \restriction \delta <_{J_\delta} 
f_1 \restriction \delta$ or $f_1 \restriction \delta =_{J_\delta} f_2
\restriction \delta$ 
\endroster
\endproclaim
\bigskip

\remark{Remark}  1) We can justify (f)$'$ by pcf theory quotation,
see below.
\nl
2) We should prove that the p.o.w.i.s. being existential holds.

Note that in proving \scite{p.3}, \scite{p.3A} the main point is the
``existential limit". This proof has affinity to the first step in the
elimination of quantifiers in the theory of $(\omega,<)$.  For this it
is better if $I_\theta = ({\Cal F},<_{J^{\text{bd}}_\theta})$ has many
cases of existence.  Toward this we ``padded it" in $(*)_0$ of the
proof - take care of successor $(f \in {\Cal F} \Rightarrow f+1 \in
{\Cal F})$, have zero $(0_\theta \in {\Cal F})$ without losing the
properties we have.
\nl
2) The demand of \scite{p.3} may seem very strong, but by pcf theory
it is $q$ natural.
\endremark
\bigskip

\demo{\stag{p.3.1} Observation}  1) Theorem \scite{p.3A} implies
Theorem \scite{p.3}.
\nl
2) If (a)-(d) of \scite{p.3} holds, \ub{then} $(f) \Rightarrow (f)'$.
\nl
3) If (a)-(d) of \scite{p.3} holds then (f) $\Rightarrow$ (e).
\enddemo
\bigskip

\demo{Proof}  1) By 2).
\nl
2) Let 

$$
\align
S =: \{\delta < \theta:&\delta \text{  is a limit ordinal and if }
f_1,f_2 \in {\Cal F}_\delta \\
  &\text{ then for some } \alpha < \delta \text{ we have }
 \beta \in (\alpha,\delta) \Rightarrow \\
  &(f_1(\beta) < f_2(\beta)
 \equiv f_1(\alpha) < f_2(\alpha))\}.
\endalign
$$
\mn
By (f) we know that $S$ is a stationary subset of $\theta$.  Let
${\Cal D}_\theta$ be the club filter on $\theta$ and $D =: {\Cal
D}_\theta + S$, it is a normal filter on $\theta$ and $S \in D$.  So
sub-clauses $(\alpha),(\beta),(\gamma)$ of $(f)'$ holds. 

Let $J_\delta = J^{\text{bd}}_\delta$ for $\delta \in S$ so $\bar J =
\langle J_\delta:\delta \in S \rangle$ satisfies sub-clauses
$(\delta),(\varepsilon)$ of $(f)'$.  To prove $(\zeta)$ assume $S'
\subseteq S,S' \in D^+$ and $w_\delta \in J_\delta$ for $\delta \in
S'$.  Then sup$(w_\delta) < \delta$ and $S'$ is a stationary subset of
$\delta$ hence by Fodor lemma for some $\beta(*) < \theta$ the set $S'' =
\{\delta \in S':\sup(w_\delta) = \beta(*)\}$ is a stationary subset of
$\theta$ and so $[\beta(*),\theta)$ is an end segment of $\theta$ and is
equal to $\cup\{[\beta(*),\delta):\delta \in S''\}$ which is included
in $\cup\{\delta \backslash w_\delta:\delta \in S'\}$, as required in
$(\zeta)$ from $(f)'$, so sub-clause $(\zeta)$ really holds.

To prove sub-clause $(\eta)$ of clause $(f)'$ note that what it says is
what is said in $(f)$. 
\nl
3) Should be clear.  Given $f_1,f_2 \in {\Cal F}$; by sub-clause
$(\eta)$ of $(f)'$ for each $\delta \in S$ there are $w_\delta \in
J_\delta$ and $\ell_\alpha < 3$ such that $\ell_0 = 0 \wedge \alpha
\in \delta \backslash w_\delta \Rightarrow f_1(\alpha) < f_2(\alpha)$
and $\ell_\delta = 1 \wedge \alpha \in \delta \backslash w_\delta
\Rightarrow f_1(\alpha) = f_2(\alpha)$ and $\ell_\delta = 2 \wedge
\alpha \in \delta \backslash w_\delta \Rightarrow f_1(\alpha) >
f_2(\alpha)$.  So for some $\ell < 2$ the set $S' := \{\delta \in
S:\ell_\delta = \ell\}$ is stationary, hence $\cup\{\delta \backslash
w_\delta:\delta \in S'\}$ include an end segment of $\theta$ and we
are easily done.  \hfill$\square_{\scite{p.3.1}}$
\enddemo
\bigskip

\demo{Proof of \scite{p.3A}}  Without loss of generality 
\mr
\item "{$(*)_0$}"  $(a) \quad (\forall f \in {\Cal F})
(\exists^\infty g \in {\Cal F})\bigl(f \restriction [1,\theta) = g
\restriction [1,\theta)\bigr)$; 
\nl

\hskip25pt moreover for $f \in {\Cal F}$ we have
\nl

\hskip25pt $\omega = \{g(0):g \in {\Cal F}$ and $g \restriction [1,\theta) -
f \restriction [1,\theta)\}$
\sn
\item "{${{}}$}"  $(b) \quad \alpha < \beta < \theta \Rightarrow 
{\Cal F}_\alpha = \{f \restriction \alpha:f \in {\Cal F}_\beta\}$;
moreover $\alpha < \theta \Rightarrow {\Cal F}_\alpha =$
\nl

\hskip25pt  $\{f \restriction \alpha:f \in {\Cal F}\}$
\sn
\item "{${{}}$}"  $(c) \quad$ if $f \in {\Cal F}$, \ub{then} $f+1 \in {\Cal F}$
\sn
\item "{${{}}$}"  $(d) \quad$ the 
$f \in{}^\theta\{0\}$, the constantly zero function, belongs to ${\Cal F}$.
\ermn
[Why?  Let ${\Cal F}' = \{f \in {}^\theta\kappa$: for some $n,(\forall
\alpha < \theta)(f(1 + \alpha)=n) \wedge f(0) < \omega$ 
or for some $f' \in{\Cal F}$ and $n
< \omega$ we have $(\forall \alpha < \theta)(f(1 + \alpha)= \omega
(1+f'(\alpha))+n) \wedge f(0) < \omega\}$ 
and for $\alpha < \theta$, replace ${\Cal F}_\alpha$
by ${\Cal F}'_\alpha =\{f \restriction \alpha:f \in {\Cal F}'\}$.  Now
check that $(a)-(e),(f)'$ of the assumption still holds.]

We define ${\frak s} = (J,\bar I,\bar \pi)$ as follows:
\mr
\item "{$(*)_1$}"  $(a) \quad J = (\theta +1;<)$
\sn
\item "{${{}}$}"  $(b)(\alpha) \quad$ let $I_\theta =
({\Cal F},<_{J^{\text{bd}}_\theta})$ and
\sn
\item "{${{}}$}"  $\quad (\beta) \quad I_\alpha = 
({\Cal F}_{1 +\alpha +1},<_{\alpha +1})$ for $\alpha < \theta$ where
$$
f_1 <_{\alpha +1} f_2 \Leftrightarrow f_1(1 + \alpha) < f_2(1 +\alpha)
$$
\sn
\item "{${{}}$}"  $\quad (c) \quad$ for $\alpha < \beta < \theta +1$ let
$\pi_{\alpha,\beta}:I_\beta \rightarrow I_\alpha$ be
$$
\pi_{\alpha,\beta}(f) = f \restriction (1 +\alpha +1).
$$
\ermn
Note that
\mr
\item "{$(*)_2$}"     $I_\alpha$ is a non-trivial (see Definition
\scite{m.1}(6)).
\ermn
[Why?  By $(*)_0(a)$ and the choice of $<_{I_\alpha}$ in $(*)_1(b)(\beta)$.] 
\mr
\item "{$(*)_3$}"  ${\frak s} = (J,I,\bar \pi)$ is a p.o.w.i.s. even
nice
\nl
[Note clause (d) of Definition \scite{p.1} holds by clause (e) of
Theorem \scite{p.3}.]
\sn
\item "{$(*)_4$}"  ${\frak s}$ is a limit of ${\frak t} =: {\frak
s} \restriction \theta = ((\theta,<),\bar I \restriction
\theta,\bar \pi \restriction \theta)$.
\nl
[Why?  Note that clause (d) of Definition \scite{p.1} holds by clause
(f) here and Fodor lemma.  Easy to check the other clauses.]
\sn
\item "{$(*)_5$}"  ${\frak t}$ is a $\kappa$-p.o.w.i.s.
\nl
[Why?  Check, as $\alpha < \theta \Rightarrow |{\Cal F}_\alpha| \le
\kappa$.]
\ermn
Now $G^{\frak s}_\theta$ is an almost $\kappa$-automorphism group by
Claim \scite{p.2}, the ``existential limit" holds by $(*)_6$
below (note: $J$ is linear).  Now rk$(I^{\frak s}_\theta)
= \gamma$ and $H^{\frak s}_\theta$ is a subgroup of $G^{\frak
s}_\theta$ of cardinality $2 \le \kappa$.  

By \scite{m.9}

$$
\tau^{\text{nlg}}_{G^{\frak s}_\theta,G^{<1,{\frak s}}_\theta} =
\text{ rk}(I^{\frak s}_\theta) = \gamma
$$
\mn
and nor$^{< \infty}_{G^{\frak s}_\theta}
(H^{\frak s}_\theta) = G^{\frak s}_\theta$ and by \scite{nx.3}(4),
$\tau^{\text{nlf}}_{G^{\frak s}_\theta,H^{\frak s}_\theta} = \gamma$.

We still have to check
\mr
\item "{$(*)_6$}"  ``${\frak s}$ is an existential limit of
${\frak t}$", see Definition \scite{p.1.7}.
\ermn
That is we have to prove clause (f) of \scite{p.1.7}, so we should
prove its conclusion, assuming its assumption which means in our case
\mr
\item "{$\circledast_1$}"  $(a) \quad k=k_1 + k_2,{\Cal E}$ is an equivalence
relation on ${\Cal S}^k$
\sn
\item "{${{}}$}"  $(b) \quad \bar f \in {}^{k_1}({\Cal F}_\theta)$ and
$\alpha(*) < \delta$
\sn
\item "{${{}}$}"  $(c) \quad \bar e = \langle e_\alpha \in
[\alpha(*),\theta)\rangle$ is such that $e_\alpha \in {\Cal S}^k/{\Cal E}$
\sn
\item "{${{}}$}"  $(d) \quad \langle \bar g^\alpha:\alpha \in
[\alpha(*),\theta)\rangle$ is such that $\bar g^\alpha \in
{}^{(k_2)}({\Cal F}_\alpha)$
\sn
\item "{${{}}$}"  $(e) \quad$ if $\alpha(*) \le \alpha < \beta$ then:
$$
e_\alpha \text{ is the } 
{\Cal E}\text{-equivalence class of tp}_{\text{qf}}
(\langle f_\ell(1+\alpha):\ell <
k_1\rangle \char 94 \langle g^\beta_\ell(1+\alpha):\ell <
k_2\rangle,\emptyset,\kappa).
$$
\ermn
Without loss of generality [recalling clause (e) of the assumption and 
$(*)_0(c)$]
\mr
\item "{$\circledast_2$}"  $(f) \quad \langle f_\ell:\ell < k_1\rangle$
is $\le_{J^{\text{bd}}_\theta}$-increasing
\sn
\item "{${{}}$}"  $(g) \quad f_0$ is constantly zero
\sn
\item "{${{}}$}"  $(h) \quad$ for each $\ell < k_1-1$ we have: 
$f_{\ell +1} = f_\ell$ mod
$J^{\text{bd}}_\theta$ or $f_{\ell +1} = f_\ell +1$ 
\nl

\hskip25pt mod $J^{\text{bd}}_\theta$ or 
$f_\ell + \omega \le f_{\ell +1}$ mod $J^{\text{bd}}_\omega$
\sn
\item "{${{}}$}"  $(i) \quad \langle f_\ell:\ell < k_1\rangle$ is
without repetition
\sn
\item "{${{}}$}"  $(j) \quad \langle f_\ell(0):\ell < k_1\rangle$ is
without repetition.
\ermn
Possibly increasing $\alpha(*) < \theta$ \wilog \,
\mr
\item "{$\circledast_3$}"  if $\alpha \in [\alpha(*),\theta)$ and
$\ell_1,\ell_2 < k_1$ then $f_{\ell_1}(\alpha) < f_{\ell_2}(\alpha)
\Leftrightarrow f_{\ell_1}(\alpha(*)) < f_{\ell_2}(\alpha(*))$.
\ermn
Hence by clause (f) of $\circledast_2$
\mr
\item "{$\circledast_4$}"   $\langle f_\ell(\alpha(*)):\ell <
k_1\rangle$ is non-decreasing.
\ermn
For notational simplicity
\mr
\item "{$\circledast_5$}"  $(a) \quad \langle f_e \restriction
\delta:\ell < k_1\rangle = \langle g^\delta_\ell:\ell < k_1\rangle$ so
$k_1 < k_2$
\sn
\item "{${{}}$}"  $(b) \quad$ if $\ell_1 < k_2,\ell_2 \in [k_1,k_2)$
then $g^\delta_{\ell_1} = g^\delta_{\ell_2} \equiv
g^\delta_{\ell_1}(0)= g^\delta_{\ell_2}(0)$.
\ermn
Next for some $p^*$
\mr
\item "{$\circledast_6$}"  $p^* \in {\Cal S}^k$ and for some $S'
\subseteq S$ from $D^+$, for every $\delta \in S'$ for the 
$J_\delta$-majority of $\alpha < \delta$, say $\alpha \in \delta
\backslash w_\alpha,w_\alpha \in J_\delta$, 
 we have $p^* = \text{ tp}_{\text{qf}}
(\langle g^\delta_\ell \restriction (1 + \alpha +1):\ell <
k_2\rangle,\emptyset,I_\alpha)$.
\ermn
[Why?  By sub-clause $(\eta)$ of clause $(f)'$, as $J_\delta$ is an
ideal (applied to $(g^\delta_{\ell_1},g^\delta_{\ell_2})$ for every
$\ell_1,\ell_2 < k_2$) for 
each $\delta \in S$ we can choose $w_\delta \in J_\delta$
and $q_\delta \in {\Cal S}^k$ such that for every $\alpha \in (\delta
\backslash w)$ we have tp$_{\text{qf}}
(\langle g^\delta_\ell(1 + \alpha):\ell < k_2\rangle,
\emptyset,I_\alpha)$ is equal to
$q_\delta$.  For each $p \in {\Cal S}^k$ let
$S_p = \{\delta \in S:q_\delta = p\}$.  So $S = \cup\{S_p:p \in {\Cal S}^k\}$,
hence for some $p$ we have $S_p \in D^+$. So let $S' = S_p,p^* = p$.]

So \wilog \, considering the way $I_\alpha$ was defined by $\circledast_5$
\mr
\item "{$\circledast_7$}"  there are $E^*_1,E^*_2,<_*$ such that
{\roster
\itemitem{ $(a)$ }  $E^*_1$ is an equivalence relation on $k_2 =
\{0,\dotsc,k_2-1\}$
\sn
\itemitem{ $(b)$ }  $E^*_2$ is an equivalence relation on $k_2$
refining $E^*_1$
\sn
\itemitem{ $(c)$ }  $<_*$ linearly order $k_2$
\sn
\itemitem{ $(d)$ }  if $\delta \in S',\alpha \in \delta \backslash
w_\delta$ so $p^* = \text{ tp}_{\text{qf}}(\langle
g^\delta_\ell(\alpha):\ell < k_2\rangle$ then:
\sn
\itemitem{ ${{}}$ }  $(\alpha) \quad \ell_1 E^*_2 \ell_2$ iff
$g^\delta_{\ell_1}(1 + \alpha) = g^\delta_{\ell_2}(1 + \alpha)$
\sn
\itemitem{ ${{}}$ }   $(\beta) \quad \ell_1 E^*_2 \ell_2$ iff
$g^\delta_{\ell_1} \restriction (1 + \alpha +1) = 
g^\delta_{\ell_2} \restriction (1 + \alpha +1)$
\sn
\itemitem{ ${{}}$ }   $(\gamma) \quad (\ell_1/E^*_1) <_* (\ell_2/E^*_1)$ iff
$g^\delta_{\ell_1}(1 + \alpha) < g^\delta_{\ell_2}(1 + \alpha)$.
\endroster}
\ermn
Let $\langle u_0,\dotsc,u_{m-1}\rangle$ list the $E^*_1$-equivalence
classes in $<_*$-increasing order.  Necessary $0 \in u_0$.

Let $\alpha_* = \text{ min}(\delta_* \backslash w_{\delta_*})$ where
$\delta_* = \text{ min}(S')$.  We now define $g_\ell \in {}^\theta
\kappa$ for $\ell < k_2$ as follows.  So necessarily 
for a unique $i=i(\ell),\ell \in u_i$ and let $i_1 =
i_1(\ell) \le i$ be maximal such that $u_{i_1} \cap
\{0,\dotsc,k_1-1\} \ne \emptyset,j_2 = j_2(\ell) = \text{ min}(\{u_1
\cap \{0,\dotsc,k_1-1\})$.  It is well defined as necessary $0 \in
u_0$ because $f_0$ is constantly zero.  Now we let
\mr
\item "{$\boxdot_0$}"  $g_\ell = (g^{\alpha_*}_\ell \restriction
\{0\}) \cup ((f_{j_2} + (i-i_1)) \restriction [1,\theta))$.
\ermn
Now
\mr
\item "{$\boxdot_1$}"  if $\ell < k_1$ then $g_\ell = f_\ell$
\nl
[Why?  Check the definition $g^{\alpha_*}_\ell(0) = f_\ell(0)$ as
$g^{\alpha_*}_\ell = f_\ell$.]
\sn
\item "{$\boxdot_2$}"   $g_\ell \in {\Cal F}$ for $\ell < k_2$
\nl
[Why?  As $f_{j_2} \in {\Cal F}$ and clauses (a)+(c) of $(*)_0$.]
\sn
\item "{$\boxdot_3$}"  if $\ell_1 E^*_2 \ell_2$ then $g_{\ell_1} =
g_{i_2}$
\nl
[Why?  First, as $\ell_1 E^*_2 \ell_2$ we have $g_\ell(0) =
g^{\alpha_*}_{\ell_1}(0) = g^{\alpha_*}_{\ell_2}(0) = g_\ell(0)$.
Second, clearly $i(\ell_1) = i(\ell_2),i_1(\ell_1) = i_1(\ell_2)$ and
$j_2(\ell_1) = j_2(\ell_2)$ hence for $\alpha \in [1,\theta)$ we have
$$
\align
g_{\ell_1}(\alpha) = (&f_{j_2(\ell_1)}(\alpha) +
(i(\ell_1)- i_1(\ell_1)) = \\
  &f_{j_2(\ell_1)}(\alpha) + (i(\ell_2) - i_1(\ell_2)) =
g_{\ell_2}(\alpha).
\endalign
$$
\mn
So we are done.]
\sn
\item "{$\boxdot_4$}"  if $\ell_1,\ell_2 < k_2$ but $\neg(\ell_1 E^*_2
\ell_2)$ then $g_{\ell_1} \ne g_{\ell_2}$
\nl
[Why?  If $\ell_1,\ell_2 < k_1$ then $g_{\ell_1} = f_{\ell_1} \ne
f_{\ell_2} = g_{\ell_2}$.  If $\ell_1 < k_2,\ell_2 \in [k_1,k_2)$ as
$\neg(\ell_1 E^*_1 \ell_2)$ by $(*)_5(b)$ we have
$g^{\alpha_*}_{\ell_1}(0) \ne g^{\alpha_*}_{\ell_2}(0)$, hence
$g_{\ell_1}(0) = g^{\alpha_*}_{\ell_1}(0) \ne g^{\alpha_*}_{\ell_2}(0)
= g_\ell(0)$ hence $g_{\ell_1} \ne g_{\ell_2}$.  Lastly, if $\ell_1
\in [k_1,k_2),\ell_2 < k_2$ the proof is similar.]
\sn
\item "{$\boxdot_5$}"  if $\ell_1,\ell_2 < k_2,\ell_1 E^*_1 \ell_2$
then $\neg(g_{\ell_1} <_{I_\theta} g_{\ell_2})$
\nl
[Why?  As $g_{\ell_1} \restriction [1,\theta) = g_{\ell_2}
\restriction [1,\theta)$, so $g_{\ell_1} = g_{\ell_2}$ mod
$J^{\text{bd}}_\theta$, so $I_\theta \models \neg(g_{\ell_1} <
g_{\ell_2})$.]
\sn
\item "{$\boxdot_6$}"  if $\ell_1,\ell_2 < k_2$ and $(\ell_1/E^*_1)
<_* (\ell_2/E^*_2)$ then $g_{\ell_1} <_{I_\theta} g_{\ell_2}$
\nl
[Why?  If $f_{j_2(\ell_1)} + \omega \le f_{j_2(\ell_2)}$ mod
$J^{\text{bd}}_\theta$ then easily $g_{\ell_1}
<_{J^{\text{bd}}_\theta} f_{j_2(\ell_1)} + w
\le_{J^{\text{bd}}_\theta} f_{j_2(\ell)} \le_{J^{\text{bd}}_\theta}
g_{\ell_2}$ so we are done.
If $j_2(\ell_1) = j_2(\ell_2)$ then as still $i(\ell_1) < i(\ell_2)$
we have $g_{\ell_1} =_{J^{\text{bd}}_\theta} f_{j_2(\ell_1)} +
(i(\ell_1)-j_2(\ell_1) < f_{j_2(\ell_1)} + (i(\ell_2) = g_2(\ell_2))
=_{J^{\text{bd}}_\theta} g_{\ell_2}$ as required.
If $j_2(\ell_1) \ne j_2(\ell_2)$ then necessarily $j_2(\ell_1) <
j_2(\ell_2),i_1(\ell_1) < i_1(\ell_2)$ moreover $i_1(\ell_1) \le
i(\ell_1) < j_2(\ell_2) \le i(\ell_2)$ but by $\circledast(h)$ we have
$f_{j_1(\ell_1)} + (j_2(\ell_1) - i_1(\ell_1))
\le_{J^{\text{bd}}_\theta} f_{j_2(\ell_2)}$ so we are easily done.]
\ermn
Together $\langle g_\ell:\ell < k_2\rangle$ is as required for proving
$(f)'$ of \scite{p.1.7}, the definition of existential limit,
i.e. $(*)_6$.     \hfill$\square_{\scite{p.3}} \quad
\square_{\scite{p.3A}}$
\enddemo
\bn
We quote
\proclaim{\stag{p.3C} Claim}  Assume {\rm cf}$(\kappa) = \theta >
\aleph_0,\alpha < \kappa \Rightarrow (\alpha)^\theta < \kappa$ and 
$\lambda = \kappa^\theta$.  \ub{Then} we can find $\langle {\Cal
F}_i:i \le \theta\rangle,S,D$ satisfying the conditions from
\scite{p.3A} with $\gamma = \lambda$ (and more).
\endproclaim
\bigskip

\demo{Proof}  By \scite{p.3B} and \cite{Sh:g}.
\hfill$\square_{\scite{p.3C}}$
\enddemo
\bigskip

\proclaim{\stag{p.3B} Claim}  Assume
\mr
\item "{$\circledast$}"  $(a) \quad \bar \lambda = \langle \lambda_i:i
< \theta \rangle$ is an increasing sequence of regular cardinals with
\nl

\hskip25pt limit $\kappa$
\sn
\item "{${{}}$}"  $(b) \quad \lambda = \text{\rm tcf}(\dsize \prod_{i
< \theta} \lambda_i,<_{J^{\text{bd}}_\theta})$
\sn
\item "{${{}}$}"  $(c) \quad$ {\rm max pcf}$\{\lambda_i:i < j_*\} < \kappa$
for every $j < \theta$.
\ermn
1) \ub{Then} there are $D,S^*,u$ such that
\mr
\item "{$(\alpha)$}"  $u \in [\theta]^\theta,S^* \subseteq \theta$ is
stationary
\sn
\item "{$(\beta)$}"  there are no $\zeta < \theta,u_\varepsilon \in
[u]^\theta$ for $\varepsilon < \theta$ such that for a club of $\delta
< \theta$ if $\delta \in S^*$ then for at least one $\varepsilon <
\delta$ we have {\rm max pcf}$\{\lambda_i:i \in\delta \cap u_\varepsilon\} 
< \text{\rm max pcf}\{\lambda_i:i \in \delta\}$ hence
\sn
\item "{$(\gamma)$}"  $D$ is a normal filter on $\theta$ where: $D$ is
$\{S \subseteq \theta$: for some sequence $\langle
u_\varepsilon:\varepsilon < \theta\rangle$ of subsets of $\theta$ each
of cardinality $\theta$ and for some club $E$ of $\theta$, if
$\delta \in E \cap S \cap S^*$ then for every $\varepsilon < \delta$
we have {\rm max pcf}$\{\lambda_i:i \in \delta \cap u_\varepsilon\} =
\text{\rm max pcf}\{\lambda_i:i \in \delta \cap u\}\}$
\sn
\item "{$(\delta)$}"  by renaming $u=\theta$ and for $\delta \in S^*$
let $J_\delta = \{u \subseteq \delta:\,\text{\rm max pcf}\{\lambda_i:i
\in \delta \backslash u\} < \text{\rm max pcf}\{\lambda_i:i <
\delta\}$.
\ermn
2) We can choose ${\Cal F}_i \subseteq \dsize \prod_{j<i} \lambda_i$
for $i \le \theta$ such that all the conditions in \scite{p.3A} holds.
\endproclaim
\bigskip

\demo{Proof}  By \cite[II,3.5]{Sh:g}, see on this \cite[\S18]{Sh:E12}.
\enddemo
\bigskip

\demo{\stag{p.5} Conclusion}  If $\kappa$ is strong limit singular of
uncountable cofinality \ub{then} $\tau^{\text{atw}}_\kappa \ge
\tau^{\text{nlg}}_\kappa \ge \tau^{\text{nlf}}_\kappa > 2^\kappa$.
\enddemo
\bigskip

\demo{Proof}  By \scite{p.3A} and \scite{p.3B}.  \hfill$\square_{\scite{p.5}}$
\enddemo
\bigskip

\remark{\stag{5.7} Remark}   1)  If $\kappa = \kappa^{\aleph_0}$ do we have
$\tau^{\text{atw}}_\kappa \ge \tau^{\text{nlg}}_\kappa \ge
\tau^{\text{nlf}}_\kappa > \kappa^+$?  But if $\kappa = \kappa^{<
\kappa} > \aleph_0$ then quite easily yes.
\nl
2) In \scite{p.5} we can weaken ``$\kappa$ is strong limit".  E.g. if
$\kappa$ has uncountable cofinality and $\alpha < \kappa \Rightarrow
|\alpha|^{\text{cf}(\kappa)} < \kappa$, then $\tau^{\text{nlf}}_\kappa
> \kappa^{\text{cf}(\kappa)}$; see more in \cite[\S18]{Sh:E12}.
\nl
3) We elsewhere will weaken the assumption in \scite{p.3},
\scite{p.3A} but deduce only that $\tau^{\text{nlg}}_\kappa$ is large.
\endremark
\newpage


\nocite{ignore-this-bibtex-warning} 
\newpage
    
REFERENCES.  
\bibliographystyle{lit-plain}
\bibliography{lista,listb,listx,listf,liste}

\def\germ{\frak} \def\scr{\cal} \ifx\documentclass\undefinedcs
  \def\bf{\fam\bffam\tenbf}\def\rm{\fam0\tenrm}\fi 
  \def\defaultdefine#1#2{\expandafter\ifx\csname#1\endcsname\relax
  \expandafter\def\csname#1\endcsname{#2}\fi} \defaultdefine{Bbb}{\bf}
  \defaultdefine{frak}{\bf} \defaultdefine{=}{\B} 
  \defaultdefine{mathfrak}{\frak} \defaultdefine{mathbb}{\bf}
  \defaultdefine{mathcal}{\cal}
  \defaultdefine{beth}{BETH}\defaultdefine{cal}{\bf} \def\bbfI{{\Bbb I}}
  \def\mbox{\hbox} \def\text{\hbox} \def\om{\omega} \def\Cal#1{{\bf #1}}
  \def\pcf{pcf} \defaultdefine{cf}{cf} \defaultdefine{reals}{{\Bbb R}}
  \defaultdefine{real}{{\Bbb R}} \def\restriction{{|}} \def\club{CLUB}
  \def\w{\omega} \def\exist{\exists} \def\se{{\germ se}} \def\bb{{\bf b}}
  \def\equivalence{\equiv} \let\lt< \let\gt>
  \def\implies{\Rightarrow}\def\mathfrak{\bf}\def\germ{\frak} \def\scr{\cal}
  \ifx\documentclass\undefinedcs
  \def\bf{\fam\bffam\tenbf}\def\rm{\fam0\tenrm}\fi 
  \def\defaultdefine#1#2{\expandafter\ifx\csname#1\endcsname\relax
  \expandafter\def\csname#1\endcsname{#2}\fi} \defaultdefine{Bbb}{\bf}
  \defaultdefine{frak}{\bf} \defaultdefine{=}{\B} 
  \defaultdefine{mathfrak}{\frak} \defaultdefine{mathbb}{\bf}
  \defaultdefine{mathcal}{\cal}
  \defaultdefine{beth}{BETH}\defaultdefine{cal}{\bf} \def\bbfI{{\Bbb I}}
  \def\mbox{\hbox} \def\text{\hbox} \def\om{\omega} \def\Cal#1{{\bf #1}}
  \def\pcf{pcf} \defaultdefine{cf}{cf} \defaultdefine{reals}{{\Bbb R}}
  \defaultdefine{real}{{\Bbb R}} \def\restriction{{|}} \def\club{CLUB}
  \def\w{\omega} \def\exist{\exists} \def\se{{\germ se}} \def\bb{{\bf b}}
  \def\equivalence{\equiv} \let\lt< \let\gt>
\begin{thebibliography}{JShT 654}
\makeatletter \renewcommand{\@biblabel}[1]{[#1]} \makeatother
\def\eprintfootnotetext{References of the form {\tt math.XX/$\cdots$}
 refer to the {\tt xxx.lanl.gov} archive}
\ifx\documentstyle\undefinedcontrolsequence
   \def\anyfootnote{\footnote{*}}
   \else\def\anyfootnote{\footnote}\fi
\def\eprintfn{\ifEprint\anyfootnote{\eprintfootnotetext}\fi\Eprintfalse }
\newif\ifEprint  \Eprinttrue

\bibitem[JShT 654]{JShT:654}Winfried Just, Saharon Shelah, and Simon Thomas.
\newblock {The automorphism tower problem revisited}.
\newblock {\em {Advances in Mathematics}}, {\bf 148}:243--265, 1999.
\newblock math.LO/0003120.

\bibitem[Sh:E12]{Sh:E12}Saharon Shelah.
\newblock {Analytical Guide and Corrections to \cite{Sh:g}.}
\newblock math.LO/9906022.

\bibitem[Sh:c]{Sh:c}Saharon Shelah.
\newblock {\em {Classification theory and the number of nonisomorphic models}},
  volume~92 of {\em {Studies in Logic and the Foundations of Mathematics}}.
\newblock {North-Holland Publishing Co., Amsterdam, xxxiv+705 pp}, 1990.

\bibitem[Sh:g]{Sh:g}Saharon Shelah.
\newblock {\em {Cardinal Arithmetic}}, volume~29 of {\em {Oxford Logic
  Guides}}.
\newblock {Oxford University Press}, 1994.

\bibitem[Sh:F579]{Sh:F579}{Shelah, Saharon}.
\newblock {More on tower of automorphism of groups}.

\bibitem[Th85]{Th85}Simon Thomas.
\newblock {The automorphism tower problem}.
\newblock {\em Proceedings of the American Mathematical Society}, {\bf
  95}:166--168, 1985.

\bibitem[Th98]{Th98}Simon Thomas.
\newblock {The automorphism tower problem II}.
\newblock {\em Israel Journal of Mathematics}, {\bf 103}:93--109, 1998.

\bibitem[Wel39]{Wel39}H.~Wielandt.
\newblock {Eine Verallgemeinerung der invarianten Untergruppen}.
\newblock {\em Math. Z.}, {\bf 45}:209--244, 1939.

\end{thebibliography}

\enddocument